\documentclass{article}

\usepackage{fullpage}

\usepackage{amsmath,amsfonts,amsthm,amssymb}
\usepackage{color,xcolor}
\usepackage{graphicx,graphics}
\usepackage{subfigure}
\usepackage{url}
\usepackage{hyperref}

\usepackage{stmaryrd}

\usepackage{multirow}

\usepackage{comment} 
\usepackage{tikz}
\usetikzlibrary{mindmap}

\usepackage{lscape}
 \usepackage{threeparttable}

\usepackage{algorithmic,algorithm}
\usepackage{todonotes}

\definecolor{mygreen}{RGB}{28,172,0} 
\definecolor{mylilas}{RGB}{170,55,241}
\definecolor{gray}{rgb}{0.5, 0.5, 0.5}
\definecolor{darktaupe}{rgb}{0.28, 0.24, 0.2}
\definecolor{arsenic}{rgb}{0.23, 0.27, 0.29}
\definecolor{azure(colorwheel)}{rgb}{0.0, 0.5, 1.0}
\definecolor{bostonuniversityred}{rgb}{0.8, 0.0, 0.0}
\definecolor{darkWhite}{rgb}{0.94,0.94,0.94}
\definecolor{arsenic}{rgb}{0.23, 0.27, 0.29}
\definecolor{azure(colorwheel)}{rgb}{0.0, 0.5, 1.0}
\definecolor{bostonuniversityred}{rgb}{0.8, 0.0, 0.0}

\newtheorem{theorem}{Theorem}[section]

\usepackage{hyperref}
\usepackage{lscape}
 \usepackage{threeparttable}

\usepackage[capitalise,noabbrev]{cleveref}

\hypersetup{
    colorlinks=true,
    linkcolor=blue,
    linkbordercolor = {white},
    urlcolor=arsenic,
    linktoc=all
}

\newcommand{\refp}[1]{(\ref{#1})}

\newcommand{\revMR}[1]{{#1}}
\newcommand{\revDT}[1]{{#1}}
\definecolor{mygreen}{RGB}{114,134,57}
\newcommand{\revsix}[1]{{#1}}

\def\restriction#1#2{\mathchoice
              {\setbox1\hbox{${\displaystyle #1}_{\scriptstyle #2}$}
              \restrictionaux{#1}{#2}}
              {\setbox1\hbox{${\textstyle #1}_{\scriptstyle #2}$}
              \restrictionaux{#1}{#2}}
              {\setbox1\hbox{${\scriptstyle #1}_{\scriptscriptstyle #2}$}
              \restrictionaux{#1}{#2}}
              {\setbox1\hbox{${\scriptscriptstyle #1}_{\scriptscriptstyle #2}$}
              \restrictionaux{#1}{#2}}}
\def\restrictionaux#1#2{{#1\,\smash{\vrule height .8\ht1 depth .85\dp1}}_{\,#2}} 

\def\R{\mathbb{R}}
\def\P{\mathbb{P}}

\def\N{\mathbb{N}}
\def\L{\mathcal{L}}
\def\M{\mathcal{M}}
\def\mass{\mathbb{M}}
\def\diag{\mathbb{D}}
\def\rhs{\mathtt{r}}
\def\UU{\underline{U}}
\def\CFL{\text{CFL}}

\newcommand{\bx}{\mathbf{x}}

\newcommand{\footremember}[2]{%
   \footnote{#2}
    \newcounter{#1}
    \setcounter{#1}{\value{footnote}}%
}
\newcommand{\footrecall}[1]{%
    \footnotemark[\value{#1}]%
}

\author{Sixtine Michel\footremember{inria}{Team CARDAMOM, Inria Bordeaux sud-ouest, - 200 av.  de la vieille tour, 33405 Talence, France}, Davide Torlo\footrecall{inria}, Mario Ricchiuto\footrecall{inria}, R\'emi Abgrall\footremember{zurich}{Institut f\"ur Mathematik, Winterthurstrasse 190, CH 8057 Z\"urich, Switzerland.}   }

\title{Spectral analysis of continuous FEM for hyperbolic PDEs:
influence of  approximation, stabilization, and time-stepping}
\date{\today}

\begin{document}
\graphicspath{{images/}}

\maketitle
\abstract{
We study 
 continuous \revMR{finite element dicretizations} 
  for one dimensional hyperbolic partial differential equations. \revMR{The main contribution of the paper is to provide a fully discrete spectral 
  analysis, which is used to suggest optimal values of the CFL number and of the stabilization parameters involved in different types of 
  stabilization operators. In  particular, we analyze the}
	%
	streamline-upwind Petrov-Galerkin (SUPG) stabilization technique, the continuous interior penalty (CIP) stabilization method and the local projection stabilization (LPS). \revMR{Three different  choices for the continuous}
	finite element space  \revMR{are compared:}
	Bernstein polynomials, Lagrangian polynomials on equispaced nodes, and Lagrangian polynomials on \revMR{Gauss-Lobatto} cubature nodes. 
	\revMR{For the last choice, we only consider inexact quadrature based on the formulas corresponding to the degrees of freedom of the element,
	which allows to obtain a fully diagonal mass matrix.}
	%
\revMR{We also} 
	compare different \revMR{time stepping strategies, namely}  Runge-Kutta (RK), strong stability preserving RK (SSPRK) and deferred correction time integration methods. \revMR{The latter  allows to alleviate the computational cost as the mass matrix inversion is replaced by the high order correction iterations.}
	
\revMR{To understand the effects of these choices, both time-continuous and fully discrete Fourier analysis are performed. These allow to compare  all the different combinations in terms of accuracy and stability, as well as to provide suggestions for optimal values discretization parameters involved. The results are thoroughly verified   numerically both on linear and non-linear problems, and  error-CPU time curves are provided. Our final conclusions suggest that cubature elements combined with SSPRK and CIP or LPS stabilization are the most promising combinations.}
}


\section{Introduction}\label{sec:intro}
In this work we compare different numerical methods that can approximate the solution of the one dimensional hyperbolic conservation laws
\begin{equation}\label{eq:conservation_law}
	\partial_t u (x,t) + \partial_x f(u(x,t)) = 0 \quad x\in \Omega \subset \R, \, t\in \R^+,
\end{equation}
where $\Omega\subset \R$ is an interval, $f:\R^D\to\R^D$ is the flux function and $u:\Omega\to \R^D$ is the unknown of the system of equations. 
For the spectral analysis \revMR{of the numerical methods we will mainly focus on}  
the particular case of a linear flux
\begin{equation}\label{eq:conservation_law1}
 f(u(x,t)) = a u(x,t)\,,\;\;a=\text{const}\,.
\end{equation}

In this work, we compare different explicit high order accurate schemes based on the continuous Galerkin (CG) approach.
\revMR{In general,} the standard Finite Element Method (FEM) derived by this approach require the inversion of a large sparse mass matrix.
This procedure can be expensive as the matrix multiplication must be iterated for all the time steps.
Various techniques have been introduced to overcome the mass matrix inversion while keeping the high order accuracy of the scheme. 

The first strategy we study is the one proposed in \cite{DeC_2017}. There, to avoid the full mass matrix, a mass lumping is introduced, transforming the mass matrix into a diagonal one. The deferred correction (DeC) iterative time integration method alters the right--hand side in order to recover the original order of accuracy.
Another approach consists of a careful choice of quadrature points and basis functions in order to automatically obtain a diagonal mass matrix. We denote such elements as \textit{cubature} elements \cite{LIU2017cubature9order}. The classical use of Runge--Kutta methods will provide the high order accuracy also for the time discretization.

The second aspect we will focus on is the stabilization technique.
We emphasize that without any special treatment on the boundaries, such as the ones in \cite{abgrall2020analysis, abgrall2021analysis}, the CG methods are not stable for hyperbolic problems and there is the need of stabilization. In particular, this is always true when using periodic boundary conditions (BC). The CG discretizations with stabilization techniques can have dissipation levels that are comparable to the ones brought by discontinuous Galerkin (DG) with upwind numerical flux of the same order of accuracy, still remaining decently stable \cite{moura2020EigenanalysisGJP,moura2020SpatialEigenanalysis}.
The stabilization terms play an important role and we will compare three of them.
The first is the streamline upwind Petrov--Galerkin (SUPG) stabilization \cite{article_supg1, article_supg2}, which is strongly consistent, but it is also introducing new terms in the mass matrix \revMR{which are necessary to retain the appropriate consistency order. This can only be 
alleviated when using DeC time stepping. }
The second \revMR{approach} is the so--called continuous interior penalty (CIP) method \cite{article_cip5, article_cip3, article_cip4}, which penalizes the jump of the derivative of the solution across cell boundaries. This stabilization does not affect the mass matrix and, therefore, can be easily combined with mass--matrix free methods.
The last is the local projection stabilization \cite{AMTX-LPS_2010}, which penalizes the $\mathbb{L}^2$ projection of the gradient of the error \revMR{within the elements}. 
This technique does not affect the mass matrix, but it requires the solution of another linear system \revMR{for the $\mathbb{L}^2$ projection.
In this respect, the choice of the finite element space  and of the quadrature  have enormous impact on the  cost of the method.}

The goal of this work is to analyze the different methods and their combinations, \revMR{and give suggestions concerning
the most convenient choices in terms of accuracy, stability, and cost. To achieve this objective an important role is played 
by a spectral analysis which we perform both in the time-continuous and fully discrete cases. The analysis} 
reveals the best parameters (stabilization and CFL coefficients) that can be stably used in  \revMR{practice}. 

Numerical simulations 
\revMR{for both linear and non-linear} scalar \revMR{problems, and for the shallow water system} 
confirm 
the theoretical results, \revMR{and allow to further investigate the impact of the discretization choices  on the performance of the schemes
and on their cost.} 

The paper is organized as follows.
In Section \ref{sec:discretization} we introduce the different discretization methods, starting from the choice of the elements, then discussing the stabilization terms and finally presenting the different time integration methods. \cref{sec:fourier,sec:spectralAnalysis} \revMR{are dedicated to the} 
Fourier stability analysis. \revMR{  In \cref{sec:nonlinear} we provide some elements concerning the extension of the stabilization methods discussed to nonlinear problems, and finally in}
 \cref{sec:simulation} we show  
 numerical results on \revMR{linear and nonlinear problems. 
The paper is ended by a summary and overlook on future perspectives in} \cref{sec:conclusion}.

\section{Numerical Discretization}\label{sec:discretization}
We are interested  in  the approximation of \revMR{ solutions of}  \eqref{eq:conservation_law}   on a \revMR{ tessellation  of non overlapping celles, which we denote by} $\Omega_h$. 
\revMR{ We denote by $K$ the generic cell of $\Omega_h$, and more precicely $\Omega_h=\bigcup K$.
We also introduce the set of internal element boundaries (cell faces in 2D and 3D, cell nodes in 1D) of $\Omega_h$, which we denote by   $\mathcal{F}_h$.} \revDT{$h$ denotes the characteristic mesh size of $\Omega_h$.}   
The discrete solution is sought in a \revMR{ continuous finite element} space  $V_h^p = \lbrace v_h \in{\color{red} \mathcal{C}^0}(\Omega_h) : \quad \restriction{v_h}{K} \in \P_p(K) \quad \forall K\in \Omega_h \rbrace $. 
We are interested  in particular nodal finite elements, and  we will denote by $\varphi_j$ the basis functions associated to the degree of freedom $j$, so that  $V_h^p=\text{span}\left\{\varphi_j\right\}_{j\in\Omega_h}$ and we can write $u_h(x)=\sum_{j\in\Omega_h} u_j \varphi_j(x)$.   

The unstabilized   approximation of  \eqref{eq:conservation_law} reads: find $u_h\in V_h^p$ such that for any $v_h\in W_h\subset \mathbb{L}_2(\Omega_h)$
\begin{equation}\label{remi:1}
\revMR{ \int_{\Omega} v_h \partial_t u_h dx  - \int_{\Omega} \partial_x v_h f(u_h)\; dx + \left[ v_h f(u_h)\right]_{\partial \Omega} =0. }
\end{equation}

The main topic of this paper is the study of \revMR{ the linear  
stability of  \eqref{remi:1} and  of several stabilized variants using  Fourier's analysis. We will therefore assume periodic boundary conditions. 
We aim  at characterizing the schemes both in terms of their stability range and their accuracy in the fully discrete case, for different choices of the stabilization
strategy and of the time stepping\revDT{.}}
The extensions of these discretization techniques to more dimensions is well known in literature, even if sometimes not uniquely defined. 
We believe that the one dimensional study can provide useful information also in that context. 

\revMR{ As already said, we will consider several stabilized variants of \eqref{remi:1} which can be all written in the generic form:
find $u_h\in V_h^p$ that satisfies  
\begin{equation}\label{remi:2}
\int_{\Omega} v_h ( \partial_t u_h   + \partial_x f(u_h)) dx + S(v_h,u_h)=0, \quad \forall v_h \in V^p_h
\end{equation}
having re-integrated by parts and used the continuity of the approximation, and the periodicity of the boundary conditions to pass to the strong form of the PDE, and
with}  $S$ \revDT{being} a bilinear operator defined on $V^p_h\times V_h^p$. \revMR{ Several different choices for $S$ exist, and are discussed in detail in the following  sections. } 

\subsection{Stabilization Terms}\label{sec:stabilization}
\subsubsection{Streamline-Upwind/Petrov-Galerkin - SUPG} \label{SUPG}
This method \revMR{was   introduced in  \cite{article_supg01} (see also \cite{hst10, article_supg2}  and references therein) and is strongly consistent in the sense that it vanishes when replacing the discrete solution with the exact one.   
It can be written as a Petrov-Galerkin method replacing $v_h$ in   \eqref{remi:1} with a test function belonging to the space 
\begin{equation}
W_h := \{ w_h:\quad w_h=v_h+\tau_K \partial_u f(u_h)\partial_x v_h; \quad v_h \in V_h^p \} . \label{eq_supg0}
\end{equation}}
Here $\tau_K$ denotes a positive definite stabilization parameter with the dimensions of a  time-step  that we will assume to be  constant for every element.
\revMR{ Although other definitions are possible, here  we will evaluate this parameter as
$$
\tau_K = \delta \dfrac{h_K}{\|\partial_u f\|_K}
$$
where \revDT{$h_K$ is the cell diameter} and the  denominator  represents a reference value of the flux Jacobian norm on the element $K$.}

The final stabilized variational formulation reads
\begin{equation}
\int_{\Omega} v_h \partial_t u_h \; dx + \int_{\Omega} v_h \partial_x f(u_h)\; dx + \underbrace{\sum_{K \in \Omega} \int_{K} \big (\partial_u f(u_h) \partial_x v_h) \tau_K
	\left(\partial_t u_h + \partial_x f(u_h) \right)\; dx}_{S(v_h,u_h)} = 0.
\label{eq_supg1}
\end{equation}

\revMR{To characterize the accuracy of the  method, we can use the consistency analysis discussed e.g.  in \cite[\S3.1.1 and \S3.2]{AR:17}. In particular, of a finite element polyomial approximation of degree $p$ we can 
easily  show that given a smooth exact solution $u^e(t,x)$, replacing formally  $u_h$ by the projection of $u^e$ on the finite element space,
we can write 
\begin{equation}
\begin{split}
\epsilon(\psi_h) :=& \Big|
\int_{\Omega} \psi_h \partial_t (u_h^e - u^e) \; dx - \int_{\Omega} \partial_x \psi_h (\partial_x f(u_h^e)-\partial_x f(u^e))\; dx \\ + &
\sum_{K \in \Omega}\sum\limits_{l,m \in K} \dfrac{\psi_l - \psi_m}{k+1}  \int_{K} \big (\partial_u f(u_h) \partial_x \varphi_i) \tau_K
	\left(\partial_t (u_h^e - u^e) + \partial_x ( f(u_h^e) -f(u^e)) \right)\; dx 
	\Big| \le C h^{p+1},
\end{split}
\label{eq_supg2}
\end{equation}
\revDT{with $C$ \revsix{a constant} independent of $h$,} for all  functions $\psi$ of class at least $\mathcal{C}^1(\Omega)$, of which  $\psi_h$  denotes the finite element projection.
A key point in this estimate is the strong consistency of the method allowing to subtract its formal  application to the exact solution (thus subtracting zero),
and obtaining the above expression featuring differences between the exact solution/flux and its evaluation on the finite element space. 
Preserving this error estimate  precludes the possibility of lumping the mass matrix, and in particular  the entries associated to the stabilization term.
This makes the scheme relatively inefficient when using standard explicit time stepping.}

\revMR{As a final note, for a linear flux \eqref{eq:conservation_law1}, which is the  main focus of the analysis of this paper, and   for exact integration with  $\tau_K = \tau$,
a classical result is obtained in the time continuous case by testing  with $v_h =u_h + \tau\, \partial_t u_h$   to obtain \cite{article_supg2}
\begin{equation}
\begin{split}
\int\limits_{\Omega_h}\partial_t\left(\dfrac{u^2_h}{2}+\tau^2\dfrac{(a\partial_x u_h)^2}{2}\right) +
\int\limits_{\Omega_h}a\partial_x\left( \dfrac{u^2_h}{2}+\tau^2\dfrac{( \partial_t u_h)^2}{2}\right) = -\int\limits_{\Omega_h}\tau (\partial_tu_h+a\partial_xu_h)^2.
\end{split}
\label{eq_supg3}
\end{equation}
With  periodic boundary conditions this easily shows that the  norm $|||u|||^2 :=\int_{\Omega_h} \dfrac{u^2_h}{2}+\tau^2\dfrac{(a\partial_x u_h)^2}{2} dx$ is non-increasing.
The interested reader can refer to  \cite{article_supg2} for  the analysis of some  (implicit) fully discrete schemes.}



\subsubsection{Continuous Interior Penalty - CIP}  \label{CIP}
An alternative, which maintains the structure of the mass matrix, is the continuous interior penalty (CIP) stabilization used in \cite{article_cip5, article_cip3, article_cip4}. This method has been develop by E. Burman and P. Hansbo in \cite{article_cip6}, but it can be seen as a variation of the method originally proposed by Douglas and Dupont \cite{inbook_cip7}. 

This method stabilizes convection-diffusion-reaction problems by adding a least-squares term based on the jump in the gradient of the discrete solution over element boundaries. With this simple concept we obtain stability for convection-reaction-diffusion problems also in the vanishing viscosity limit.

The method reads
\begin{equation}
\revMR{\int_{\Omega_h} v_h \partial_t u_h \; dx + \int_{\Omega_h} v_h \partial_x f(u_h)\; dx+  \underbrace{ \sum_{{\sf f} \in\mathcal{F}_h}  \int_{\sf f} \tau_{\sf f} [\partial_x v_h] \cdot [\partial_x u_h] \; d\Gamma}_{S(v_h,u_h)}} = 0,    
\label{eq_cip0}
\end{equation}
\revMR{ with $[\cdot]$ denoting the jump of a quantity \revDT{across a face $\sf f$}, and   where we recall that  $\mathcal{F}_h$ is the collection of internal boundaries (points in 1D),  and ${\sf f}$ are its elements. 
In one space dimension the last integral   reduces to a point evaluation. 
Although other definitions are possible, we evaluate the scaling parameter in the  stabilization as 
\begin{equation}
\tau_{\sf f} =  \delta   \,h_{\sf f}^2 \| \partial_uf\|_{\sf f}
\label{eq_cip1}
\end{equation}
with $\| \partial_uf\|_{\sf f}$ a reference value of the norm of the flux Jacobian on ${\sf f}$ \revDT{and $h_{\sf f}$ a characteristic size of the mesh neighboring $\sf f$.}}

\revMR{ The advantage of this method is that the formulation remains symmetric, and that the mass matrix can be lumped for efficient time marching if the finite element space allows it. 
The drawback is a slight  increase in the stencil associated to the use of the gradients in all neighboring elements. 
Note that for higher order approximations  \cite{Burman2020ACutFEmethodForAModelPressure,larson2019stabilizationHighOrderCut}  suggest the use of jumps in 
higher derivatives 
to improve the stability of the method. In this work, we only focus on 
the gradient jump stabilization. For orders up to 4  this seems to be enough \revDT{to get $\mathbb{L}_2$ stability} and allows the  study in more detail the impact of the coefficient $\delta$ in the stabilization.  


As  before, we can easily characterize the  accuracy of the method following e.g. \cite[\S3.1.1 and \S3.2]{AR:17}, and show that 
for all  functions $\psi$ of class at least $\mathcal{C}^1(\Omega)$, of which  $\psi_h$  denotes the finite element projection, we have the truncation error estimate
\begin{equation}
\begin{split}
\epsilon(\psi_h) := \Big|
\int_{\Omega} \psi_h \partial_t (u_h^e - u^e) \; dx -& \int_{\Omega} \partial_x \psi_h (\partial_x f(u_h^e)-\partial_x f(u^e)\; dx \\ 
              +&  \sum\limits_{{\sf f}\in\mathcal{F}_h} \int\limits_{\sf f}\tau_{\sf f} [\partial_x \psi_h] \cdot [\partial_x (u_h^e-u^e)] 
	\Big| \le C h^{p+1},
\end{split}
\label{eq_cip2}
\end{equation}
\revDT{with $C$ \revsix{a constant} independent of $h$}. The estimate  is  again a direct consequence of standard approximation results applied to $u^e_h-u^e$ and to its derivatives.

The symmetry of the stabilization  makes is rather easy to derive a linear stability estimate. In particular, for a linear flux with periodic boundary conditions 
we can easily show that 
\begin{equation}
\begin{split}
\int\limits_{\Omega_h}\partial_t\dfrac{u^2_h}{2}= - \sum\limits_{{\sf f}\in\mathcal{F}_h}\int\limits_{\sf f} \tau_{\sf f} [\partial_x u_h]^2
\end{split}
\label{eq_cip3}
\end{equation}
which can be integrated in time to obtain a bound on  the $\mathbb{L}_2$ norm of the solution. }

\subsubsection{Local Projection Stabilization - LPS} \label{LPS}
Another symmetric stabilization approach is the  Local Projection Stabilization (LPS) method. 
Its original formulation  was presented in \cite{BB-LPS_2001} for Stokes equations. Then, the LPS was successfully extended to transport problems in \cite{BB_LPS-2004} and applications of local projection methods to 
Oseen and Navier-Stokes equations  were studied in \cite{BB_LPS-2006,AMTX-LPS_2010}. The local projection method  also aims at providing some   control on the fluctuations of the gradient of the discrete solution. 
\revMR{The method can be written as follows: find $u_h\in V_h^p$ such that $\forall v_h \in V_h^p$}

\revsix{
\begin{equation}
    \left  \{
    \begin{array}{ll}
    	& \int_{\Omega_h} v_h \partial_t u_h \; dx  + \int_{\Omega_h} v_h \partial_x f(u_h) \; dx  + \underbrace{\sum_{K \in \Omega_h} \int\limits_{K} \tau_K   \partial_x v_h (\partial_x  u_h - w_h) \; dx}_{S(v_h ,u_h)}= 0,  \\
    	& \int_{\Omega_h} v_h   w_h\; dx - \int_{\Omega_h} v_h \partial_x  u_h\; dx = 0. 
	\end{array}	
    \right .\label{eq_lps0}
\end{equation}
}
\revMR{
For this method, the stabilization parameter is evaluated as 
\begin{equation}
\label{eq_lps1}
\tau_K = \delta  h_K  \|\partial_u f\|_K  .
\end{equation}
Compared to the CIP approach this method has the drawback of requiring the mass matrix inversion in the gradient \revDT{$\mathbb{L}_2$} projection represented by the second \revDT{equation} in \eqref{eq_lps0}.
So the possibility of simplifying this operator, and, more precisely, to lump the mass matrix, appear as essential  elements for its efficient implementation.}

As before we can easily characterize the accuracy of this method.  The truncation error estimate for a polynomial approximation of degree $p$ reads in this case
\begin{equation}
\begin{split}
\epsilon(\psi_h) := \Big|
\int_{\Omega} \psi_h \partial_t (u_h^e - u^e) \; dx -& \int_{\Omega} \partial_x \psi_h (\partial_x f(u_h^e)-\partial_x f(u^e))\; dx \\ 
              +&  \sum\limits_{K \Omega_h} \int\limits_{K} \partial_x \psi_h ( \partial_xu^e_h -\partial_xu^e_h  )
                          +  \sum\limits_{K \Omega_h} \int\limits_{K} \partial_x \psi_h (  \partial_x u^e - w_h^e  )
	\Big| \le C h^{p+1},
\end{split}
\label{eq_lps2}
\end{equation}
where the last term is readily estimated using
$$
\int_{\Omega_h} \psi_h  ( w^e_h-\partial_x u^e)\; dx = \int_{\Omega_h} v_h (\partial_x  u_h^e -    \partial_x u^e ) \le  \mathcal{O}(h^p).
$$

Finally, for a linear flux and taking $\tau_K=\tau$, as for the SUPG, we can  test with  $v_h=u_h$ in the first of \eqref{eq_lps0}, and $v_h=\tau w_h$ in the second and sum up the result to get
(using the periodicity)
\begin{equation}
\begin{split}
\int\limits_{\Omega_h}\partial_t\dfrac{u^2_h}{2}= - \sum\limits_{K} \int\limits_{K} \tau ( \partial_x u_h -  w_h)^2,
\end{split}
\label{eq_lps3}
\end{equation}
which can be integrated in time to obtain a bound on  the $\mathbb{L}_2$ norm of the solution.

%

\subsection{Finite Element Spaces and Quadrature Rules}\label{sec:FEMquad}
We  describe the \revMR{one-dimensional finite element spaces we consider in the Fourier analysis. References to the corresponding multi-dimensional extensions are suggested for completeness where appropriate.}

\revMR{\revDT{In a one dimensional discretized space  $\Omega_h$} 
an element K is a    segment, i.\,e.,  $K=[x_i,x_{i+1}]$ for some $i$.}
We  \revMR{define in this section} the restriction of the basis functions of $V_h^p$ on each element $K$, which are polynomials of degree at most $p$. We denote with $\{\varphi_1, \ldots, \varphi_N\}$  the basis functions of $\P^p(K)$, and their definitions amounts to describe the degrees of freedom, i.e., the dual basis. In one dimension, $N=p+1$. \revMR{
We  consider two families of polynomials:}
\begin{enumerate}
	\item Lagrange polynomials.  They are uniquely defined by the interpolation points \revDT{$\xi_j$ with  $\xi_1=x_i<\ldots <\xi_j<\ldots <\xi_N=x_{i+1}$}. We study two cases
	\begin{itemize}
		\item Equidistant points: $\xi_j=x_i+j\frac{x_{i+1}-x_i}{p}$ for $j=0,\dots ,p$,
		\item Gauss--Lobatto points: \revMR{ the roots of } Legendre polynomial of degree $p+1$  mapped onto $[x_{i},x_{i+1}]$.
	\end{itemize}
	\item Bernstein polynomials. Linearly mapping $K$ onto  $[0,1]$ they are defined for $j=0,\dots ,p$ by
	$$B_{j}(x)= \begin{pmatrix}
		p\\ j
	\end{pmatrix} x^{p-j}(1-x)^j.$$
	Bernstein polynomials verify the following properties
	$$\sum_{j=0}^p B_{j}(x)\equiv 1, \qquad B_{j}(x)\geq 0 \quad \forall x \in [0,1].$$
	Even if the degrees of freedom associated to this approximation have no physical meaning, we identify them geometrically  with the Greville  points $\xi_j=\tfrac{j}{p}$. 
\end{enumerate}

The use of different polynomial basis functions leads to different properties. Let us remark that the evaluation of integrals is done by Gaussian quadrature formulae, because of their efficiency. If Gauss points are used in the discretization of the polynomials, the same points will be used in the quadrature formula. 
Thanks to this, we see that for Lagrange polynomials defined on Gauss quadrature points
$$\int_{x_i}^{x_{i+1}} \varphi_l(x)\varphi_j(x)\; dx=(x_{i+1}-x_i)\omega_l \delta_l^j\quad \text{with } \,\omega_l:=\frac{1}{(x_{i+1}-x_i)}\int_{x_i}^{x_{i+1}} \varphi_l^2(x)\ dx>0.$$ This leads to a diagonal local mass matrix 
$$\mass^i_{l,j}=\begin{pmatrix}
\int_{x_i}^{x_{i+1}} \varphi_l(x)\varphi_j(x)\; dx
\end{pmatrix}.$$
This does not hold for Lagrange polynomials defined on equidistant points or the Bernstein polynomials.\\

Another \revMR{important} property that we need 
to effectively apply the DeC method of \cite{paola_svetlana}  is the positivity of the lumped mass matrix \revMR{entries}, i.e., $\mathbb{D}_{k,k}:=\sum_{j=0}^N \int_{x_i}^{x_{i+1}} \varphi_j \varphi_k\; dx= \int_{x_i}^{x_{i+1}} \varphi_k \; dx >0$. The positivity of these values is trivially verified for Bernstein polynomials and for Lagrange polynomials with matching quadrature formulae. 
In the case of equispaced points Lagrangian polynomials, the lowest degree ($p\leq 7$ in one dimension) they also verify the positivity of the lumped matrix.
This is not true in the case of two dimensional problems and triangular meshes, where already for degree $p=2$ we have nonpositive values in the diagonal of the lumped matrix. 
This mainly motivated the choice of Bernstein polynomials, as well as the Lagrange interpolation with the Gauss--Lobatto points.

In the following we will use the wording
\begin{itemize}
	\item \textit{basic} elements for Lagrangian polynomials on equispaced points with Gauss--Legendre quadrature;
	\item \textit{cubature} elements for Lagrangian polynomials on on Gauss--Lobatto points and quadrature rule using the same points;
	\item \textit{Bernstein} elements for Bernstein polynomials with Gauss--Legendre quadrature.
\end{itemize}

\subsection{Time Integration}\label{sec:timeIntegration}

\revMR{The finite element semi-discrete equations constitute  a coupled system of ordinary differential equations which can be written as }
\begin{equation}\label{eq:linear_system}
\revMR{  \mass  \dfrac{dU}{dt}  = \rhs(t) }
\end{equation}
where $U$ 
is the collection of all the degrees of freedom, $\mass$ and $\rhs$ are the global mass matrix and right-hand side term defined in the previous sections through the element definition and stabilization terms.
We must remark that $\mass$ is diagonal only in the case of the \textit{cubature} elements without the SUPG stabilization, while, for all other choices, it is a sparse non--diagonal matrix.

In the following, we describe \revMR{two different time integration  strategies: }
explicit Runge--Kutta (RK) methods and their strong stability preserving (SSP) variant;  
 Deferred Correction, which allows to \revMR{\revDT{avoid} the mass matrix inversion \revDT{through} the correction iterations.} 
\subsubsection{Explicit Runge--Kutta and Strong Stability Preserving Runge--Kutta schemes} 
\revDT{Runge--Kutta time integration methods} can be described by the following one step procedure 
\begin{equation}
\begin{split}
	&U^{(0)}:=U^n,\\
	&U^{(s)}:=U^n + \Delta t \sum_{j=0}^{s-1}\alpha_j^s \mass^{-1} \rhs(U^{(j)})\quad s=1,\dots, S,\\
	&U^{n+1}:= U^n +  \Delta t  \sum_{s=0}^S \beta_s \mass^{-1} \rhs(U^{(s)}).
	\end{split}\label{eq:RK}
\end{equation}
Here, we use the superscript $n$ to indicate the timestep and the superscript in brackets $(s)$ to denote the stage of the method.
In particular, we will refer to Heun's method with RK2, to Kutta's method with RK3 and the original Runge--Kutta fourth order method as RK4. The respective Butcher's tableau can be found in \cref{app:timeCoefficients} in \cref{tab:Butcher}.

\revMR{
A particular case is that of  SSPRK methods  introduced in \cite{shu-1988}. They 
are essentially convex combinations of forward Euler steps,  and 
can be rewritten as follows}
\begin{equation}\label{eq:SSPRKformula}
	\begin{split}
	&U^{(0)}:=U^n,\\
	&U^{(s)}:=\sum_{j=0}^{s-1} \left( \gamma_j^s U^{(j)} + \Delta t \mu_j^s \mass^{-1} \rhs(U^{(j)}) \right) \quad s=1,\dots, S,\\
	&U^{n+1}:= U^{(S)} ,
	\end{split}
\end{equation}
with $\gamma_j^s, \mu_j^s\geq 0$ for all $j,s=1,\dots, S$. 
We will consider here 
 the second order
 3 stages SSPRK(3,2) presented by Shu and Osher in \cite{shu-1988}, the third order 
 SSPRK(4,3) presented in \cite[Page 189]{Ruuth-2006}, 
 and the  fourth order 
 SSPRK(5,4) defined in \cite[Table 3]{Ruuth-2006}.
For complete reproducibility of the results, we put all their Butcher' tableaux in \cref{app:timeCoefficients} in \cref{tab:ButcherSSPRK}.
%

\subsubsection{The \textit{Deferred Correction} scheme}
\revMR{\revDT{Deferred} correction methods}  were originally introduced in \cite{ecdd148d37dd401ca8d415ae18d3ecf2} as explicit solvers of ODEs, but soon implicit \cite{DeC_ODE} or positivity preserving  \cite{oeffner_torlo_2019_DeCPatankar} versions and extensions to PDE solvers \cite{DeC_2017} were studied.
 In \cite{DeC_2017,DeC_AT,paola_svetlana} the method is also used to avoid the inversion of the mass matrix, applying a mass lumping and adding correction iterations to regain the order of convergence.
This is only achievable when the lumped matrix have only positive values on its diagonal. Hence, the use of \textit{Bernstein} polynomials is recommended in \cite{DeC_2017}, but also the \textit{cubature} elements can serve the purpose.

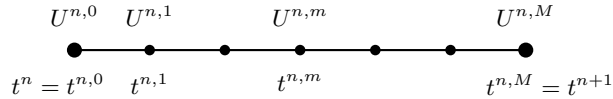
\begin{figure}[h]
	\small	\centering
	\begin{tikzpicture}
	\draw [thick]   (0,0) -- (6,0) node [right=2mm]{};
	\fill[black]    (0,0) circle (1mm) node[below=2mm] {$t^n=t^{n,0} \,\, \quad$} node[above=2mm] {$U^{n,0}$}
	(1,0) circle (0.7mm) node[below=2mm] {$t^{n,1}$} node[above=2mm] {$U^{n,1}$}
	(2,0) circle (0.7mm) node[below=2mm] {}
	(3,0) circle (0.7mm) node[below=2mm] {$t^{n,m}$} node[above=2mm] {$U^{n,m}$}
	(4,0) circle (0.7mm) node[below=2mm] {}
	(5,0) circle (0.7mm) node[below=2mm] {}
	(6,0) circle (1mm) node[below=2mm] {$\qquad t^{n,M}=t^{n+1}$} node[above=2mm] {$U^{n,M}$};
	\end{tikzpicture}
	\caption{Subtimesteps inside the time step $[t^n,t^{n+1}]$} \label{fig_DeC-time_disc}
\end{figure}


Consider a discretization of each timestep into $M$ \revDT{subtimesteps}	
as in  \cref{fig_DeC-time_disc}.
For each subtimestep the goal is to find the solution of the integral form of the semidiscretized ODE \eqref{eq:linear_system} as 
\begin{equation}\label{eq:L2}
	\mass \left( U^{n,m} -U^{n,0} \right) - \int_{t^{n,0}}^{t^{n,m}} \rhs (U(s)) ds \approx \L^2(\UU)^m:= \mass \left( U^{n,m} -U^{n,0} \right) - \Delta t \sum_{z \in \llbracket 0, M \rrbracket } \rho_{z}^m \rhs (U^{n,z}) = 0,
\end{equation}
with $\UU=\left( U^{n,0},\dots, U^{n,M} \right)$ and  \revMR{having used} 
 high order quadrature 
 with 
 points $t^{n,0},\dots, t^{n,M}$ and weights $\rho^{m}_z$ for every different subinterval   
 (see  \cite{DeC_2017,DeC_AT,paola_svetlana} for details).  \revMR{The algebraic system  }
$\L^2(\UU^*)=0$ \revMR{is in general implicit and nonlinear} 
and may not be easy to solve. The DeC procedure \revMR{approximates iteratively} 
this solution 
by successive corrections \revMR{relying on a} 
a low order easy--to--invert operator $\L^1$.  This operator is \revMR{typically} obtained 
using an explicit timestepping and a lumped mass matrix, i.e.,
\begin{equation}\label{eq:L1}
\mass \left( U^{n,m} -U^{n,0} \right) - \int_{t^{n,0}}^{t^{n,m}} \rhs (U(s)) ds \approx \L^1(\UU)^m:= \diag \left( U^{n,m} -U^{n,0} \right) - \Delta t \beta^m \rhs (U^{n,0}) = 0.
\end{equation}
Here, $\diag$ denotes a diagonal matrix obtained from the lumping of $\mass$, i.e., $\diag_{ii}:=\sum_{j} \mass_{ij}$, and $\beta^m:= \frac{t^{n,m}-t^{n,0}}{t^{n+1}-t^n}$. The values of the coefficients $\beta^m$ and $\rho^m_z$ for equispaced subtimesteps can be found in \cref{app:timeCoefficients}. 
Denoting with the superscript $(k)$ index the iteration step, we describe the DeC algorithm as
\begin{subequations}
\begin{align}
	&U^{n,m,(0)}:=U^n & m=0,\dots,M,&\\
	&U^{n,0,(k)}:=U^n  & k=0,\dots, K,&\\
	&\L^1(\UU^{(k)})=\L^1(\UU^{(k-1)})-\L^2(\UU^{(k-1)})& k=1,\dots, K,& \label{eq:DeCUpdate}\\
	&U^{n+1}:=U^{n,M,(K)}.&&
\end{align}
\end{subequations}

It has been proven \cite{DeC_2017} that if $\L^1$ is coercive, $\L^1-\L^2$ is Lipschitz with a constant $\alpha_1 \Delta t >0$ and the solution of $\L^2(\UU^*)=0$ exists and is unique, then, the method converges with an error of $\mathcal{O}(\Delta t^K)$. Hence, choosing $K=M+1$ we obtain a $K$-th order accurate scheme.

\revMR{Relying only on the inversion of the  the low order operator, the method has for each iteration a cost equivalent essentially to the assembly of the right hand side, whatever the complexity of the mass matrix
appearing in $\L^2$.}
The only requirement that is necessary for the DeC approach is the invertibility of the lumped mass matrix, 
which limits its application  to  equispaced Lagrange elements only to the  degrees for which this is the case, and to other choices as the  
 \textit{Bernstein}  and \textit{cubature} elements \revMR{introduced earlier}. 

Finally,  \revMR{for the following analysis we note that} 
the DeC method can be \revMR{cast in a form similar to} 
a Runge--Kutta method \revMR{by  rewriting} 
 \eqref{eq:DeCUpdate} 
 as
\begin{equation}\label{eq:DeCasSSPRK}
	U^{n,m,(k+1)}=U^{n,m,(k)} - \diag^{-1} \mass \left(U^{n,m,(k)}-U^{n,0,(k)}\right) +\sum_{j=0}^M \Delta t \rho_{j}^m \diag^{-1}\rhs(U^{n,j,(k)}). 
\end{equation}
\revMR{Comparing with }
\eqref{eq:SSPRKformula}, \revMR{we can immediately define the SSPRK coefficients associated to DeC as} 
$\gamma^{m,(k+1)}_{m,(k)}=\mathbb{I}-\diag ^{-1} \mass$ with $\mathbb{I}$ the identity matrix,  $\gamma^{m,(k+1)}_{0,(0)}=\diag ^{-1} \mass$, $\mu^{m,(k+1)}_{r,(k)}=\rho^m_r$ for $m,r=0,\dots,M$ and $k=0,\dots,K-1$ and instead of the mass matrix, we use the diagonal one. 
\section{Fourier Analysis}\label{sec:fourier}
The dispersion and the stability properties of numerical methods can be shown \revMR{by means of a spectral analysis. We will focus on the linear case \eqref{eq:conservation_law1}  with periodic boundary conditions:}
\begin{equation}
    \partial_t u + a \partial_x u = 0, \quad x\in[0,1].
    \label{eq_disp_1}
\end{equation}
\revMR{The main idea is to  investigate the semi and fully discrete evolution of periodic waves represented by the 
  the ansatz}
\begin{align}
    u &= Ae^{i(kx - \xi t)} = Ae^{i(kx-\omega t)}e^{\epsilon t} \qquad \mbox{with} \quad \xi = \omega + i \epsilon, \quad i=\sqrt{-1}. \label{eq_disp_uex1} 
\end{align}
Here, $\epsilon$ denotes the damping rate, while \revMR{ the wavenumber is denoted by   $k=2\pi/L$  with} $L$ the wavelength. 
\revMR{We recall that the phase velocity defined as} 
\begin{equation}
    C =  \frac{\omega}{k}
    \label{eq_PDE_velocity}
\end{equation}
\revMR{represents the celerity with which waves propagate in space, and it is in general a function of the wavenumber.}
Substituting 
\refp{eq_disp_uex1} in the advection equation \refp{eq_disp_1} leads to \revMR{the well known result} 
\begin{align}
    &  \quad C = a \quad \mbox{and} \quad \epsilon = 0.
\end{align}
\revMR{The objective of the next sections is to provide the semi and fully discrete equivalents of the above relations for the finite element methods introduced earlier.
We will consider  polynomial degrees up to 3, for all combinations  of different stabilization methods and time integration. This will also allow to investigate the parametric stability with respect to
the time step (\CFL number) and stabilization parameter $\delta$. In practice, for each choice we will evaluate the accuracy of the discrete
approximation   of  $\omega$ and $\epsilon$,  and we will provide conditions for the non-positivity of    the damping $\epsilon$. For completeness, the study is performed first  in the semi-discrete time continuous case  in
 \cref{sec:fourier_space}. We the consider the fully discrete schemes in  \cref{sec:fourier_space_time}.}%

\subsection{Preliminaries and time continuous analysis}
\label{sec:fourier_space}


The Fourier analysis for numerical schemes on the periodic domain  is based on Parseval theorem. 
\begin{theorem}[Parseval]
	Let $\hat{u}(k):= \int_{0}^{1} u(x)e^{-i2\pi kx} dx$ for $k\in \mathbb{Z}$ be the Fourier modes of the function $u$. The $\mathbb{L}_2$ norms of the function $u$ and of the Fourier modes coincide, i.e.,
	\begin{equation}
		\int_{0}^{1} u^2(x)dx = \sum_{k\in \mathbb{Z}} |\hat{u}(k)|^2.
	\end{equation}
\end{theorem}
Thanks to this theorem, we can study the amplification and the dispersion of the basis functions of the Fourier space. 
The key ingredient of this study is the repetition of the stencil of the scheme from one cell to another one. \revMR{In particular, }
 using the ansatz \eqref{eq_disp_uex1} 
 we can  \revMR{write local equations coupling} 
  degrees of freedom belonging to neighbouring cells through a multiplication by \revMR{the factor of  $e^{i\theta}$ representing the shift in space along the oscillating solution. 
  The dimensionless  coefficient
  	\begin{equation}\label{eq_theta}
\theta:=  k\Delta x
	\end{equation}
	is a discrete reduced wave number which naturally appears all along the analysis.
Formally replacing the ansatz in the scheme we end up with a dense  algebraic problem of dimension $p$ (the polynomial degree)  
reading in the time continuous case }
\begin{equation}
 \eqref{eq_disp_1} \text{ \revDT{and} }\eqref{eq_disp_uex1} \quad \Rightarrow \quad - i\xi \mass \mathbf{U}    + a \mathcal{K}_x \mathbf{U} = 0 \label{eq_sys_mat}
\end{equation}
\begin{equation}
\hspace*{-1cm}
    \mbox{with} \quad (\mass)_{ij} = \int_{\Omega} \phi_i \phi_j dx, \qquad (\mathcal{K}_{x})_{ij} = \int_{\Omega} \phi_i \partial_x \phi_j dx + S(\phi_i,\phi_j), 
\end{equation} 
\revMR{with $\phi_j$ the finite element basis functions and $ \mathbf{U} $ the array of all the degrees of freedom. Although system \eqref{eq_sys_mat}   is in general a global eigenvalue problem, we can reduce its complexity by 
exploiting more explicitly the ansatz \eqref{eq_disp_uex1}. More exactly, we can introduce elemental vectors of unknowns $\widetilde{\mathbf{U}}_K$, which,  for continuous finite elements, are a arrays of $p$ degrees of freedom
including only one of the two boundary nodes. Using the periodicity of the solution and   denoting by $K\pm1$ the neighboring elements, we have
\begin{equation}
\widetilde{\mathbf{U}}_{K\pm1} = e^{\pm\theta}\widetilde{\mathbf{U}}_{K}.
\label{eq_ad_fourier0}
\end{equation}
This allows to show that  \eqref{eq_sys_mat}  is equivalent to a  compact system (we drop the subscript $_K$ as they system is equivalent for all cells)
\begin{equation}
-i\xi \widetilde{\mass}   \widetilde{\mathbf{U}}  + a \widetilde{\mathcal{K}}_x \widetilde{\mathbf{U} }=0,
\label{eq_ad_fourier}
\end{equation}
where   the matrices $\widetilde{\mass}$ and $ \widetilde{\mathcal{K}}$ are readily obtained from the elemental discretization matrices by using  \eqref{eq_ad_fourier0}.

 As shown in \cite{sherwin_dispersion} some particular cases  can be easily studied analytically.
 For example for the semidiscretized $\mathbb{P}_1$ CG scheme without stabilization one easily finds that}
\begin{equation} 
\frac{\omega}{k} = a \frac{\sin(\theta)}{\theta}\frac{3}{2+\cos(\theta)} \quad \mbox{and} \quad \epsilon = 0.
\end{equation}
\revMR{As the degree of the approximation increases, so does the size of the eigenvalue problem. For the non stabilized CG $\P_2$ scheme we can still find an analytical solution associated to the quadratic equation (cf also \cite{sherwin_dispersion})} reading 
\begin{equation} \hspace{-0.5cm}
\frac{\omega_{1,2}}{k} = a \frac{4\sin(\theta)\pm 2 \sqrt{40\sin^2 (\frac{\theta}{2})-\sin^2(\theta)}}{\theta(\cos(\theta)-3)}.
\end{equation}
For more general cases, the study needs to be performed numerically.  \\

Defining with $\lambda_i(\theta)$ the eigenvalues of \refp{eq_ad_fourier}, $\omega_i(\theta) = \text{Im}(\lambda_i(\theta))$ and $\epsilon_i(\theta) =- \text{Re}(\lambda_i(\theta))$ are the respective phase and damping coefficients of each mode of the solution.  In practice, we solve numerically the eigenvalue problem \eqref{eq_ad_fourier} for $\theta=k \Delta x_p = \frac{2\pi}{N_x}$ varying in $[0,\pi]$, where $N_x$ is the number of the nodes in each wavelength and $\Delta x_p=\Delta x/p$ is the average distance between degrees of freedom.  However, to satisfy the Nyquist stability criterion, it is necessary to have $\Delta x_p \leq \frac{L}{2}$, with $L$ the wavelength.\\

\begin{figure}[H]
     \centering
     \includegraphics[width=0.9\textwidth]{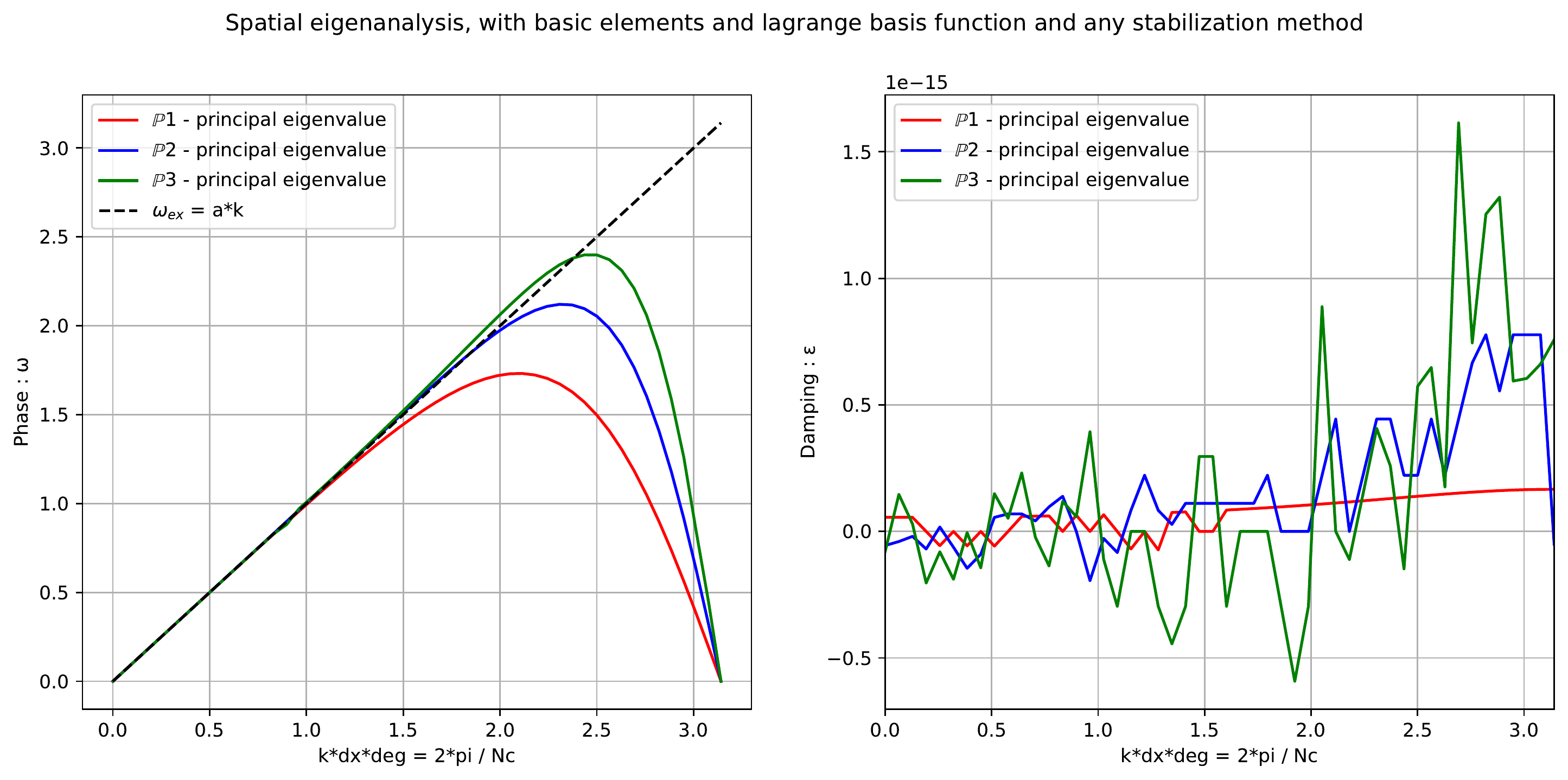}
     \caption{Phase $\omega$ (left) and amplification $\epsilon$ (right) with \textit{basic} elements without stabilization for $\mathbb{P}_1, \mathbb{P}_2$ and $\mathbb{P}_3$.}
     \label{fig_disp_lag_CT_DS-any}
\end{figure}

\revMR{As an example,} in \cref{fig_disp_lag_CT_DS-any} we plot $\omega$ and $\epsilon$ and we see that CG scheme does not have diffusive terms, or, in other words, there is no damping ($\epsilon =0$) in the CG scheme.  
For clarity of the pictures, we plot in \cref{fig_disp_lag_CT_DS-any} only the principal eigenvalue of each system ($p=1,2,3$), \revDT{i.\,e., the one that minimizes $|\omega_i-ak|$.} 
As expected, with $\mathbb{P}_1$ elements, the scheme is more dispersive than with $\mathbb{P}_2$ or $\mathbb{P}_3$ elements, while, for all of them, there is no dissipation, since the scheme is not stabilized and there is no time discretization. \\



\begin{figure}[H]
     \begin{center}
     \includegraphics[width=0.9\textwidth]{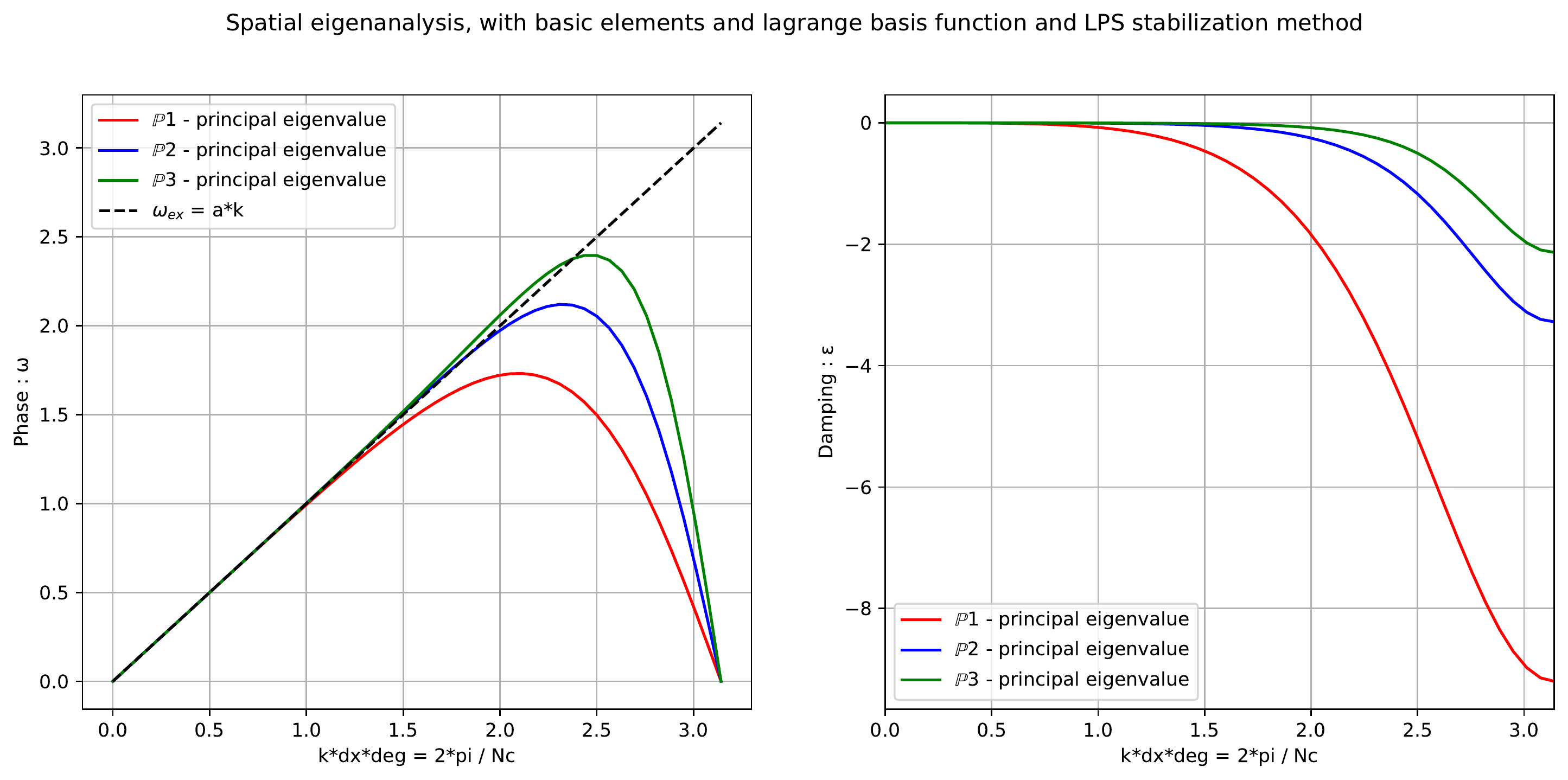}
     \end{center}
     \caption{Phase $\omega$ (left) and amplification $\epsilon$ (right) with \textit{basic} elements with LPS stabilization for $\mathbb{P}_1, \mathbb{P}_2$ and $\mathbb{P}_3$.}
     \label{fig_disp_lag_CT_DS-lps}
\end{figure}

\revMR{We apply the same analysis to stabilized methods. The results obtained with SUPG, CIP and LPS stabilizations   lead to an almost identical result shown in \cref{fig_disp_lag_CT_DS-lps}
(reporting the LPS data). The interested reader can  access all the other plots online  \cite{TorloMichel2020git}.  From the  plot  we can see that the increase in polynomial degree
provides the   expected large reduction in dispersion error,   while retaining a small  amount of }
 numerical dissipation, which permits the damping of \textit{parasite} modes. 

\subsection{Fully discrete analysis}\label{sec:fourier_space_time}



\subsubsection{Methodology}
We \revMR{analyze now  
the fully discrete  schemes obtained  using the RK, SSPRK and DeC time marching methods  presented in \cref{sec:timeIntegration}. }
Let us consider as an example the SSPRK schemes \eqref{eq:SSPRKformula}. If we define as $A:=\mass^{-1} \mathcal{K}_x$ we can write the schemes as follows
\begin{equation}
\left  \{
    \begin{array}{ll}
    	\mathbf{U}^{(0)}:= & \mathbf{U}^n \\
        \mathbf{U}^{(s)} := & \sum_{j=0}^{s-1} \left( \gamma_{sj} \mathbf{U}^{(j)} + \Delta t \mu_{sj} A \mathbf{U}^{(j)} \right), \quad s \in \llbracket 1,S \rrbracket ,\\
	    \mathbf{U}^{n+1}:=&\mathbf{U}^{(S)}.
	\end{array}
    \right .
    \label{eq:discreteSSPRK}
\end{equation}

Expanding all the stages, we can obtain the following formulation:
\begin{equation}
\mathbf{U}^{n+1} = \mathbf{U}^{(0)} + \sum_{j=1}^{S} \nu_{j} \Delta t^jA^j  \mathbf{U}^{(0)} = \left(\mathcal{I} + \sum_{j=1}^{S} \nu_{j}\Delta t^j A^j  \right) \mathbf{U}^n, \label{RK_gamma0}
\end{equation}
where coefficients $\nu_j$ in \eqref{RK_gamma0} are obtained as combination of coefficient $\gamma_{sj}$ and $\mu_{sj}$ in \eqref{eq:discreteSSPRK} and $\mathcal{I} $ is the identity matrix. For example, coefficients of the fourth order of accuracy scheme \textit{RK4} are $\nu_1=1$, $\nu_2 = 1/2$, $\nu_3=1/6$ and $\nu_4 = 1/24$.\\ 
%

\revMR{We can now compress the problem proceeding as in the time continuous case. In  particular, using \eqref{eq_ad_fourier0} 
one easily shows that the problem can be written in terms of the local $p\times p$ matrices $\widetilde{A}:= a \widetilde{\mass}^{-1} \widetilde{\mathcal{K}_x}$ 
and in particular that
%
\begin{align*}
\widetilde{\mathbf{U}}^{n+1} = G \widetilde{\mathbf{U}}^{n}\quad\text{with}\quad
G:= e^{\epsilon \Delta t } e^{-i\omega \Delta t} & \approx  \left(\widetilde{\mathcal{I}} + \sum_{j=1}^{S} \nu_{j} \Delta t^j \widetilde{A}^j \right) ,
\end{align*}
where $G\in \R^{p\times p}$ }is the amplification matrix depending on $\theta, \Delta t$ and $ \Delta x$. Considering each eigenvalue $\lambda_i$ of $G$, we can write the following formulae for the corresponding phase $\omega_i$ and damping coefficient $\epsilon_i$
\revDT{\begin{align*}
\begin{cases}
   	e^{\epsilon_i \Delta t } \cos(\omega_i \Delta t)  = \text{Re}(\lambda_i) ,\\
   	- e^{\epsilon_i \Delta t } \sin(\omega_i \Delta t)  = \text{Im}(\lambda_i),
\end{cases} 
\Leftrightarrow \, \begin{cases} 
\omega_i\Delta t = \arctan \left( \frac{-\text{Im}(\lambda_i)}{\text{Re}(\lambda_i)} \right)  ,\\
(e^{\epsilon_i \Delta t })^2 = \text{Re}(\lambda)^2 + \text{Im}(\lambda)^2,
\end{cases} \Leftrightarrow \, \begin{cases} 
\dfrac{\omega_i}{k} = \arctan \left( \frac{-\text{Im}(\lambda_i)}{\text{Re}(\lambda_i)} \right)  \frac{1}{k \Delta t},\\
 \epsilon_i = \log\left( | \lambda_i | \right) \frac{1}{\Delta t}.
\end{cases}
\end{align*}
For the DeC }method we can proceed with the same analysis transforming also the other involved matrices into their Fourier equivalent ones. 
Using \eqref{eq:DeCasSSPRK} these terms would contribute to the construction of $G$ not only in the $\widetilde{A}$ matrix, but also in the coefficients $\nu_j$, which become matrices as well.
At the end we just study the final matrix $G$ and its eigenstructure, whatever process was needed to build it up.

The matrix $G$ represents the evolution in one timestep of the Fourier modes for all the $p$ different types of degrees of freedom. The damping coefficients $\epsilon_i$ tell if the modes are increasing or decreasing in amplitude and the phase coefficients $\omega_i$ describe the phases of such modes.

We remark that a necessary condition for stability of the scheme is that $ |\lambda_i | \leq 1$ or, equivalently, $\epsilon_i \leq 0 $ for all the eigenvalues. The goal of our study is to find the largest CFL number for which the stability condition is fulfilled and such that the dispersion error is not too large. 
Furthermore, we notice that the matrix $G$ depends not only on $\theta, \Delta x$ and $\Delta t$, but also on at the stabilization coefficients $\tau_K$. 
Hence, the proposed analysis should contain an optimization process also along the stabilization parameter. 
\revMR{With the notation of section \S2, we will in particular set} 
\begin{equation*}
\begin{split}
\quad \text{SUPG :}\;\;& \tau_K =\delta  \Delta x/|a|,\\[5pt]
\quad \text{LPS :} \;\;& \tau_K =\delta  \Delta x  | a |,\\[5pt]
\quad\text{CIP :}   \;\;& \tau_f = \delta\Delta x^2  | a |.
\end{split} 
\end{equation*}

\revMR{One of our objectives is to explore the space of parameters (CFL,$\delta$), and to propose criteria allowing to set these parameters to}
provide the most stable, least dispersive and least expensive methods.
\revMR{A clear and natural criterion is to exclude all parameter values for which we obtain  } 
a positive damping coefficient $\epsilon(\theta)>10^{-12}$ for any  value of the reduced wavenumber $\theta$ (taking into account the machine precision errors that might occur). \revDT{Doing so, we obtain what we will denote as \textit{stable area} in $(\CFL,\theta)$ space. }
For all the other points we \revDT{propose 3 strategies to minimize the product between error and computational cost. 
In the following we describe the 3 strategies to find the best parameters couples (\CFL,$\delta$):}\revMR{
\begin{enumerate}
\item  {\it maximize the {\CFL} in the stable area};
\item {\it minimize a global solution error, denoted by $\eta_u$,  while maximizing the {\CFL} in the stable area.} In particular, we start from  the relative square error of $u$ 
\begin{align}
	\left[\frac{u(t)-u_{ex}(t)}{u_{ex}(t)}\right]^2= &\left[e^{\epsilon t - i t(\omega-\omega_{ex})}-1\right]^2\\
	=&\left[e^{\epsilon t}\cos(t(\omega-\omega_{ex}))-1\right]^2+\left[e^{\epsilon t}\sin(t(\omega-\omega_{ex}))\right]^2\\
	=&e^{2\epsilon t} - 2 e^{\epsilon t} \cos (t(\omega-\omega_{ex})) +1.
\end{align}
Here, we denote with $\epsilon$ and $\omega$ the damping and phase of the \textit{principal }mode.
For a small enough  dispersion error $|\omega-\omega_{ex} |\ll 1$, we can expand the cosine in the previous formula in a truncated Taylor series as 
\begin{align}
	\left[\frac{u(t)-u_{ex}(t)}{u_{ex}(t)}\right]^2\approx&\underbrace{\left[e^{\epsilon t} -1\right]^2}_{\text{Damping error}} + \underbrace{e^{\epsilon t}t^2 \left[\omega- \omega_{ex}\right]^2}_{\text{Dispersion error}}.
\end{align}
We then compute an   error at the final time $T=1$, over the whole phase domain, using   at least 3 points per wave $0\leq k \Delta x_p \leq \frac{2\pi}{3}$, with $\Delta x_p=\frac{\Delta x}{p}$, and   $p$ the degree of the polynomials. 
We  obtain the following $\mathbb{L}_2$ error definition, 
\begin{equation}
	\eta_u(\omega,\epsilon)^2:= \frac{3}{2\pi} \left[\int_{0}^{\frac{2\pi}{3}} (e^{\epsilon}-1 )^2 dk + \int_{0}^{\frac{2\pi}{3}} e^\epsilon(\omega-\omega_{ex})^2 dk \right].
\end{equation}
Recalling that $\epsilon=\epsilon(k\Delta x,\CFL,\delta)$ and $\omega=\omega(k,\Delta x,\CFL,\delta)$ and $\omega_{ex}=ak$, \revDT{we need to further set the parameter $\Delta x_p$. We choose it to be large $\Delta x_p=1$, with the hope that for finer grids the error will be smaller.} Finally,  we seek the  couple $(\text{CFL}^*,\delta^*)$ allowing to solve
\begin{equation}
\label{cfl_d}
	(\text{CFL}^*,\delta^*):=\arg \max_{\text{CFL}} \left\lbrace \eta(\omega(\text{CFL},\delta)),\epsilon(\text{CFL},\delta)) < \mu \min_{(\text{CFL},\delta) \text{stable}} \eta(\omega(\text{CFL},\delta),\epsilon(\text{CFL},\delta))\right\rbrace.
\end{equation}
\item {\it minimize the dispersion error $\eta_\omega$ while maximizing the CFL in the stable area.} In particular we set in this case 
\begin{equation}
	\eta_\omega ^2(\omega):= \int_{0}^{\frac{2\pi}{3}} \left( \frac{\omega-\omega_{ex}}{\omega_{ex}}\right)^2 dk.
\end{equation}
As before we choose the optimal parameters from \eqref{cfl_d}.
\end{enumerate}
For the second and third strategies, the parameter $\mu$ must be chosen in order to balance the  requirements on stability and accuracy. 
After having tried different values,  we  have  set $\mu$ to $1.3$ providing a sufficient flexibility to obtain results of practical usefulness, which we verified in numerical computations as we will see later.}


In the following we will compare all the methods with these error measures, in order to suggest the best possible schemes between the proposed ones.

\section{Results of the  fully discrete spectral analysis}\label{sec:spectralAnalysis}
\begin{figure}
	\centering
	\includegraphics[width=\textwidth]{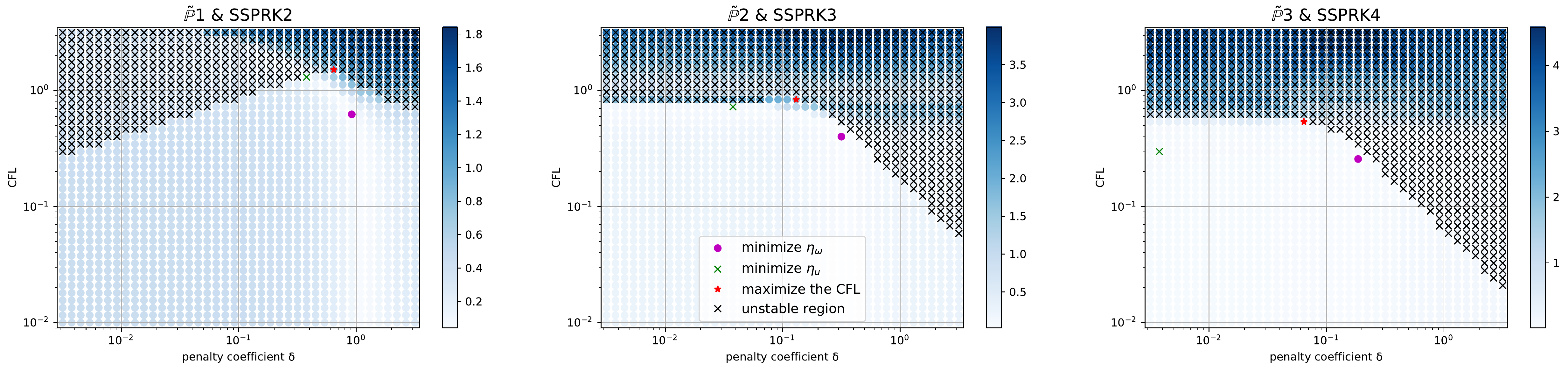}
	\caption{Computation of optimal parameters according to errors $\eta_\omega$ and $\eta_u$. $(\CFL,\delta)$ plot of $\eta_u$ (blue scale) and instability area (black crosses) for cubature elements SSPRK scheme with SUPG stabilization method. From left to right $\mathbb{P}_1$, $\mathbb{P}_2$, $\mathbb{P}_3$. The purple circle is the optimizer of $\eta_u$, the green cross is the optimizer of $\eta_\omega$, the red star is the maximum stable CFL. }
	\label{fig_disp_cohen_FD_SSPRK-cfl_vs_tau-SUPG}
\end{figure}

\begin{figure}
	\centering
	\includegraphics[width=\textwidth]{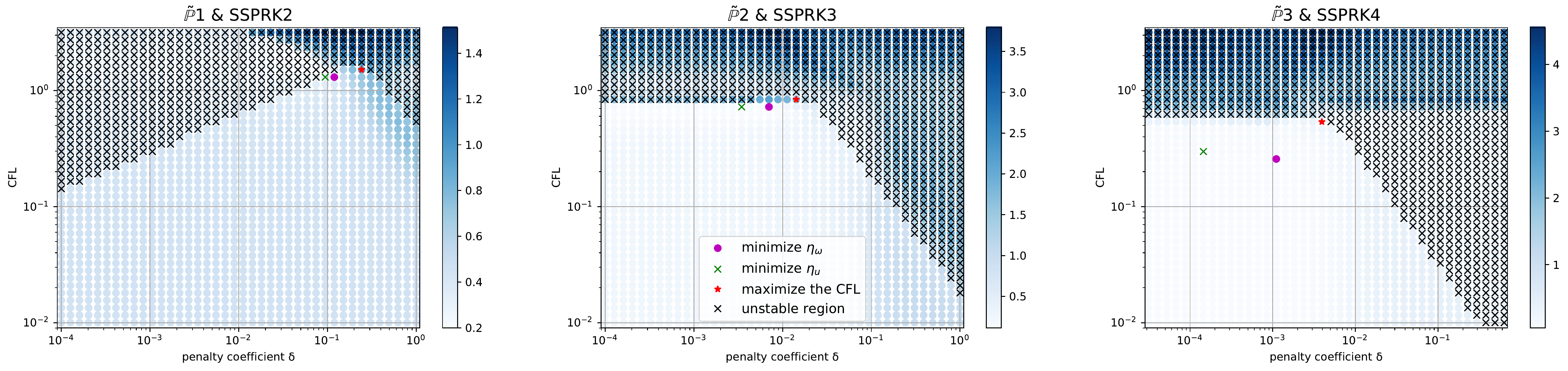}
	\caption{Computation of optimal parameters according to errors $\eta_\omega$ and $\eta_u$. $(\CFL,\delta)$ plot of $\eta_u$ (blue scale) and instability area (black crosses) for cubature elements SSPRK scheme with CIP stabilization method. From left to right $\mathbb{P}_1$, $\mathbb{P}_2$, $\mathbb{P}_3$. The purple circle is the optimizer of $\eta_u$, the green cross is the optimizer of $\eta_\omega$, the red star is the maximum stable CFL. }
	\label{fig_disp_cohen_FD_SSPRK-cfl_vs_tau-CIP}
\end{figure}

\begin{figure}
	\centering
	\includegraphics[width=\textwidth]{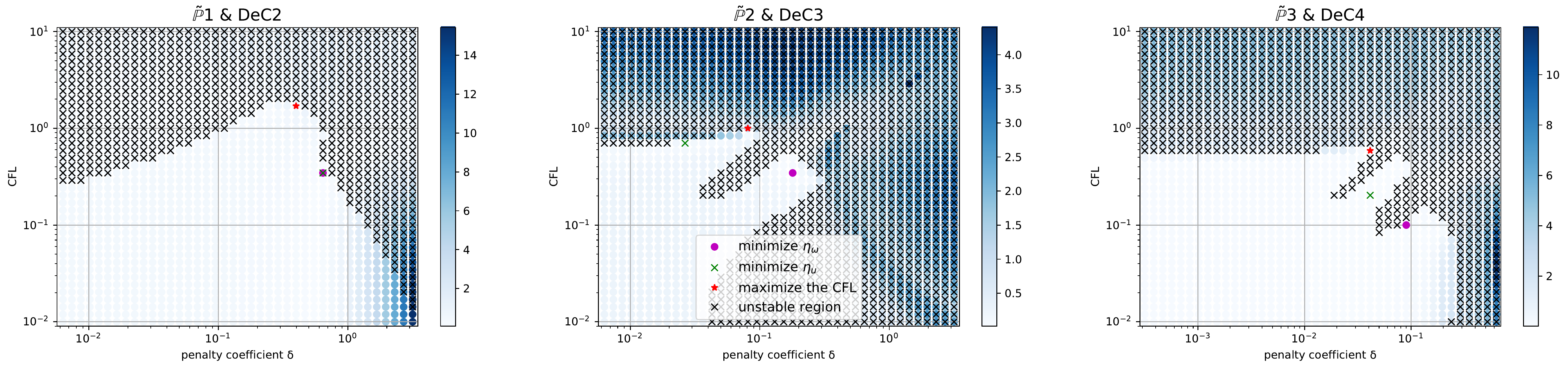}
	\caption{Computation of optimal parameters according to errors $\eta_\omega$ and $\eta_u$. $(\CFL,\delta)$ plot of $\eta_u$ (blue scale) and instability area (black crosses) for cubature elements DeC scheme with SUPG stabilization method. From left to right $\mathbb{P}_1$, $\mathbb{P}_2$, $\mathbb{P}_3$. The purple circle is the optimizer of $\eta_u$, the green cross is the optimizer of $\eta_\omega$, the red star is the maximum stable CFL. }
	\label{fig_disp_cohen_FD_DeC-cfl_vs_tau-SUPG}
\end{figure}

\begin{figure}
	\centering
	\includegraphics[width=\textwidth]{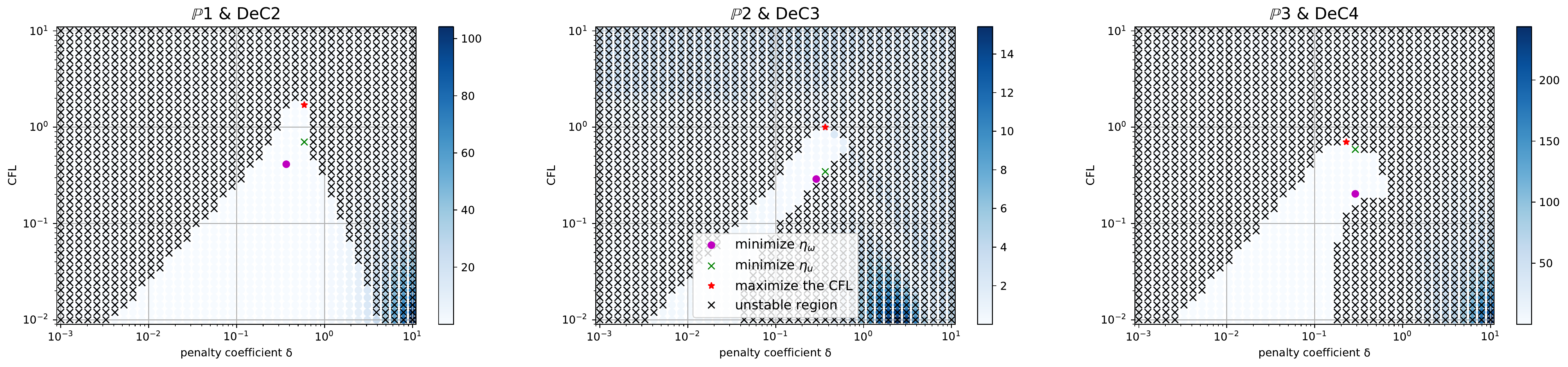}
	\caption{Computation of optimal parameters according to errors $\eta_\omega$ and $\eta_u$. $(\CFL,\delta)$ plot of $\eta_u$ (blue scale) and instability area (black crosses) for Bernstein elements DeC scheme with SUPG stabilization method. From left to right $\mathbb{P}_1$, $\mathbb{P}_2$, $\mathbb{P}_3$. The purple circle is the optimizer of $\eta_u$, the green cross is the optimizer of $\eta_\omega$, the red star is the maximum stable CFL. }
	\label{fig_disp_bezier_FD_DeC-cfl_vs_tau-SUPG}
\end{figure}

\revMR{The typical   results  reported	in   \cref{fig_disp_cohen_FD_SSPRK-cfl_vs_tau-SUPG,fig_disp_cohen_FD_SSPRK-cfl_vs_tau-CIP,fig_disp_cohen_FD_DeC-cfl_vs_tau-SUPG,fig_disp_bezier_FD_DeC-cfl_vs_tau-SUPG,fig_disp_lagrange_FD_DeC-cfl_vs_tau-LPS} show in the plane $(\delta,\CFL)$ the unstable (crossed) and stable regions, and with colored symbols the optimal points corresponding to the three strategies introduced earlier.}
\revMR{In case of ambiguity, the point with maximum $\delta$  is marked in the figures.} 
A summary of the results for all combinations of schemes is    provided in \cref{tab:dispersion_cfl-RES,tab:dispersion_cfl-REW,tab:dispersion_cfl-best}.

\revMR{Before commenting these results we remark that some of the schemes are equivalent. For example without mass lumping \textit{Bernstein} and \textit{basic} elements are the same  up to an orthogonal change of variable. 
This is not the case when using DeC due to the difference in lumped mass matrices.
%
 Similarly, the mass matrix used for \textit{cubature }elements is already diagonal, which makes the DeC procedure entirely equivalent to the RK scheme with Butcher  tableau corresponding to the quadrature weights of the DeC.
Only for SUPG  a  difference is observed due to the contributions to the mass matrix of the stabilization.}\\

\revMR{Concerning the plots, it is interesting to remark the appearance of four different structures which have an impact on the practical usefulness of some of the combinations.
	\begin{itemize}
\item The first kind of structures  are associated to schemes presenting V-shaped stability regions. We can observe these on  \cref{fig_disp_cohen_FD_SSPRK-cfl_vs_tau-SUPG,fig_disp_cohen_FD_SSPRK-cfl_vs_tau-CIP}, for 
 $p=1$. This shape requires a very careful choice of the stability parameter as  \revDT{ small perturbations of $\delta$} may lead, for a given CFL, to  an unstable behavior.   Generally, lowering the CFL increases somewhat the
 robustness allowing more flexibility in the choice of $\delta$.  
We highlight that this type of topology is common to   all the second order schemes, as well as to all DeC schemes with \textit{basic }and \textit{Bernstein} elements for degree $p\geq 2$.

\item Another  structure typically observed is an L-shaped stability region as in  \cref{fig_disp_cohen_FD_SSPRK-cfl_vs_tau-SUPG,fig_disp_cohen_FD_SSPRK-cfl_vs_tau-CIP} for $p=2,3$.
This shape is characterized by a CFL bound $\CFL\leq C_1$ and a one--sided bound on the stabilization coefficient $\delta\leq C_2 \CFL^{C_3}$, and it much more robust concerning the choice of  the stability parameter
as all values below a certain maximum are stable. Most of the schemes with $p\geq 2$, besides those listed in the first group, belong to this category.


\begin{figure}
	\centering
	\includegraphics[width=\textwidth]{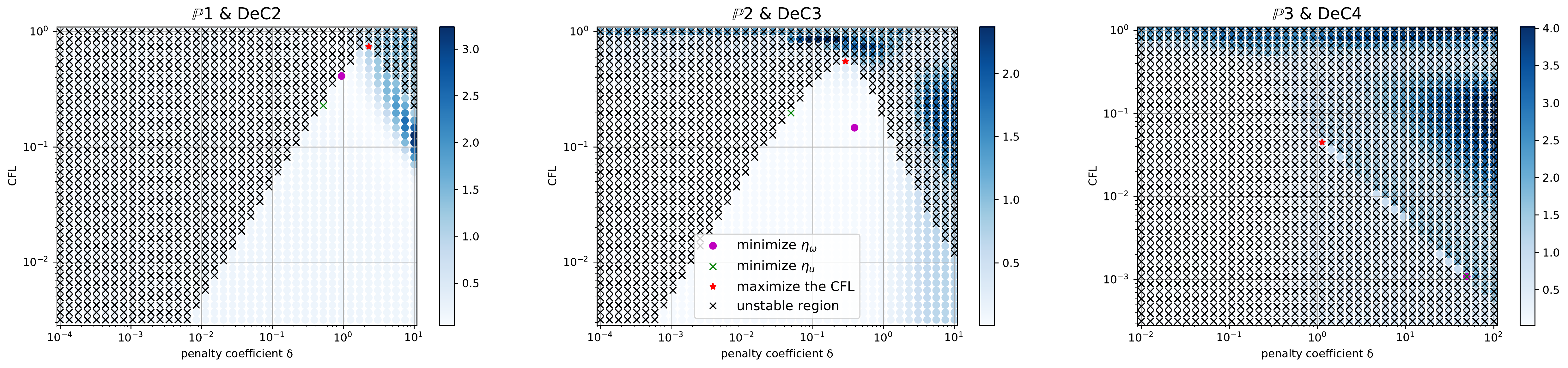}
	\caption{Computation of optimal parameters according to errors $\eta_\omega$ and $\eta_u$. $(\CFL,\delta)$ plot of $\eta_u$ (blue scale) and instability area (black crosses) for basic elements DeC scheme with LPS stabilization method. From left to right $\mathbb{P}_1$, $\mathbb{P}_2$, $\mathbb{P}_3$. The purple circle is the optimizer of $\eta_u$, the green cross is the optimizer of $\eta_\omega$, the red star is the maximum stable CFL. }
	\label{fig_disp_lagrange_FD_DeC-cfl_vs_tau-LPS}
\end{figure}

\item The third kind of structures involve  ``broom"- or ``box''-shaped stability domains. In the first case  we observe  
 two clear bounds $\delta\geq C_1\CFL^{C_2}$ and $\delta <C_3$ plus   a small stable stripe with higher $\CFL>(C_3/C_1)^{1/C_2}$ and $\delta > C_3$. This is for example visible 
 in  \cref{fig_disp_bezier_FD_DeC-cfl_vs_tau-SUPG}.  In the second case, see for example  \cref{fig_disp_cohen_FD_DeC-cfl_vs_tau-SUPG}, we also have  two bounds of the type $\CFL\geq C_1$ and $\delta <C_2$,
 with an additional   stable stripe outside these bounds. The problem with this type of methods is that the optimal parameters, {\it viz.} those involving the highest CFL, 
 are within a stripe which means that instability may be introduced by lowering the CFL\footremember{stripe}{These values do not allow to decrease the CFL}. For applications involving multiscale problems, or  variable mesh sizes this is clearly unacceptable in practice.
 Schemes showing this sort of  behaviors are   all the SUPG schemes with DeC time stepping, and with $p\geq 2$\revDT{, for which we indicate good values $(\CFL,\delta)$ in \cref{tab:dispersion_cfl_box}.}

\item Finally, the DeC scheme with \textit{basic} elements and $p=3$ shows essentially everywhere instability for CIP and LPS stabilization. The study finds some very thin oblique stripes of stability, but they are not wide enough to find stable regions.
See \cref{fig_disp_lagrange_FD_DeC-cfl_vs_tau-LPS} for an example.
\end{itemize}}

\begin{table}[h] 
	\small
	\begin{center} 
		\begin{tabular}{| c | c || c | c | c || c |c | c | }  
	     \hline 
	     \multicolumn{2}{|c||}{Element $\&$ }  & \multicolumn{3}{|c||}{No stabilization }  & \multicolumn{3}{|c|}{SUPG}  \\ \hline 
	     \multicolumn{2}{|c||}{ Time scheme }  & $p=1$ & $p=2$ & $p=3$   & $p=1$ & $p=2$ & $p=3$  \\ \hline \hline 
        \parbox[t]{2mm}{\multirow{3}{*}{\rotatebox[origin=c]{90}{Basic}}}              &  RK  &   /   &  0.389  &  0.389  &  0.624 (0.464)   &  0.492 (0.07)   &  0.389 (0.027)  \\ 
               &  SSPRK  &   /   &  0.492  &  0.389  &  0.889 (0.464)   &  0.554 (0.089)   &  0.438 (0.027)  \\ 
               &  DeC  &  /  &  /  &  /  &  1.701 (0.588)    &  0.492 (0.229)\footrecall{stripe}   &  0.492 (0.089)\footrecall{stripe}  \\
         \hline 
       \parbox[t]{2mm}{\multirow{3}{*}{\rotatebox[origin=c]{90}{\centering  Cub.}}}              &  RK  &   /   &  0.492  &  0.492  &  0.971 (0.767)   &  0.624 (0.13)   &  0.464 (0.064)  \\ 
               &  SSPRK  &   /   &  0.624  &  0.492  &  1.512 (0.642)   &  0.838 (0.13)   &  0.538 (0.064)  \\ 
               &  DeC  &   /   &  0.492  &  0.492  &  1.701 (0.398)   &  1.0 (0.081)\footrecall{stripe}   &  0.588 (0.041)\footrecall{stripe}  \\ 
         \hline 
       \parbox[t]{2mm}{\multirow{3}{*}{\rotatebox[origin=c]{90}{\centering  Bern.}}}              &  RK  &   /   &  0.389  &  0.389  &  0.624 (0.464)   &  0.492 (0.07)   &  0.389 (0.027)  \\ 
               &  SSPRK  &   /   &  0.492  &  0.389  &  0.889 (0.464)   &  0.554 (0.089)   &  0.438 (0.027)  \\ 
               &  DeC  &  /  &  /  &  /  &  1.701 (0.588)    &  1.0 (0.367)\footrecall{stripe}   &  0.702 (0.229)\footrecall{stripe}  \\  
         \hline 
        \hline 
        \end{tabular} 
    \end{center} 
 \begin{center} 
		\begin{tabular}{| c | c || c | c | c || c | c | c | }  
	     \hline 
	     \multicolumn{2}{|c||}{Element $\&$ }  & \multicolumn{3}{|c||}{LPS}  & \multicolumn{3}{|c|}{CIP}  \\ \hline 
	     \multicolumn{2}{|c||}{ Time scheme }  & $p=1$ & $p=2$ & $p=3$   & $p=1$ & $p=2$ & $p=3$  \\ \hline \hline 
        \parbox[t]{2mm}{\multirow{3}{*}{\rotatebox[origin=c]{90}{Basic}}}              &  RK  &  0.681 (0.767)   &  0.478 (0.077)   &  0.378 (0.032)   &  0.838 (0.094)   &  0.538 (5.54e-03)   &  0.4 (8.38e-04)   \\ 
               &  SSPRK  &  1.093 (0.767)   &  0.605 (0.109)   &  0.425 (0.038)   &  1.125 (0.119)   &  0.624 (7.02e-03)   &  0.464 (6.61e-04)   \\ 
               &  DeC  &  0.744 (2.29)   &  0.554 (0.289)   &  /   &  0.838 (0.289)   &  0.588 (0.02)   &  /   \\ 
         \hline 
       \parbox[t]{2mm}{\multirow{3}{*}{\rotatebox[origin=c]{90}{\centering  Cub.}}}              &  RK  &  1.093 (0.702)   &  0.681 (0.143)   &  0.538 (0.049)   &  0.971 (0.191)   &  0.723 (0.011)   &  0.538 (1.84e-03)   \\ 
               &  SSPRK  &  1.557 (1.0)   &  0.863 (0.17)   &  0.605 (0.049)   &  1.512 (0.242)   &  0.838 (0.014)   &  0.538 (3.93e-03)   \\ 
               &  DeC  &  1.093 (0.702)   &  0.681 (0.143)   &  0.538 (0.049)   &  0.971 (0.191)   &  0.723 (0.011)   &  0.538 (1.84e-03)   \\ 
         \hline 
       \parbox[t]{2mm}{\multirow{3}{*}{\rotatebox[origin=c]{90}{\centering  Bern.}}}              &  RK  &  0.681 (0.767)   &  0.478 (0.077)   &  0.378 (0.032)   &  0.838 (0.094)   &  0.538 (5.54e-03)   &  0.4 (8.38e-04)   \\ 
               &  SSPRK  &  1.093 (0.767)   &  0.605 (0.109)   &  0.425 (0.038)   &  1.125 (0.119)   &  0.624 (7.02e-03)   &  0.464 (6.61e-04)   \\ 
               &  DeC  &  0.744 (2.29)   &  0.052 (0.215)   &  0.109 (0.215)   &  0.838 (0.289)   &  0.059 (0.016)   &  0.119 (7.02e-03)   \\ 
         \hline 
        \end{tabular} 
    \end{center} 
     \caption{Optimized CFL and penalty coefficient $\delta $ in parenthesis, only maximizing CFL } \label{tab:dispersion_cfl-best}
\end{table}%

\begin{table}[h] 
	\small
	\begin{center} 
		\begin{tabular}{| c | c || c | c | c || c |c | c | }  
	     \hline 
	     \multicolumn{2}{|c||}{Element $\&$ }  & \multicolumn{3}{|c||}{No stabilization }  & \multicolumn{3}{|c|}{SUPG}  \\ \hline 
	     \multicolumn{2}{|c||}{ Time scheme }  & $p=1$ & $p=2$ & $p=3$   & $p=1$ & $p=2$ & $p=3$  \\ \hline \hline 
        \parbox[t]{2mm}{\multirow{3}{*}{\rotatebox[origin=c]{90}{Basic}}}              &  RK  &   /   &  0.151  &  0.191  &  0.389 (0.089)   &  0.17 (2.57e-03)   &  0.215 (8.38e-03)  \\ 
               &  SSPRK  &   /   &  0.191  &  0.242  &  0.492 (0.089)   &  0.215 (2.57e-03)   &  0.273 (5.22e-03)  \\ 
               &  DeC  &  /  &  /  &  /  &  0.702 (0.588)   &  0.143 (0.022)   &  0.024 (0.013)  \\ 
         \hline 
       \parbox[t]{2mm}{\multirow{3}{*}{\rotatebox[origin=c]{90}{\centering  Cub.}}}              &  RK  &   /   &  0.492  &  0.242  &  0.971 (0.538)   &  0.624 (0.045)   &  0.222 (0.019)  \\ 
               &  SSPRK  &   /   &  0.624  &  0.307  &  1.304 (0.378)   &  0.723 (0.038)   &  0.298 (3.78e-03)  \\ 
               &  DeC  &   /   &  0.492  &  0.242  &  0.346 (0.642)    &  0.702 (0.026)   &  0.203 (0.041)   \\ 
         \hline 
       \parbox[t]{2mm}{\multirow{3}{*}{\rotatebox[origin=c]{90}{\centering  Bern.}}}              &  RK  &   /   &  0.151  &  0.191  &  0.389 (0.089)   &  0.17 (2.57e-03)   &  0.215 (8.38e-03)  \\ 
               &  SSPRK  &   /   &  0.191  &  0.242  &  0.492 (0.089)   &  0.215 (2.57e-03)   &  0.273 (5.22e-03)  \\ 
               &  DeC  &  /  &  /  &  /  &  0.702 (0.588)    &  0.346 (0.367)\footrecall{stripe}    &  0.588 (0.289)\footrecall{stripe}   \\ 
         \hline 
         \hline 
        \hline 
        \end{tabular} 
    \end{center} 
 \begin{center} 
		\begin{tabular}{| c | c || c | c | c || c | c | c | }  
	     \hline 
	     \multicolumn{2}{|c||}{Element $\&$ }  & \multicolumn{3}{|c||}{LPS}  & \multicolumn{3}{|c|}{CIP}  \\ \hline 
	     \multicolumn{2}{|c||}{ Time scheme }  & $p=1$ & $p=2$ & $p=3$   & $p=1$ & $p=2$ & $p=3$  \\ \hline \hline 
        \parbox[t]{2mm}{\multirow{3}{*}{\rotatebox[origin=c]{90}{Basic}}}              &  RK  &  0.335 (0.077)   &  0.165 (3.78e-03)   &  0.209 (0.013)   &  0.4 (0.011)   &  0.165 (1.60e-04)   &  0.222 (2.03e-04)   \\ 
               &  SSPRK  &  0.478 (0.077)   &  0.209 (3.78e-03)   &  0.265 (9.15e-03)   &  0.624 (0.011)   &  0.191 (2.03e-04)   &  0.257 (3.26e-04)   \\ 
               &  DeC  &  0.229 (0.522)   &  0.197 (0.049)   &  /   &  0.346 (0.077)   &  0.203 (2.42e-03)   &   /   \\ 
         \hline 
       \parbox[t]{2mm}{\multirow{3}{*}{\rotatebox[origin=c]{90}{\centering  Cub.}}}              &  RK  &  0.863 (0.492)   &  0.605 (0.041)   &  0.235 (0.012)   &  0.971 (0.119)   &  0.624 (3.46e-03)   &  0.257 (1.13e-04)   \\ 
               &  SSPRK  &  1.23 (0.412)   &  0.767 (0.041)   &  0.298 (4.12e-03)   &  1.304 (0.094)   &  0.723 (3.46e-03)   &  0.298 (1.45e-04)   \\ 
               &  DeC  &  0.863 (0.492)   &  0.605 (0.041)   &  0.235 (0.012)   &  0.971 (0.119)   &  0.624 (3.46e-03)   &  0.257 (1.13e-04)   \\ 
         \hline 
       \parbox[t]{2mm}{\multirow{3}{*}{\rotatebox[origin=c]{90}{\centering  Bern.}}}              &  RK  &  0.335 (0.077)   &  0.165 (3.78e-03)   &  0.209 (0.013)   &  0.4 (0.011)   &  0.165 (1.60e-04)   &  0.222 (2.03e-04)   \\ 
               &  SSPRK  &  0.478 (0.077)   &  0.209 (3.78e-03)   &  0.265 (9.15e-03)   &  0.624 (0.011)   &  0.191 (2.03e-04)   &  0.257 (3.26e-04)   \\ 
               &  DeC  &  0.229 (0.522)   &  0.052 (0.215)   &  0.109 (0.215)   &  0.346 (0.077)   &  0.059 (0.016)   &  0.119 (7.02e-03)   \\ 
         \hline 
        \end{tabular} 
    \end{center} 
     \caption{Optimized CFL and penalty coefficient $\delta $ in parenthesis, minimizing  $\eta_u$ } \label{tab:dispersion_cfl-RES}
\end{table}%

\begin{table}[h] 
\small	\begin{center} 
		\begin{tabular}{| c | c || c | c | c || c |c | c | }  
	     \hline 
	     \multicolumn{2}{|c||}{Element $\&$ }  & \multicolumn{3}{|c||}{No stabilization }  & \multicolumn{3}{|c|}{SUPG}  \\ \hline 
	     \multicolumn{2}{|c||}{ Time scheme }  & $p=1$ & $p=2$ & $p=3$   & $p=1$ & $p=2$ & $p=3$  \\ \hline \hline 
        \parbox[t]{2mm}{\multirow{3}{*}{\rotatebox[origin=c]{90}{Basic}}}              &  RK  &   /   &  0.191  &  0.307  &  0.059 (0.289)   &  0.191 (0.027)   &  0.307 (0.044)  \\ 
               &  SSPRK  &   /   &  0.242  &  0.307  &  0.084 (0.289)   &  0.242 (0.027)   &  0.346 (0.035)  \\ 
               &  DeC  &   /  & /  &  /  &  0.412 (0.367)  &  0.242 (0.089)\footrecall{stripe}    &  0.017 (0.113)\footrecall{stripe}   \\ 
         \hline 
       \parbox[t]{2mm}{\multirow{3}{*}{\rotatebox[origin=c]{90}{\centering  Cub.}}}              &  RK  &   /   &  0.492  &  0.389  &  0.538 (0.767)   &  0.298 (0.316)   &  0.165 (0.156)  \\ 
               &  SSPRK  &   /   &  0.624  &  0.492  &  0.624 (0.915)   &  0.4 (0.316)   &  0.257 (0.186)  \\ 
               &  DeC  &   /   &  0.492  &  0.389  &  0.346 (0.642)     &  0.346 (0.179)\footrecall{stripe}    &  0.1 (0.09)\footrecall{stripe}  \\
         \hline 
       \parbox[t]{2mm}{\multirow{3}{*}{\rotatebox[origin=c]{90}{\centering  Bern.}}}              &  RK  &   /   &  0.191  &  0.307  &  0.059 (0.289)   &  0.191 (0.027)   &  0.307 (0.044)  \\ 
               &  SSPRK  &   /   &  0.242  &  0.307  &  0.084 (0.289)   &  0.242 (0.027)   &  0.346 (0.035)  \\ 
               &  DeC  &  /  &  /  &  /  &  0.412 (0.367)   &  0.289 (0.289)\footrecall{stripe}    &  0.203 (0.289)\footrecall{stripe}  \\ 
         \hline 
        \hline 
        \end{tabular} 
    \end{center} 
 \begin{center} 
		\begin{tabular}{| c | c || c | c | c || c | c | c | }  
	     \hline 
	     \multicolumn{2}{|c||}{Element $\&$ }  & \multicolumn{3}{|c||}{LPS}  & \multicolumn{3}{|c|}{CIP}  \\ \hline 
	     \multicolumn{2}{|c||}{ Time scheme }  & $p=1$ & $p=2$ & $p=3$   & $p=1$ & $p=2$ & $p=3$  \\ \hline \hline 
        \parbox[t]{2mm}{\multirow{3}{*}{\rotatebox[origin=c]{90}{Basic}}}              &  RK  &  0.478 (0.186)   &  0.13 (0.265)   &  0.116 (0.13)   &  0.464 (0.037)   &  0.123 (0.011)   &  0.165 (3.46e-03)   \\ 
               &  SSPRK  &  0.605 (0.378)   &  0.165 (0.265)   &  0.335 (0.026)   &  0.624 (0.046)   &  0.143 (0.014)   &  0.346 (5.22e-04)   \\ 
               &  DeC  &  0.412 (0.943)   &  0.147 (0.389)   &  /   &  0.588 (0.13)   &  0.143 (0.016)   &   /   \\ 
         \hline 
       \parbox[t]{2mm}{\multirow{3}{*}{\rotatebox[origin=c]{90}{\centering  Cub.}}}              &  RK  &  0.971 (0.492)   &  0.538 (0.119)   &  0.425 (0.024)   &  0.971 (0.119)   &  0.538 (0.011)   &  0.4 (4.00e-04)   \\ 
               &  SSPRK  &  1.23 (0.492)   &  0.681 (0.119)   &  0.478 (1.43e-03)   &  1.304 (0.119)   &  0.723 (7.02e-03)   &  0.257 (1.11e-03)   \\ 
               &  DeC  &  0.971 (0.492)   &  0.538 (0.119)   &  0.425 (0.024)   &  0.971 (0.119)   &  0.538 (0.011)   &  0.4 (4.00e-04)   \\ 
         \hline 
       \parbox[t]{2mm}{\multirow{3}{*}{\rotatebox[origin=c]{90}{\centering  Bern.}}}              &  RK  &  0.478 (0.186)   &  0.13 (0.265)   &  0.116 (0.13)   &  0.464 (0.037)   &  0.123 (0.011)   &  0.165 (3.46e-03)   \\ 
               &  SSPRK  &  0.605 (0.378)   &  0.165 (0.265)   &  0.335 (0.026)   &  0.624 (0.046)   &  0.143 (0.014)   &  0.346 (5.22e-04)   \\ 
               &  DeC  &  0.412 (0.943)   &  0.052 (0.215)   &  0.109 (0.215)   &  0.588 (0.13)   &  0.059 (0.016)   &  0.119 (7.02e-03)   \\ 
         \hline 
        \end{tabular} 
    \end{center} 
     \caption{Optimized CFL and penalty coefficient $\delta $ in parenthesis, minimizing $\eta_\omega$ } \label{tab:dispersion_cfl-REW}
\end{table}%

\begin{table}[h] 
\small	\begin{center} 
		\begin{tabular}{| c || c |c  | }  
	     \hline 
	     DeC  & \multicolumn{2}{|c|}{SUPG} \\ \hline 
	      Element &  $p=2$ & $p=3$ \\ \hline \hline 
        Basic               &  0.08 (0.025)   &  0.059 (0.035)  \\ 
        Cubature     & 0.346 (0.025) & 0.242 (2.22 e-03) \\
        Bernstein  & 0.03 (0.025) & 0.1 (0.1) \\
         \hline 
        \end{tabular} 
    \end{center} 
     \caption{Optimized CFL and penalty coefficient $\delta $ in parenthesis, stable for all smaller CFLs} \label{tab:dispersion_cfl_box}
\end{table}%

\subsection{Dispersion and damping}

\begin{figure}
     \centering
     \includegraphics[width=\textwidth]{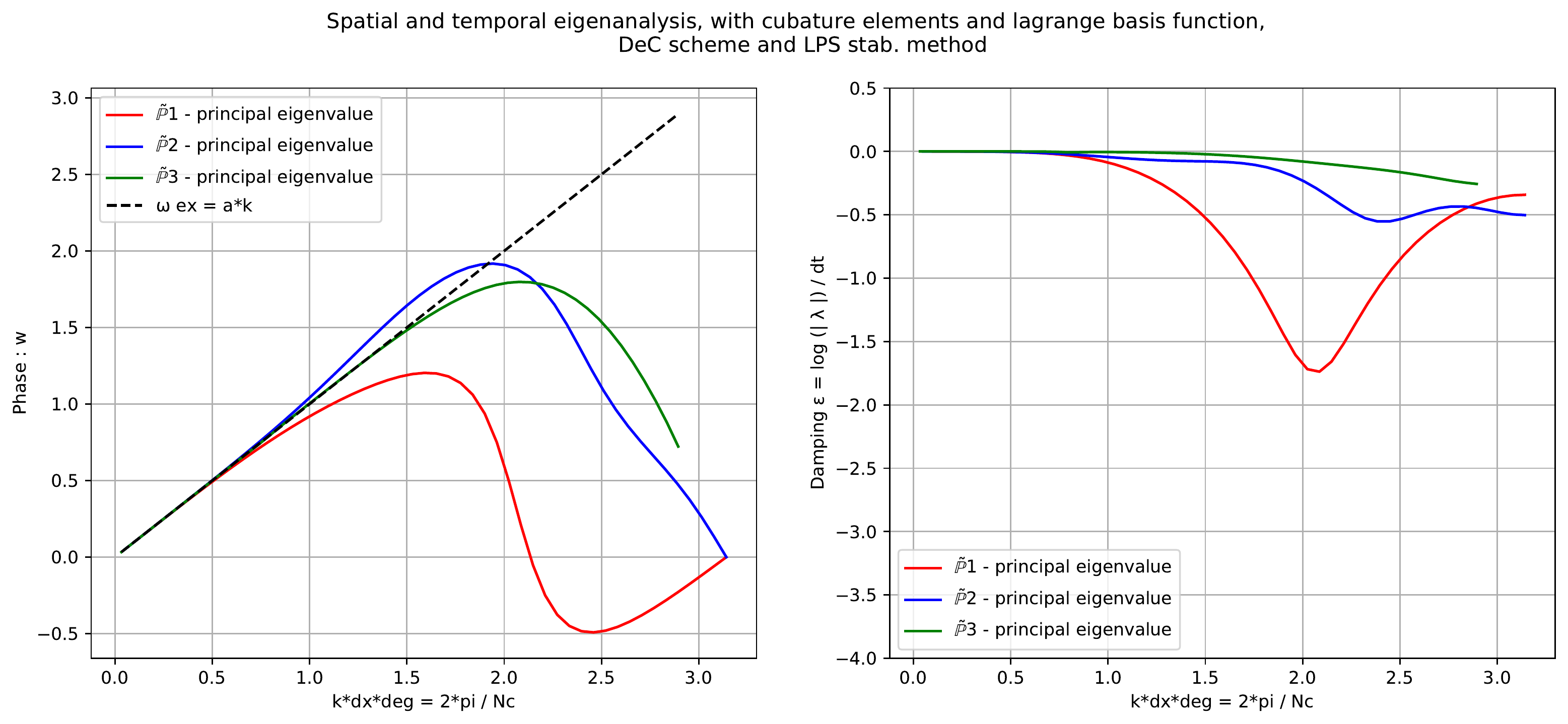}
     \caption{Comparison of dispersion in the fully discrete case, using \revsix{coefficients from \ref{tab:dispersion_cfl-RES}}, \textit{cubature} elements, DeC scheme and LPS stabilization method. $\mathbb{P}_1$ elements in red, $\mathbb{P}_2$ elements in blue and $\mathbb{P}_3$ elements in green. The phase $\omega$ of the principal eigenvalues is on the left and the damping $\epsilon_i$ on the right}
     \label{fig_disp_cohen_FD_DeC-LPS}
\end{figure}

\begin{figure}
     \centering
     \includegraphics[width=\textwidth]{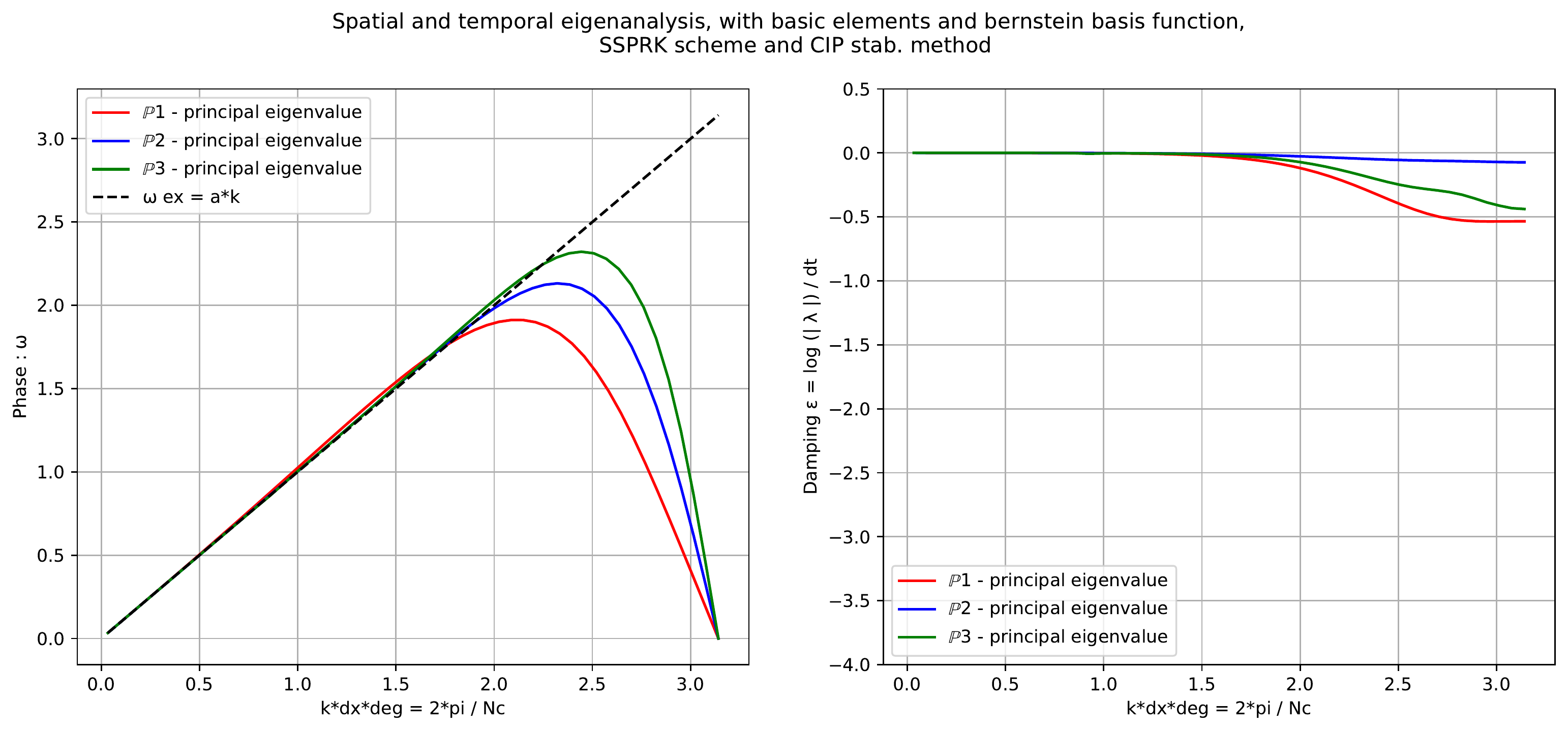}
     \caption{Comparison of dispersion in the fully discrete case, using \revsix{coefficients from \ref{tab:dispersion_cfl-RES}}, \textit{Bernstein} elements, SSPRK scheme and CIP stabilization method. $\mathbb{B}_1$ elements in red, $\mathbb{B}_2$ elements in blue and $\mathbb{B}_3$ elements in green. The phase $\omega$ of the principal eigenvalues is on the left and the damping $\epsilon_i$ on the right.}
     \label{fig_disp_bezier_FD_SSPRK-CIP}
\end{figure}

In \cref{fig_disp_cohen_FD_DeC-LPS,fig_disp_bezier_FD_SSPRK-CIP} are represented the phase and the damping of the principal eigenvalue depending on $\theta=k \Delta x=\frac{2 \pi}{N_x}$ for few schemes\revDT{ (\textit{cubature} DeC LPS and \textit{Bernstein} SSPRK CIP),} using the best parameters $(\text{CFL},\delta)$ found in the previous analysis with the optimization of $\eta_u$. 
As before, we notice that the mode for $p=1$ is particularly dispersive.
Nevertheless, the frequencies on which the scheme is dispersive are also much damped as we see in the right plots. 
For higher order methods, the phase $\omega$ of the principal mode is closer to the exact phase $\omega_{ex}=ak$ in the left figures. 
We observe that the principal mode of higher order methods is much more precise in terms of dispersion than the first order one, but also less damped in the low frequency area $\theta \geq \frac{2\pi}{3}$. 

For completeness, a comparison of damping and phase coefficients for DeC and SSPRK for all the stabilization techniques and elements can be found in \cref{app:fourier_full_disc}. 
There we used the (CFL,$\delta$) coefficients found by minimizing $\eta_u$ in \cref{tab:dispersion_cfl-RES}, and we try also to compare the obtained results. Nevertheless, we must remark that the different CFLs used for different schemes do not allow a direct comparison.\\
 


The different strategies lead to different values of best CFL and $\delta$. In general, the most reliable is the one that optimizes $\eta_u$.  Looking at \cref{tab:dispersion_cfl-RES}, we can compare the different elements, stabilization terms and time integration techniques and obtain some conclusions.
\begin{itemize}
	\item \revMR{In general  SSPRK time integration methods allow to use higher CFL}  with respect to both classical RK methods and DeC.
	\item With \textit{cubature} elements we can use \revDT{larger} CFLs conditions than with \textit{basic} and \textit{Bernstein} elements.
	\item \revMR{Concerning efficiency,  we do not observe any impact of the choice of the stabilization approach on the magnitude of the allowed CFL. Other factors are much more relevant in this respect.
	For example, for SUPG we need to stress the advantage of using DeC w.r.t. the possibility of avoiding the 
	  inversion of the non-diagonal mass matrix required by the full consistency of the method. 
	    For   CIP the larger stencil and  non-local data structure gives a small overhead, and, for LPS, the gradient projection favors clearly  \textit{cubature}  elements for which this phase requires no matrix inversion.}
	\item Some combinations produce very unstable schemes.
	As  remarked also before, DeC with high order \textit{basic }elements may have problems in the mass lumping, and we can see an example with the LPS and CIP stabilization. 
	\item DeC with SUPG stabilization leads to stability regions that are not comprehending all the CFLs smaller than the one inside the region, for a fixed $\delta$. This is very dangerous, for instance when doing mesh adaptation algorithms, hence, we marked with an asterisk in \cref{tab:dispersion_cfl-RES,tab:dispersion_cfl-REW,tab:dispersion_cfl-best} such schemes and we put in \cref{tab:dispersion_cfl_box} reliable values of (CFL,$\delta)$.
\end{itemize}

\section{A note on nonlinear stability}\label{sec:nonlinear}
The stability analysis performed before holds only for linear problems.
For nonlinear ones the original ansatz of supposing that the solutions can be decomposed orthogonally into waves that propagate at constant speed does not hold anymore.
Nevertheless, the stabilization methods presented \revMR{also introduces  some  nonlinear stabilization. To show it we will briefly consider
their potential for dissipating entropy.}
 In order to test so, we neglect the time discretization, the used elements and the quadrature and the discrete differentiation formulae.

Consider any convex smooth entropy $\rho(u)$, i.e., $\rho_{uu}(u)>0$, the respective entropy variables $\nu:=\rho_u(u)$ and the entropy flux $g(u)$ such that $\rho_u f_u = g_u$. In the following discussion, we consider the entropy variable $\nu_h=\rho_u(u)_h$ to be in the finite element space, while $u_h$ will be defined as the projection onto the finite element space of the uniquely defined function $\nu\to u=u(\nu)$, as proposed in \cite{abgrallEntropyConservative}.

When substituting $v_h=\nu_h$, the Galerkin discretization of the conservation law becomes
\begin{equation}
	\sum_K \int_{K} \nu_h \left( \partial_t u_h + \partial_x f(u_h) \right) dx = \sum_K \int_{K} \partial_t \rho_h + \partial_x g_h dx = \int_{\Omega} \partial_t \rho_h + \left[ g_h \right]_{\partial K},
\end{equation} 
which, according to the boundary conditions, gives us a measure of the variation of the entropy.

The CIP stabilization must be slightly modified for nonlinear equations with nontrivial entropies, so that it reads
\begin{equation}\label{eq:CIPforEntropy}
	s(v,u):=\sum_{K,{\sf f}\in K} \!\int_{\sf f} [\partial_x v^T]\rho_{uu}(u)^{-1}[\partial_x \nu(u)] d\Gamma,
\end{equation}
where the inverse of the hessian of the entropy must be added for unit of measure reasons and it is positive definite and invertible. So that when we substitute $v=\nu_h$ in the stabilization term, we obtain
\begin{equation}
	\begin{split}
	s(\nu,u_h)= &\!\sum_{K,{\sf f}\in K} \!\int_{\sf f} \underbrace{[\partial_x \nu_h^T]\rho_{uu}(u_h)^{-1}[\partial_x \nu_h]}_{>0} d\Gamma. 
	\end{split}
\end{equation}
It would guarantee a decrease in the discrete total entropy. Moreover, this formulation coincide with \eqref{eq_cip0} when we are dealing with the energy as entropy.

For the LPS we modify, similarly the formulation \eqref{eq_lps0} into
\begin{equation}
	\begin{cases}
		s(v,u):= \sum_K \tau_K\int_K \partial_x v^T \rho_{uu}(u)^{-1} (\partial_x \nu(u)-w) dx, \text{  with}\\
		\int_K z^T (w-\partial_x \nu(u)), \quad \forall z \in V_h
	\end{cases}
\end{equation}
As in the linear case, we can  take $\tau_K=\tau$, and  test with  $v_h=\nu_h$ in the stabilization term and we substitute $z=  \tau \rho_{uu}(u)^{-1,T} w  $ in the previous equation and we sum this 0 contribution to the stabilization term, we obtain
\begin{equation}
	\begin{split}
		s(\nu_h,u_h)= &\sum_K \tau\int_K \partial_x\nu_h^T \rho_{uu}(u_h)^{-1}  (\partial_x \nu_h - w_h)  + \rho_{uu}(u_h)  w^T_h \rho_{uu}(u_h)^{-1}(w_h-\partial_x \nu_h) dx =\\
		&\sum_K \tau\int_K (\partial_x\nu_h - w_h)^T\rho_{uu}(u_h)^{-1}(\partial_x \nu_h-w_h)  dx \geq 0.
	\end{split}
\end{equation}
As for the CIP we can say that the LPS stabilization reduces entropy. Anyway, this analysis does not guarantee that the fully discrete method will be entropy stable, as all the other discretizations (time, quadrature, differentiation and interpolation) are not taken into consideration.

\revMR{For the SUPG stabilization, as the linear analysis of \cref{SUPG}  shows, the spatial and temporal derivatives need to be properly combined.
This can be done easily for space-time discretizations (see  e.g. in \cite{Barth98}), context in which SUPG and  least squares stabilization  coincide. 
In simple cases   with constant convexity entropy, namely the energy, one can bound other types of energy norm in time, but not the entropy itself.  
For explicit methods, and general convex entropies, the non-symmetric nature of the method requires ad-hoc analysis which we leave out of this paper. }
More elaborated analysis are possible with other types of stabilization, as the ones proposed in \cite{abgrallEntropyConservative,kuzmin2020Algebraic,guermond2011Entropy}, and they will be the object of future research.\\

In the next sections, we perform also some nonlinear tests, where we use the coefficients we found in the stability analysis for the linear case, in order to understand if this 
information is also relevant for   nonlinear problems.

\section{Numerical Simulations}\label{sec:simulation}

\revMR{We perform numerical tests to check  the validity of our theoretical findings. 
We will use elements of degree $p$, with $p$ up to 3,  with time integration schemes of the corresponding order to ensure an overall error of $\mathcal{O}(\Delta x ^{p+1})$, under the CFL conditions presented earlier \revDT{in \cref{tab:dispersion_cfl-RES}. The integral formulae are performed with high order quadrature rules, for \textit{cubature} elements they are associated with the definition points of the elements themselves, for \textit{basic} and \textit{Bernstein} we use Gauss--Legendre quadrature formulae with $p+1$ points per cell.  } } 
\subsection{Linear advection equation}
We start with the one dimensional initial value problem for the linear advection  equation \eqref{eq_disp_1} on the domain $\Omega = [0,2]$ using periodic boundary conditions:
\begin{equation}
\begin{cases}
\partial_t u (x,t) + a \partial_x u (x,t)  = 0  \qquad & \quad (x,t) \in \Omega\times [0,5] , \quad a \in \mathbb{R}, \label{lin_adv_1D} \\
u (x,0) = u_0(x), & \\
u (0,t) = u(2,t), & \quad t\in [0,5],
\end{cases}
\end{equation}
where $u_0(x) = 0.1 \sin(\pi x) $. Clearly the exact solution is $u_{ex}(x,t)=u_0(x-at)$ for all $x\in\Omega$.
We discretize the mesh with uniform intervals of length $\Delta x$.
In particular, we will use different discretization scales to test the convergence: $\Delta x_1 = \{ 0.05, 0.025, 0.0125, 0.00625 \}$ for $\mathbb{P}_1$ elements, $\Delta x_2 = 2\Delta x_1 $ for $\mathbb{P}_2$ elements and $\Delta x_3 = 3 \Delta x_1 $ for $\mathbb{P}_3$ elements. \revMR{This allows to guarantee the  use ot the same  number of degrees of freedom for   different $p$.}
We will compare the errors obtained with SSPRK and DeC time integration method, with all the stabilization methods (SUPG, LPS and CIP) and with \textit{basic}, \textit{cubature} and \textit{Bernstein} elements.

\begin{figure}
	\begin{minipage}{0.48\textwidth}
		\centering
		\includegraphics[width=0.99\textwidth]{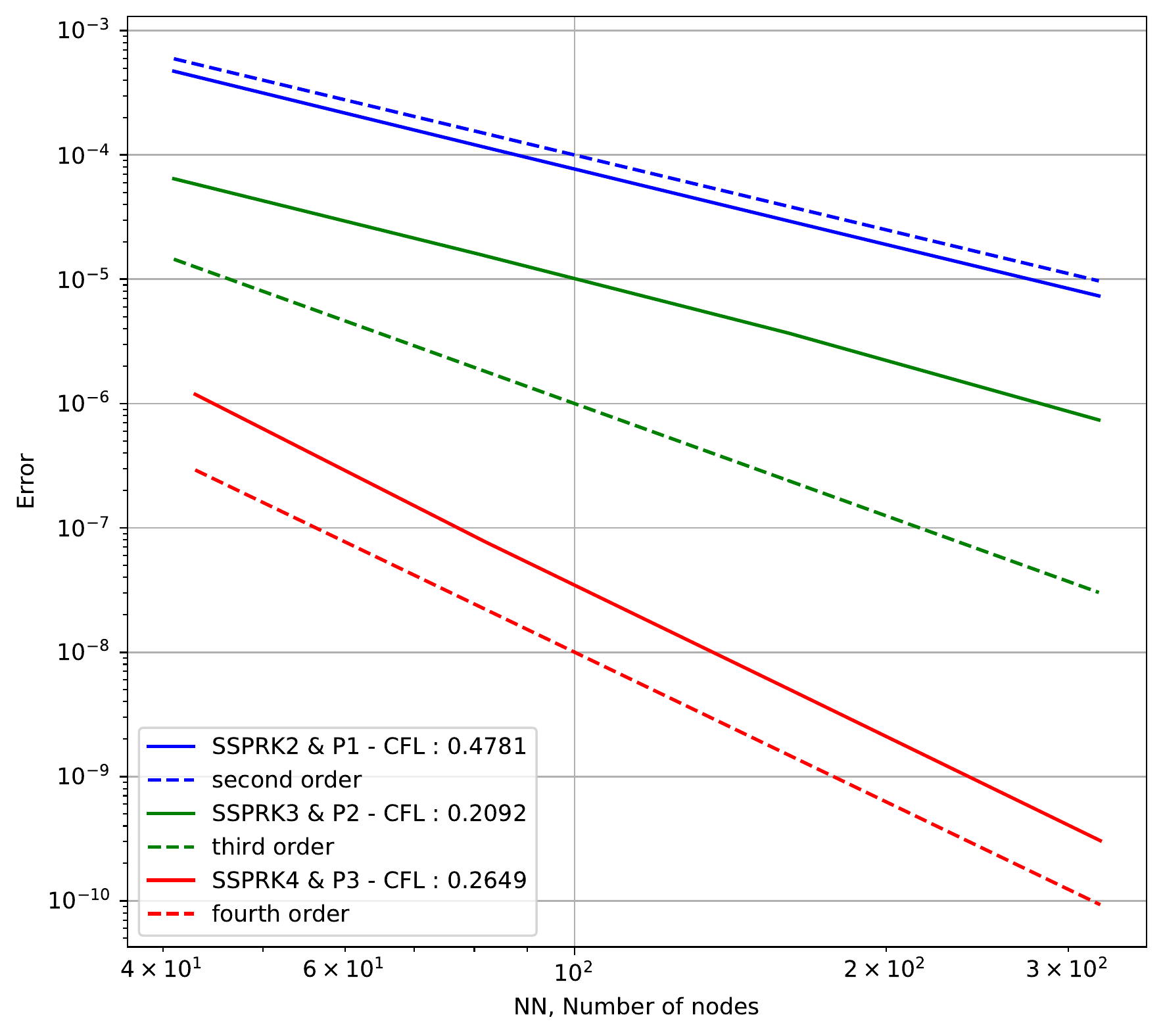}
		\caption{Error decay for linear advection with \textit{basic} elements, LPS stabilization and SSPRK. $\mathbb{P}_1,\,\mathbb{P}_2$ and $\mathbb{P}_3$ elements are, respectively, in blue green and red.}
		\label{fig_err_1d_EDP-SSPRK_lag-LPS}
	\end{minipage}\hfill
	\begin{minipage}{0.48\textwidth}
		\begin{center}
			\includegraphics[width=0.99\textwidth]{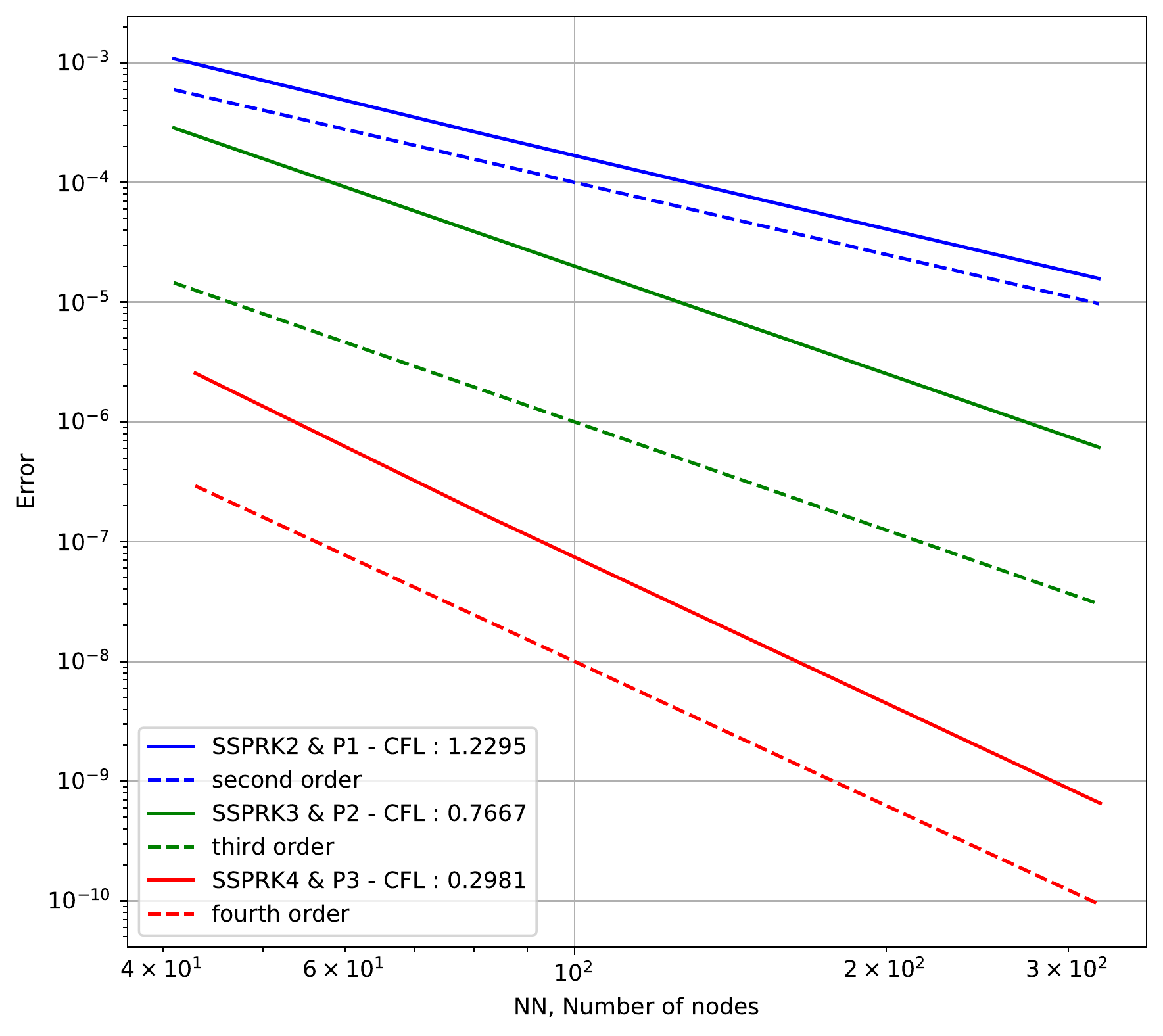}
		\end{center}
		\caption{Error decay for linear advection with \textit{cubature} elements, LPS stabilization and SSPRK. $\mathbb{P}_1,\,\mathbb{P}_2$ and $\mathbb{P}_3$ elements are, respectively, in blue green and red.}
		\label{fig_err_1d_EDP-SSPRK_cohen-LPS}
	\end{minipage}
\end{figure}

\revMR{A representative result is provided as an example in  \cref{fig_err_1d_EDP-SSPRK_lag-LPS,fig_err_1d_EDP-SSPRK_cohen-LPS}: it shows a } comparison between \textit{cubature }and \textit{basic} elements with LPS stabilization and SSPRK time integration. As we can see, the two schemes have very similar error behavior, but the \textit{basic} elements require stricter CFL conditions, see \cref{tab:dispersion_cfl-RES}, and have larger computational costs because of the full mass matrix.
A summary table with the order of accuracy reached by each simulations in \cref{tab:conv_order_LinearAdvection-RES}. The plots and all the errors are available at the repository \cite{TorloMichel2020git}.

\begin{table}[H] 
\small  
 \begin{center} 
		\begin{tabular}{| c | c || c | c | c || c | c | c || c | c | c || c | c | c | }  
	     \hline 
	     \multicolumn{2}{|c||}{Element $\&$ }  & \multicolumn{3}{|c||}{No stabilization }  & \multicolumn{3}{|c||}{SUPG}  & \multicolumn{3}{|c||}{LPS}  & \multicolumn{3}{|c|}{CIP}  \\ \hline 
	     \multicolumn{2}{|c||}{ Time scheme }  & $\mathbb{P}_1$ & $\mathbb{P}_2$ & $\mathbb{P}_3$   & $\mathbb{P}_1$ & $\mathbb{P}_2$ & $\mathbb{P}_3$   & $\mathbb{P}_1$ & $\mathbb{P}_2$ & $\mathbb{P}_3$   & $\mathbb{P}_1$ & $\mathbb{P}_2$ & $\mathbb{P}_3$  \\ \hline \hline 
       \parbox[t]{2mm}{\multirow{2}{*}{\rotatebox[origin=c]{90}{\centering  Cub.}}}              &  SSPRK &  /  & 1.98 & 3.98 & 2.04 & 2.93 & 3.98 & 2.03 & 2.95 & 3.98 & 2.05 & 2.94 & 3.98 \\ 
               &  DeC &  /  & 1.98 & 3.98 & 2.0 & 2.88 & 3.97 & 2.03 & 2.95 & 3.98 & 2.12 & 2.96 & 3.98 \\ 
         \hline 
        \parbox[t]{2mm}{\multirow{2}{*}{\rotatebox[origin=c]{90}{Basic}}}              &  SSPRK &  /  & 3.84 & 3.97 & 2.0 & 2.08 & 3.98 & 2.0 & 2.14 & 3.98 & 2.0 & 2.07 & 3.97 \\ 
               &  DeC &  /  &  /  &  /  & 2.02 & 2.72 & 2.05 & 1.95 & 2.93 &  /  & 1.98 & 2.82 &  /  \\ 
         \hline 
       \parbox[t]{2mm}{\multirow{2}{*}{\rotatebox[origin=c]{90}{\centering  Bern.}}}              &  SSPRK &  /  & 3.84 & 3.97 & 2.0 & 2.08 & 3.98 & 2.0 & 2.14 & 3.98 & 2.0 & 2.07 & 3.97 \\ 
               &  DeC &  /  &  /  &  /  &  /  &  /  &  /  & 1.98 & 3.05 & 2.04 & 1.98 & 3.0 & 2.0 \\
         \hline 
        \end{tabular} 
    \end{center} 
     \caption{Summary table of convergence orders, using coefficients obtained by minimizing $\eta_u$ in \cref{tab:dispersion_cfl-RES} } \label{tab:conv_order_LinearAdvection-RES}
\end{table}%

\begin{figure}[h!]
	\begin{center}
		\textit{Cubature} elements\\\vspace{2mm}
		\includegraphics[height=0.21\textheight,trim={0 0 62mm 0}, clip]{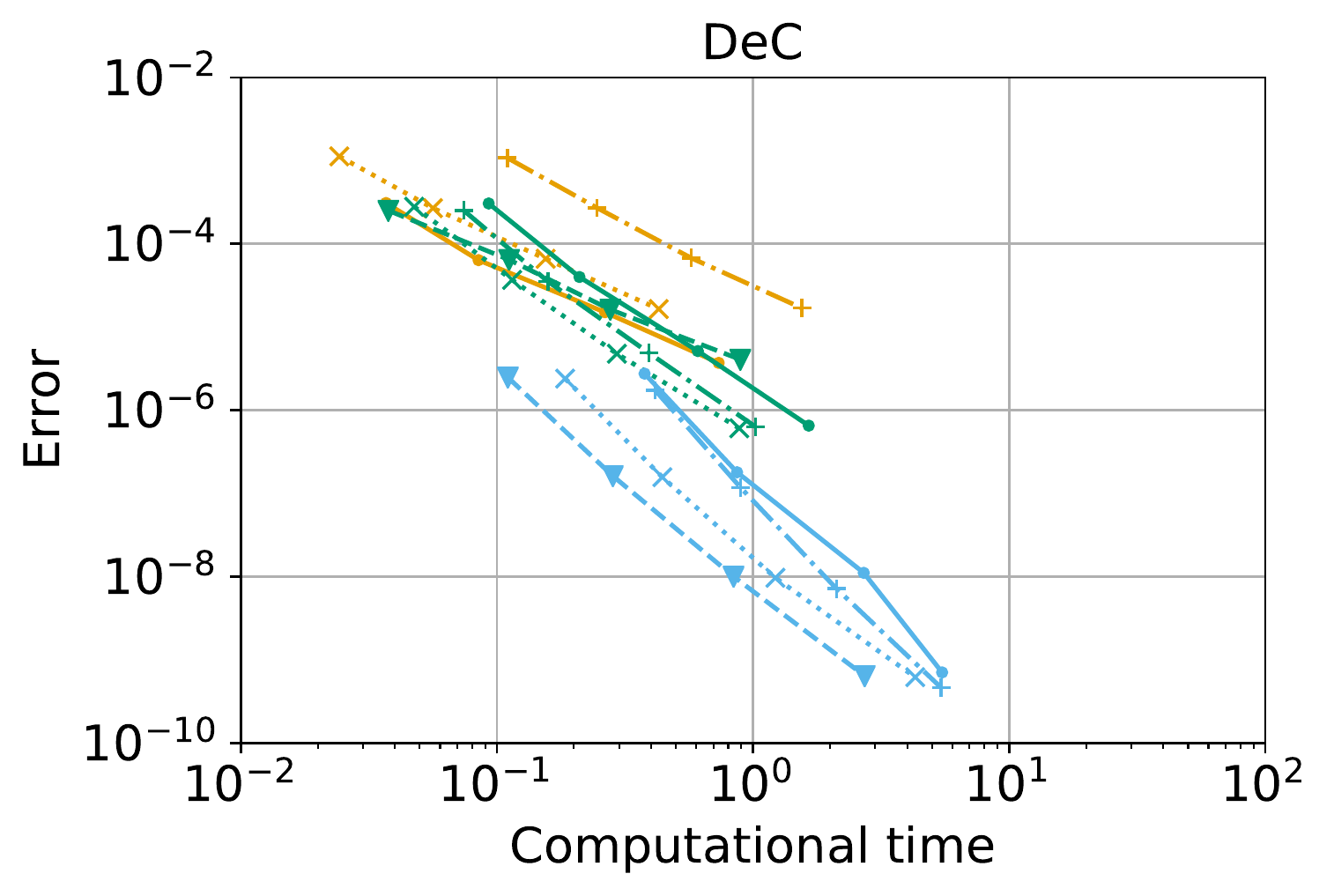} \hfill
		\includegraphics[height=0.21\textheight,trim={10mm 0 0 0}, clip]{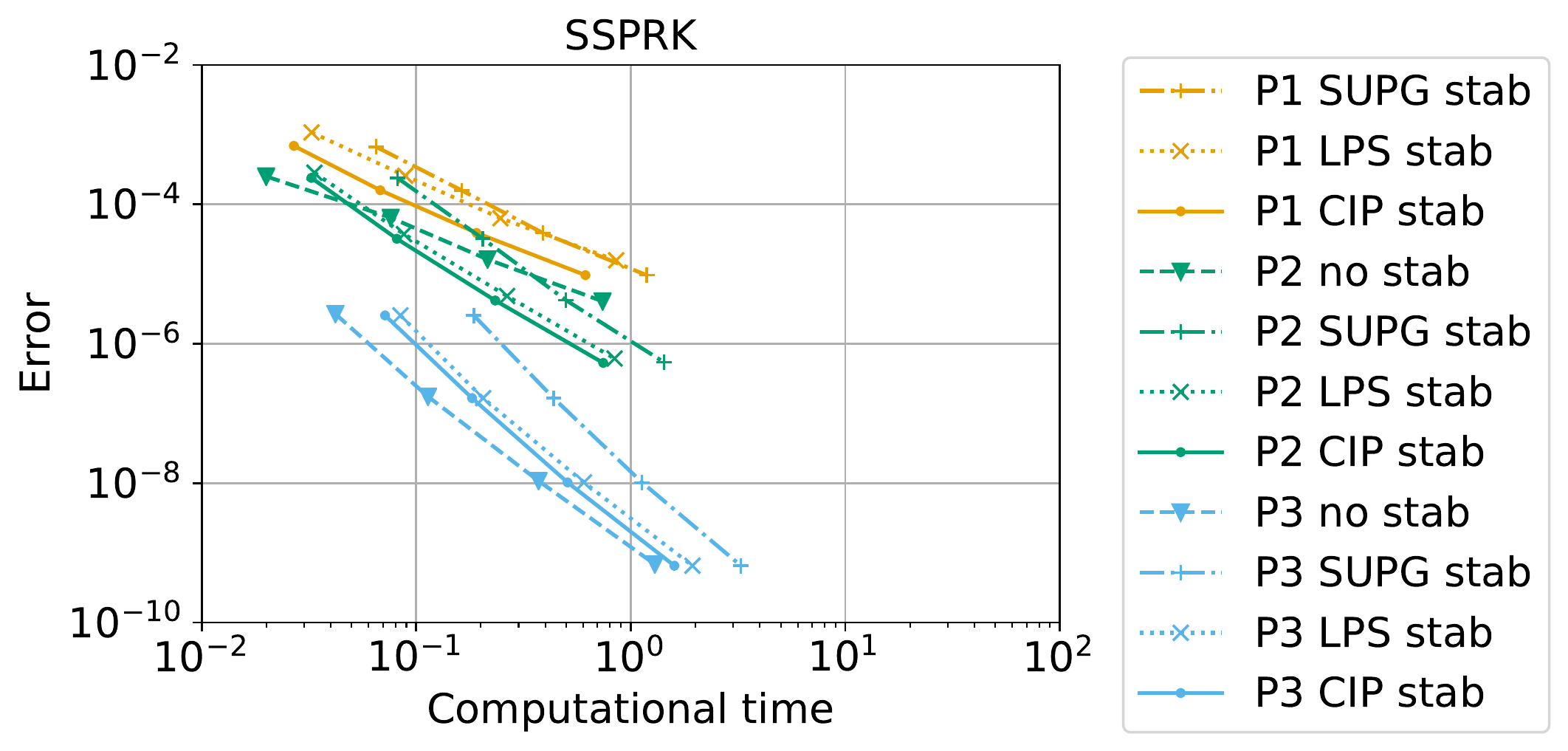}\\\vspace{2mm}
		\textit{Basic} elements\\ \vspace{2mm}
		\includegraphics[height=0.21\textheight,trim={0 0 62mm 0}, clip]{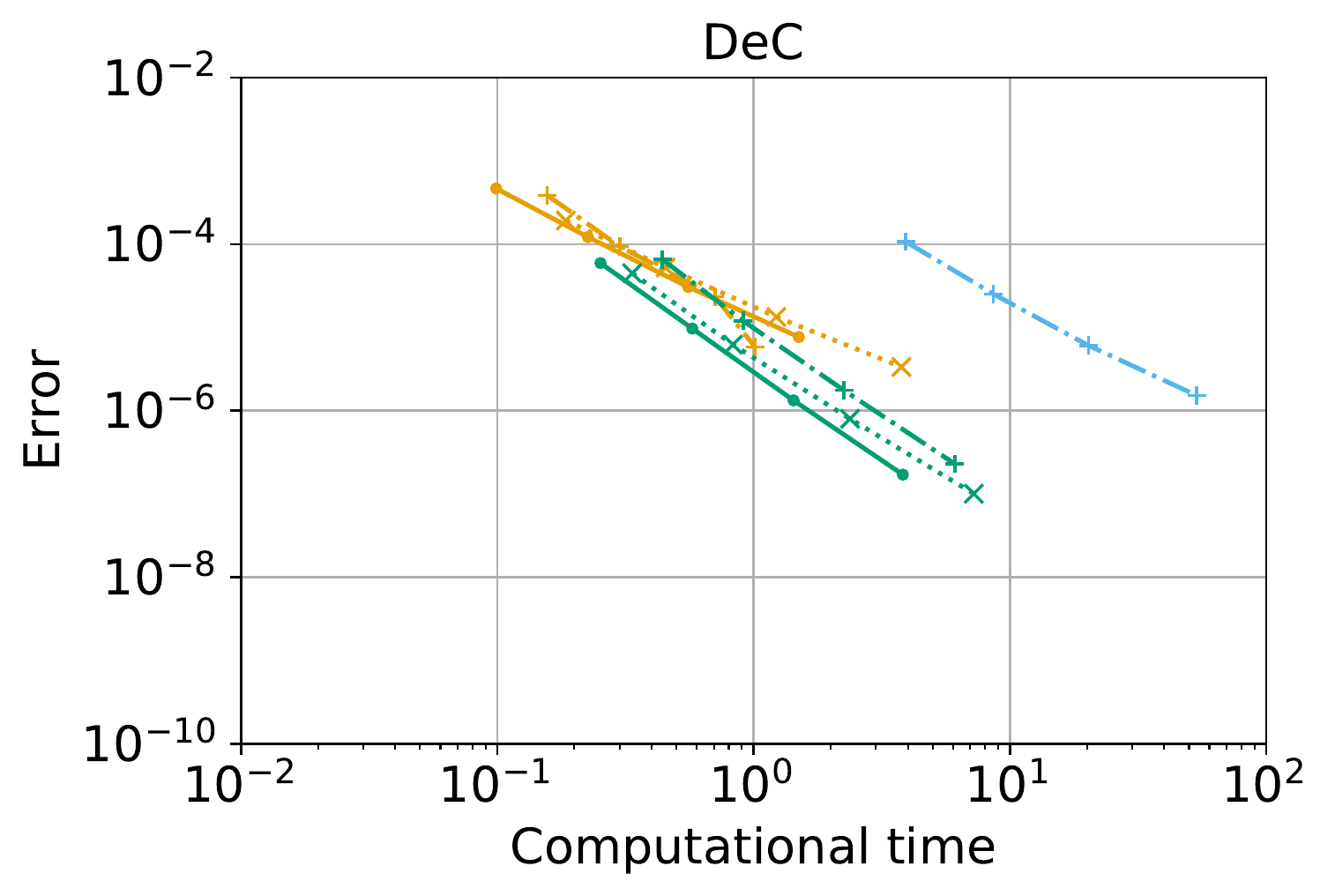}\hfill
		\includegraphics[height=0.21\textheight,trim={10mm 0 0 0}, clip]{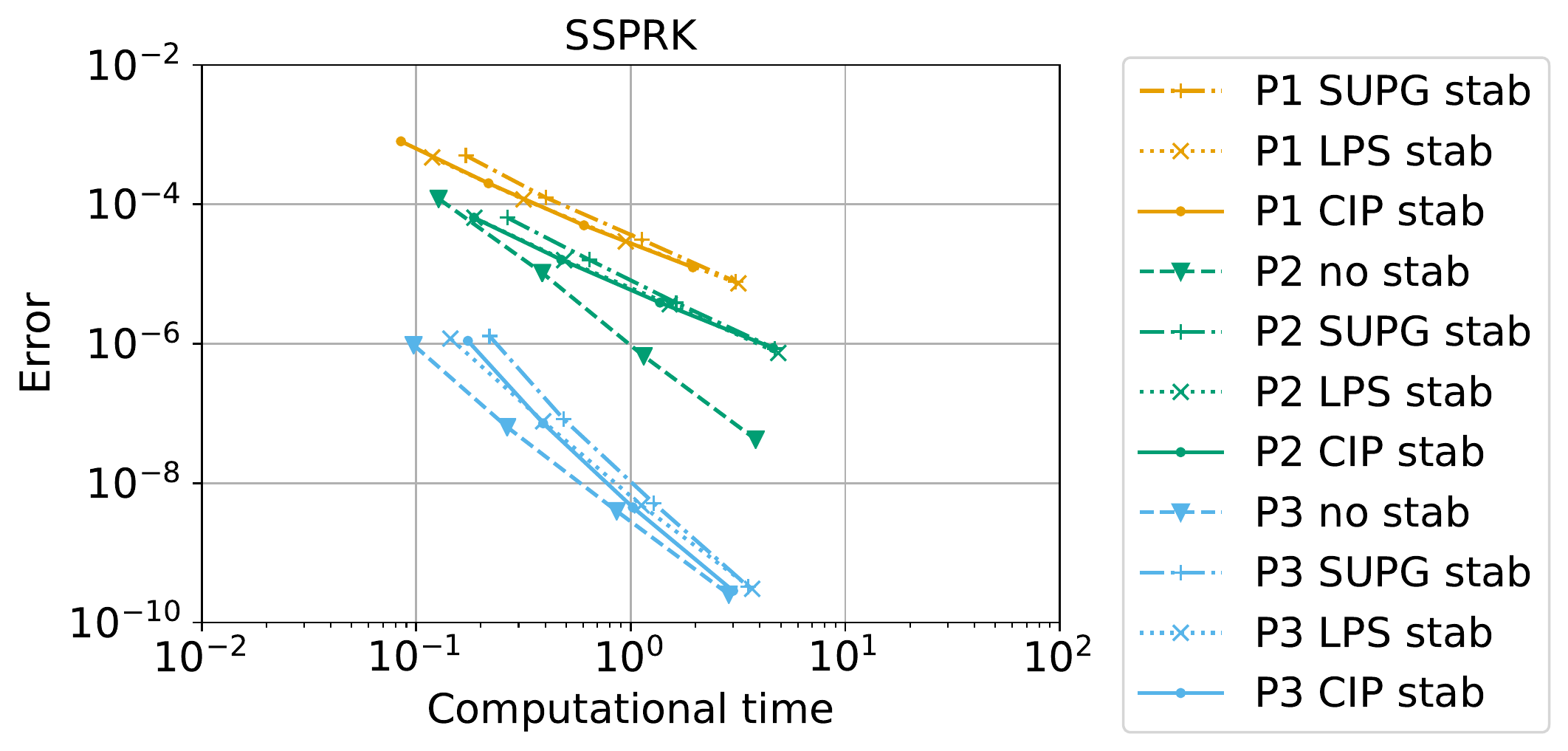}\\\vspace{2mm}
		\textit{Bernstein} elements\\\vspace{2mm}
		\includegraphics[height=0.21\textheight,trim={0 0 61mm 0}, clip]{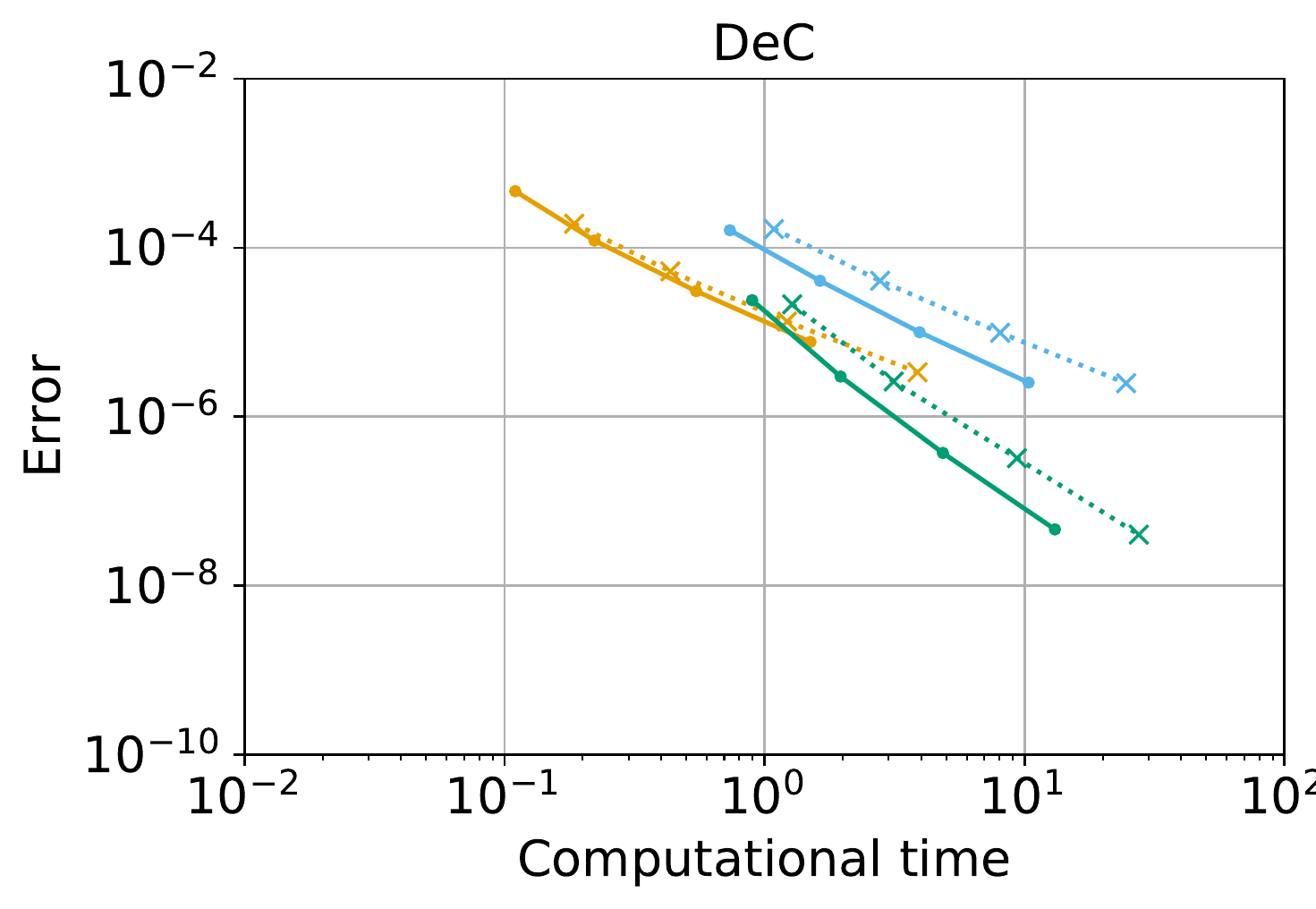}\hfill
		\includegraphics[height=0.21\textheight,trim={10mm 0 0 0}, clip]{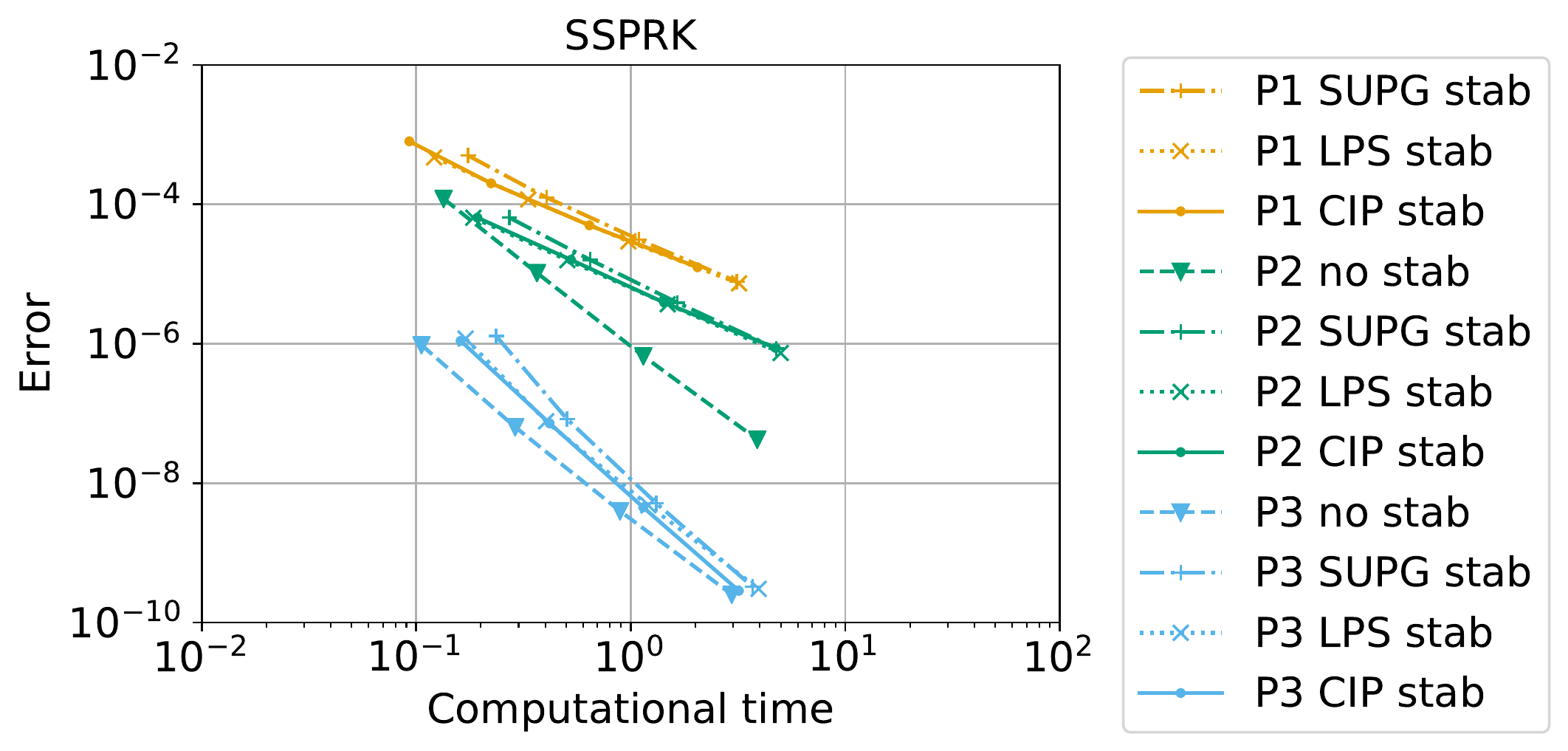}
	\end{center}
	\caption{Error for linear advection problem \eqref{lin_adv_1D} with respect to computational time for all elements and stabilization techniques: DeC on the left, SSPRK on the right}
	\label{fig:timeVsErrorLinAdv}
\end{figure}
Looking at the table   we can make the following observations. 
First of all, we remark that \revMR{ despite the  weak stability obtained in the spectral analysis, in practice the absence of damping makes it difficult to obtain converging results with a fixed CFL and for all $p$. 
For this reason, in the following we will only focus on stabilized methods.}
 
We observe otherwise  that almost all the stabilized scheme provide the expected order of accuracy. When the order is correct there are minor differences in the errors. 
There are however  few cases that \revMR{ fail in doing so and deserve some comments}. 
In particular, we notice the failure of DeC for \textit{basic} $\mathbb{P}_3$ and \textit{Bernstein} $\mathbb{B}_3$ polynomials and the SSPRK with \textit{basic} and \textit{Bernstein} $\mathbb{P}_2$ elements.  
\revMR{While disappointing, this negative  result is not completely new.  Indeed, in \cite{paola_svetlana,} obtaining correct convergence with DeC for some orders required both increasing the number of substebs,
thus making the method more expensive than the corresponding RK scheme, as well as including penalty terms on the jumps of higher order derivatives. Finally, note that this is in line with these methods falling in the 
family of ``broom'',  ``box'', and  thin striped shaped stability regions which we expect to be  difficult to use in practice}. 
\revMR{Concerning the stabilization of high order derivatives  this is also something a few authors advocate, see for instance the work by} 
Burman, Hansbo and collaborators 
\cite{Burman2020ACutFEmethodForAModelPressure,larson2019stabilizationHighOrderCut}.  
\revMR{ While this mayor explains the behavior observed, since we did not observe the need of including these terms for other cases than the DeC, we decided to 
focus on the simplest and most efficient approaches.}

\revMR{An interesting comparison is the one in \cref{fig:timeVsErrorLinAdv} where we plot }
the error of each method against computational time. \revMR{ Note that the simulations are all obtained  using }
the {\CFL} reported in \cref{tab:dispersion_cfl-RES}.  
In general, we can state that the \textit{cubature} elements obtain the best computational time as they are mass matrix free. On the other side, \textit{Bernstein} elements are slightly more expensive than \textit{basic} elements for DeC, because of the CFL restrictions that \cref{tab:dispersion_cfl-RES} requires.

Comparing time discretizations, \revMR{we see that despite the inversion of the mass matrix, SSPRK converges more rapidly than DeC. We think this is related to several reasons.
First of all,  the DeC CFL conditions are stricter,  and also DeC requires more stages. Even though not explicitly inverted, the mass matrix still needs to be assembled and multiplied to the solutions in the correction terms. 
Note however that the situation might radically change in the multidimensional case in which the mass matrix inversion in the SSPRK will  provide a much larger overhead.}

On the stabilization side, LPS and CIP behave very similarly (also their CFL do), but overall, the CIP is a little faster as it does not require the inversion of the mass matrix, for example, in DeC. 
\revMR{As expected, the SUPG stabilization requires more computational time, even if it often has larger CFs conditions. This is even clearer when using } \textit{cubature} elements, 
where \revMR{SUPG is the only case in which we still need to invert 
 the mass matrix with RK time stepping.} 

\subsection{Burgers' equation}

\revMR{We consider here application to a simple nonlinear problem to verify the applicability of the conditions obtained in the linear case.}
We test the  numerical schemes on the solution of the Burgers' equation
\begin{equation}
	\begin{cases}
		\partial_t u(x,t) + \partial_x \frac{u^2(x,t)}{2} =0 & (x,t)\in \Omega \times [0,t_f],\\
		u(x,0)=u_0(x), & x\in \Omega\\
		u(x_D,t)=g(x_D,t), & x_D\in \partial \Omega,
	\end{cases}\label{eq:Burgers}
\end{equation}
where $\Omega=[0,2]$ and $u_0(x) = -\tanh (4(x-1)) $ and $g(x,t)=u_{ex}(x,t)$ is the boundary condition. The exact solution is obtained using the method of   characteristics and reads $u_{ex}(x,t)=u_0(\chi)$ where
\begin{equation} \label{eq:characteristics}
	\chi = x - u_0(\chi)t
\end{equation}  
for all $(x,t)\in\Omega\times [0,t_f]$, solving the nonlinear equation \eqref{eq:characteristics} for $\chi$ at every point $(x,t)$. To obtain the exact solution we employed the Broyden method implemented in SciPy library \cite{2020SciPy-NMeth}.
Note that the  analytical solution shows a shock at time
\begin{equation}
t_s = -\frac{1}{\min\limits_{ x \in \Omega} u_0'(x)}=\frac{1}{4}.
\end{equation}
\revMR{This knowledge allows tho set  for this study  $t_f=0.5 t_s=0.125$, at which the solution is still smooth and the convergence of the higher order approximations can be investigated.
As before, in doing this we perform conformal refinement of the 1D grid, while paying attention to guarantee to use the same  number of degrees of freedom for   different $p$, and in particular taking:}
$\Delta x_2 = 2\Delta x_1 $ for $\mathbb{P}_2$ elements and $\Delta x_3 = 3 \Delta x_1 $ for $\mathbb{P}_3$ elements. 

Using the CFL and $\delta$ obtained in \cref{tab:dispersion_cfl-RES} we obtain the experimental order of convergence in \cref{tab:conv_order_Burger-RES}.

\begin{table}[H] 
\small  
 \begin{center} 
		\begin{tabular}{| c | c || c | c | c || c | c | c || c | c | c | }  
	     \hline 
	     \multicolumn{2}{|c||}{Element $\&$ }  & \multicolumn{3}{|c||}{No stabilization }  & \multicolumn{3}{|c||}{LPS}  & \multicolumn{3}{|c|}{CIP}  \\ \hline 
	     \multicolumn{2}{|c||}{ Time scheme }  &  $\mathbb{P}_1$ & $\mathbb{P}_2$ & $\mathbb{P}_3$   & $\mathbb{P}_1$ & $\mathbb{P}_2$ & $\mathbb{P}_3$   & $\mathbb{P}_1$ & $\mathbb{P}_2$ & $\mathbb{P}_3$ \\ \hline \hline 
       \parbox[t]{2mm}{\multirow{2}{*}{\rotatebox[origin=c]{90}{\centering  Cub.}}}              &  SSPRK &  /  & 1.99 & 3.71 & 2.05 & 2.85 & 3.67 & 2.05 & 2.85 & 3.68 \\ 
               &  DeC &  /  & 1.99 & 3.71 & 2.06 & 2.85 & 3.57 & 2.06 & 2.85 & 3.69 \\ 
         \hline 
        \parbox[t]{2mm}{\multirow{2}{*}{\rotatebox[origin=c]{90}{Basic}}}              &  SSPRK &  /  & 1.99 & 3.82 & 2.07 & 2.56 & 3.66 & 2.06 & 2.48 & 3.66 \\ 
               &  DeC &  /  &  /  &  /  & 2.7 & 2.92 &  /  & 2.59 & 2.85 &  /  \\ 
         \hline 
       \parbox[t]{2mm}{\multirow{2}{*}{\rotatebox[origin=c]{90}{\centering  Bern.}}}              &  SSPRK &  /  & 1.99 & 3.82 & 2.07 & 2.56 & 3.66 & 2.06 & 2.48 & 3.66 \\ 
               &  DeC &  /  &  /  &  /  & 2.7 & 2.9 & 1.41 & 2.59 & 2.87 & 1.37 \\ 
         \hline 
        \end{tabular} 
    \end{center} 
     \caption{Summary table of convergence order, using coefficients obtained in \cref{tab:dispersion_cfl-RES} } \label{tab:conv_order_Burger-RES}
\end{table}%

\begin{figure}[h!]
	\begin{center}
		\textit{Cubature} elements\\\vspace{2mm}
		\includegraphics[height=0.21\textheight,trim={0 0 61mm 0}, clip]{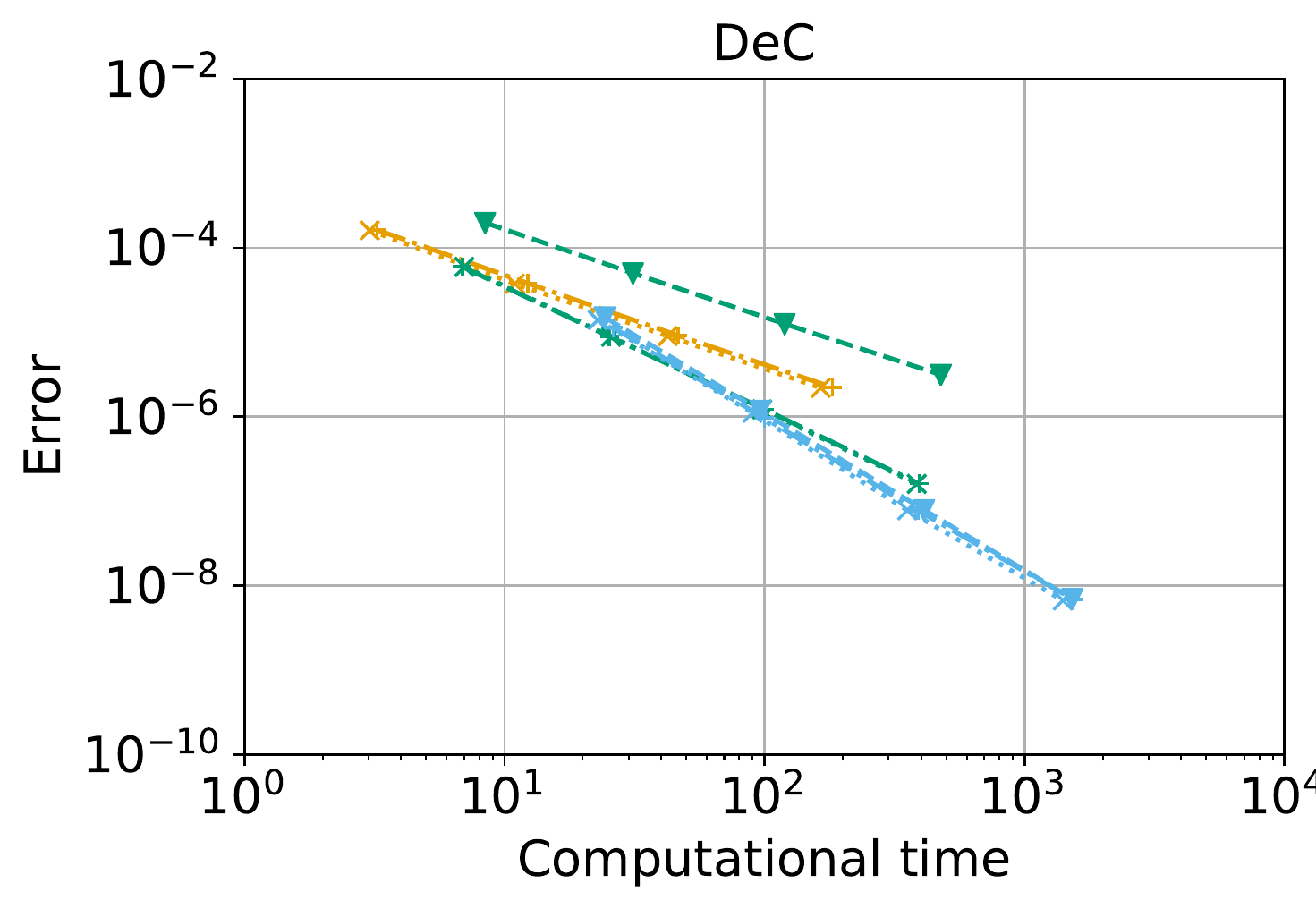} \hfill
		\includegraphics[height=0.21\textheight,trim={10mm 0 0 0}, clip]{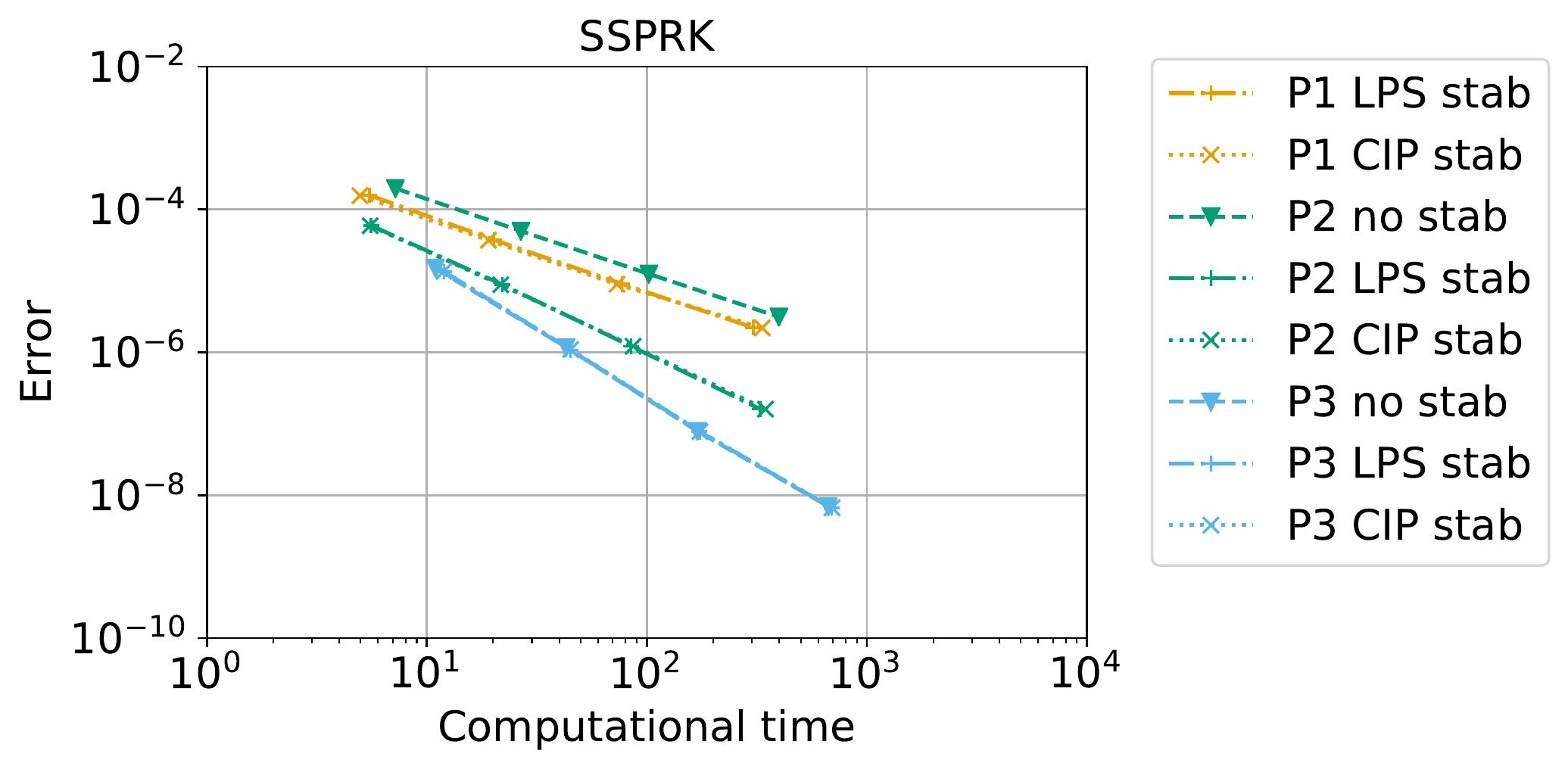}\\\vspace{2mm}
		\textit{Basic} elements\\ \vspace{2mm}
		\includegraphics[height=0.21\textheight,trim={0 0 61mm 0}, clip]{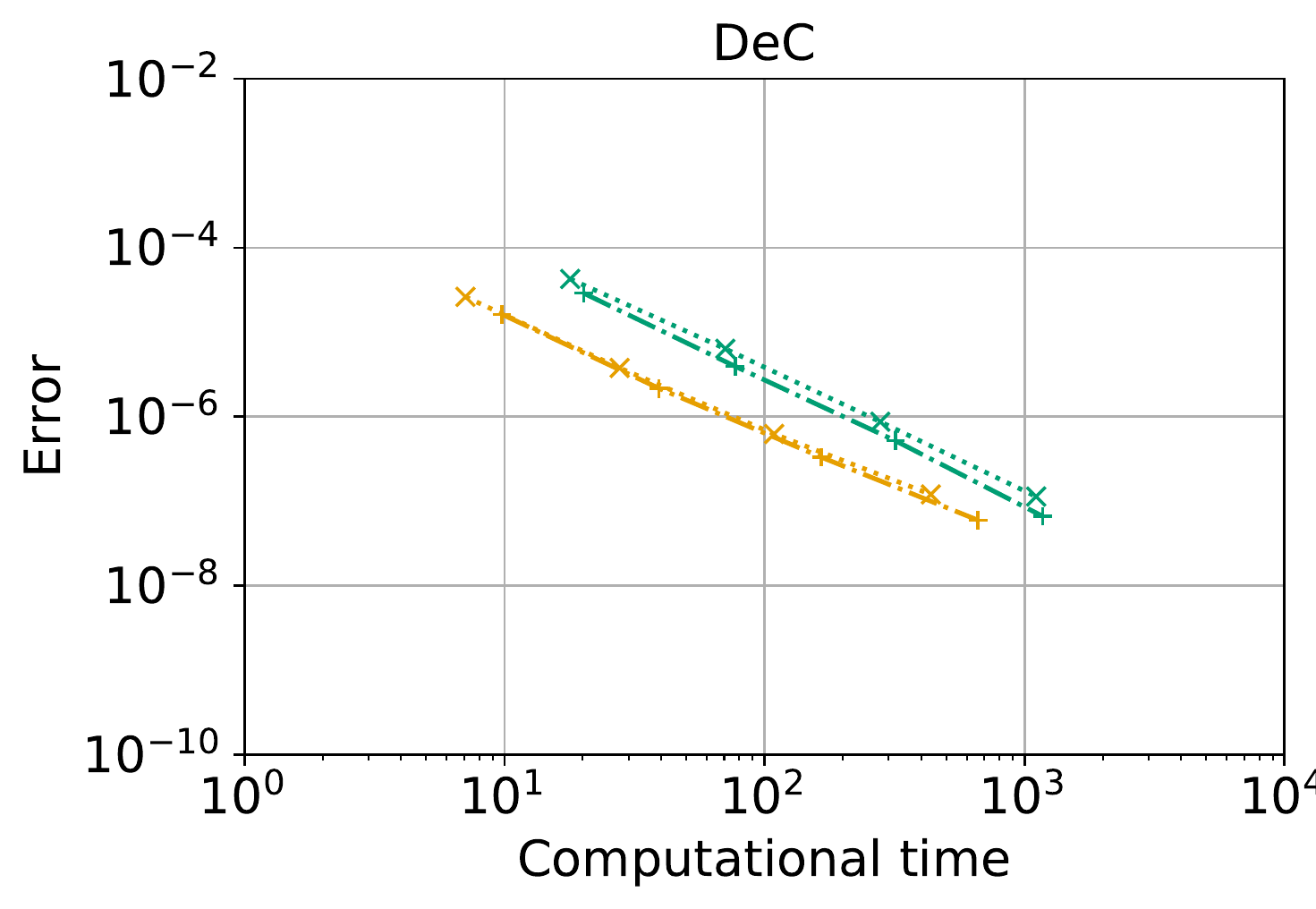}\hfill
		\includegraphics[height=0.21\textheight,trim={10mm 0 0 0}, clip]{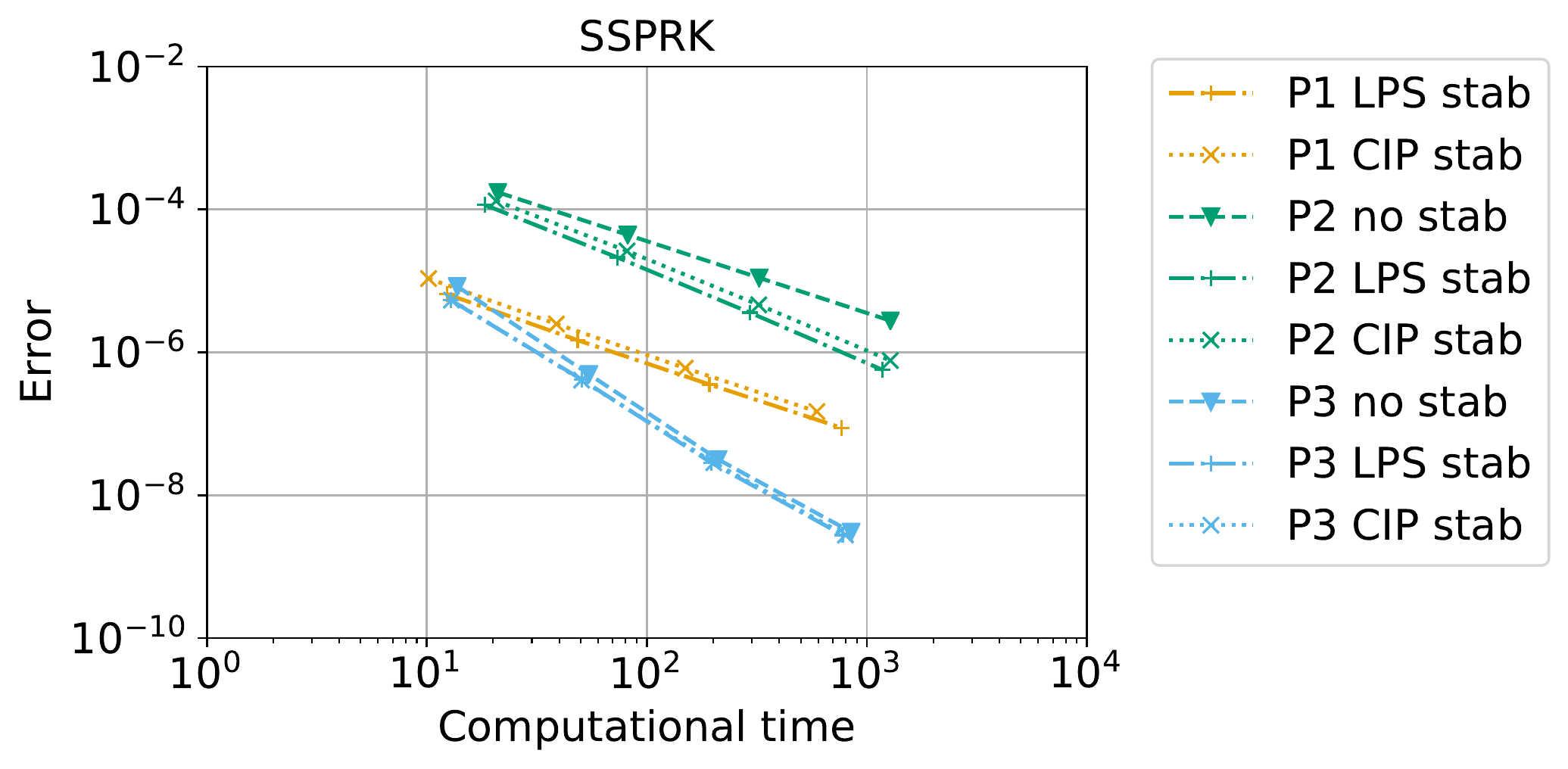}\\\vspace{2mm}
		\textit{Bernstein} elements\\\vspace{2mm}
		\includegraphics[height=0.21\textheight,trim={0 0 61mm 0}, clip]{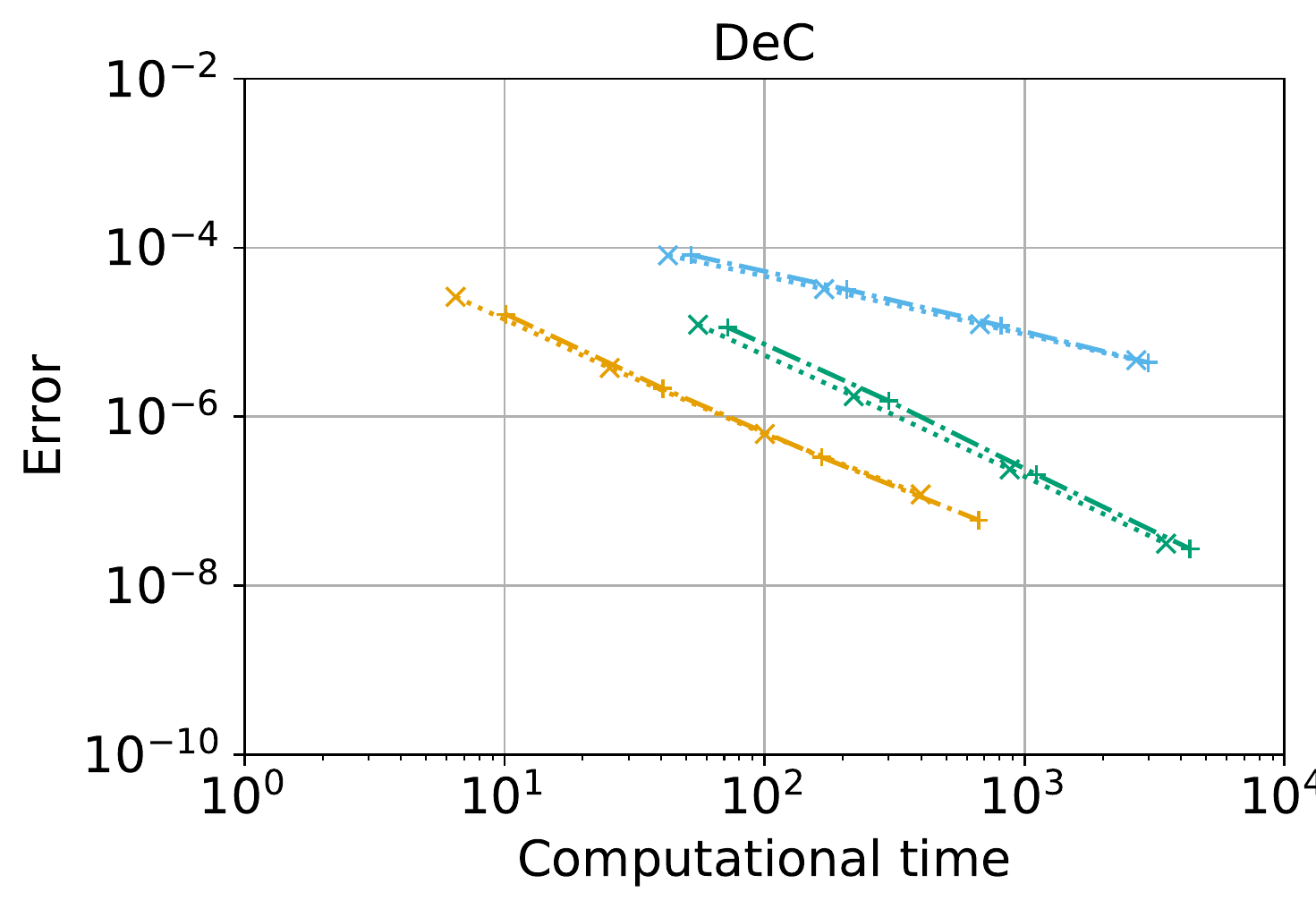}\hfill
		\includegraphics[height=0.21\textheight,trim={10mm 0 0 0}, clip]{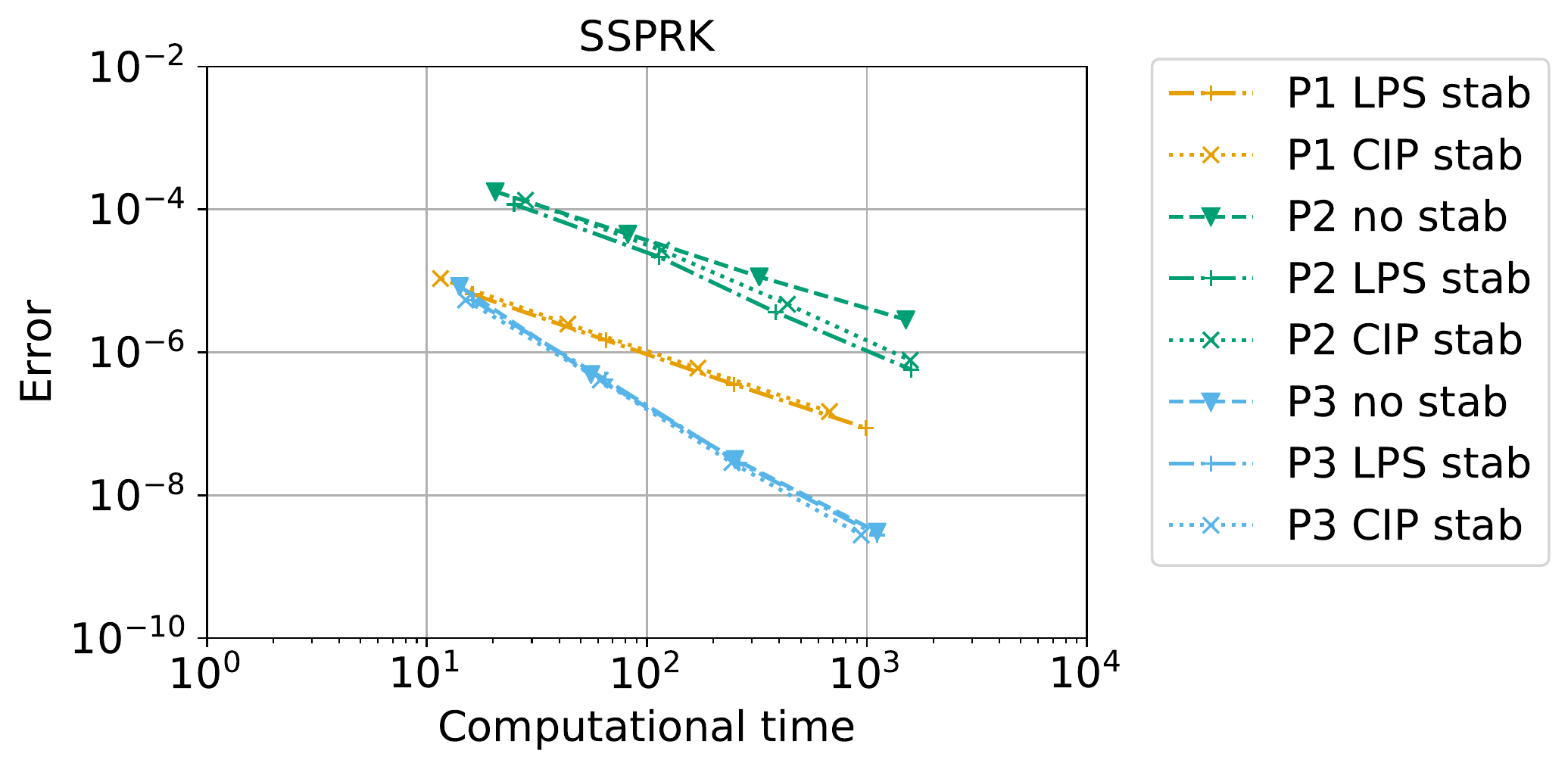}
	\end{center}
	\caption{Error for Burgers' equation \eqref{eq:Burgers} with respect to computational time for all elements and stabilization techniques: DeC on the left, SSPRK on the right}
	\label{fig:timeVsErrorBurgers}
\end{figure}

The results are very similar to the ones obtained for the linear advection case. There is a small improvement in \textit{basic }and \textit{Bernstein }$\mathbb{P}_2$ SSPRK cases, while the DeC \textit{basic }and \textit{Bernstein }$\mathbb{P}_3$ cases are even worse than the linear advection ones. 
The DeC $\mathbb{P}_1$ \textit{basic }and \textit{Bernstein} cases show a super--convergent behavior. 
The interested reader will find the convergence plots for all the combinations on the repository \cite{TorloMichel2020git}. Here we 
focus  on the comparison between error and computational time, reported in \cref{fig:timeVsErrorBurgers}. 
Again for \textit{cubature} elements it is clear the advantage in using high order methods, in particular for SSPRK methods, which has less stages than DeC. 
For this test, we only compare CIP and LPS and they \revMR{systematically out-perform  SUPG}. For these two, the difference in computational time is very minimal \revMR{ for all element choices.
This may change in the multidimensional case where the LPS may be penalized on elements requiring the inversion of the full mass matrix.}

For DeC \textit{basic} and \textit{Bernstein} $\mathbb{P}_1$ elements, the superconvergence of 
the second order schemes  makes them the best in their category, see \cref{tab:conv_order_Burger-RES}. For SSPRK the expected order of convergence of fourth order scheme shows how the high order accurate methods can provide the fastest and most precise solutions.

\subsection{Shallow water equations}

As a final application we consider the  non linear shallow water equations:
\begin{equation}
    \left  \{
    \begin{array}{ll}
    	\partial_t h + \partial_x (hu) & = 0,  \\
    	\partial_t (hu) + \partial_x (hu^2 +g\frac{h^2}{2} ) + \Phi & =0,  \label{eq_GN2}
	\end{array}	\quad x\in \Omega, \, t \in [0,5].
    \right .
\end{equation}
Here, $h$ is the water elevation, $u$ the velocity field, $g$ the gravitational acceleration. 
\revMR{We will solve the system  on the domain $\Omega = [0,200]$, and  add the source term  
  $\Phi = \Phi(x,t)$ in order to impose the  solution to be equal to
  \begin{equation}
\begin{cases}
    	h_{ex}(x,t) =  h_0 +\epsilon h_0 \text{sech}^2(\kappa(x-ct)), \\
    	u_{ex}(x,t) = c \left( 1-\frac{h_0}{h_{ex}(x,t)} \right), \\
    	\kappa = \sqrt{\frac{3\epsilon}{4h_0^2 (1+\epsilon)}} , \quad c=\sqrt{gh_0(1+\epsilon)}.
\end{cases}
    \label{eq_GN_ex}
\end{equation}
Following the classical manufactured solution method, we set}
%
\begin{align*}
    	\Phi(x,t) & = -\left[ \partial_t \left(h_{ex}(x,t)u_{ex}(x,t) \right)  + \partial_x \left( h_{ex}(x,t)u^2_{ex}(x,t) +g\frac{h^2_{ex}(x,t)}{2}   \right)  \right] \\
    			  & = -\left[ h_{ex} (\partial_t u_{ex} +u_{ex} \partial_x u_{ex} +g\partial_x h_{ex}   )  \right].
\end{align*}
For our study, we set $\epsilon = 1.2$, $h_0=1$ and the initial and Dirichlet boundary condition given by the exact solution at time 0 and at the borders of the domain.

We discretize the mesh with uniform intervals of length $\Delta x$, and \revMR{as before we perform a grid convergence by respecting the constraint
 $\Delta x_2 = 2\Delta x_1 $ for $\mathbb{P}_2$ elements and $\Delta x_3 = 3 \Delta x_1 $ for $\mathbb{P}_3$ elements.} 
In \cref{tab:conv_order_nlsw-RES} we show the convergence orders for this shallow water problem with the CFL and $\delta$ coefficients found in \cref{tab:dispersion_cfl-RES}.
\begin{table}[H] 
\small  
 \begin{center} 
		\begin{tabular}{| c | c || c | c | c || c | c | c || c | c | c | }  
	     \hline 
	     \multicolumn{2}{|c||}{Element $\&$ }  & \multicolumn{3}{|c||}{No stabilization }  & \multicolumn{3}{|c||}{LPS}  & \multicolumn{3}{|c|}{CIP}  \\ \hline 
	     \multicolumn{2}{|c||}{ Time scheme }  & $\mathbb{P}_1$ & $\mathbb{P}_2$ & $\mathbb{P}_3$   & $\mathbb{P}_1$ & $\mathbb{P}_2$ & $\mathbb{P}_3$   & $\mathbb{P}_1$ & $\mathbb{P}_2$ & $\mathbb{P}_3$  \\ \hline \hline 
       \parbox[t]{2mm}{\multirow{2}{*}{\rotatebox[origin=c]{90}{\centering  Cub.}}}              &  SSPRK &  /  & 1.96 & 5.17 & 2.26 & 2.69 & 5.02 & 2.39 & 2.68 & 5.05 \\ 
               &  DeC &  /  & 1.97 & 5.17 & 2.28 & 2.65 & 4.79 & 2.7 & 2.66 & 5.07 \\ 
         \hline 
        \parbox[t]{2mm}{\multirow{2}{*}{\rotatebox[origin=c]{90}{Basic}}}              &  SSPRK &  /  & 1.98 & 5.54 & 1.94 & 2.31 & 4.93 & 1.95 & 2.29 & 4.98 \\ 
               &  DeC &  /  &  /  &  /  & 2.23 & 2.74 &  /  & 2.01 & 2.58 &  /  \\ 
         \hline 
       \parbox[t]{2mm}{\multirow{2}{*}{\rotatebox[origin=c]{90}{\centering  Bern.}}}              &  SSPRK &  /  & 1.97 & 2.44 & 1.94 & 2.07 & 2.19 & 1.95 & 2.09 & 2.21 \\ 
               &  DeC &  /  &  /  &  /  & 2.23 & 2.0 & 2.0 & 2.01 & 2.0 & 1.98 \\ 
         \hline 
        \end{tabular} 
    \end{center} 
     \caption{Summary tab of convergence order, using coefficients obtained by minimizing $\eta_u$ } \label{tab:conv_order_nlsw-RES}
\end{table}%

\begin{figure}[h!]
	\begin{center}
		\textit{Cubature} elements\\\vspace{2mm}
		\includegraphics[height=0.21\textheight,trim={0 0 61mm 0}, clip]{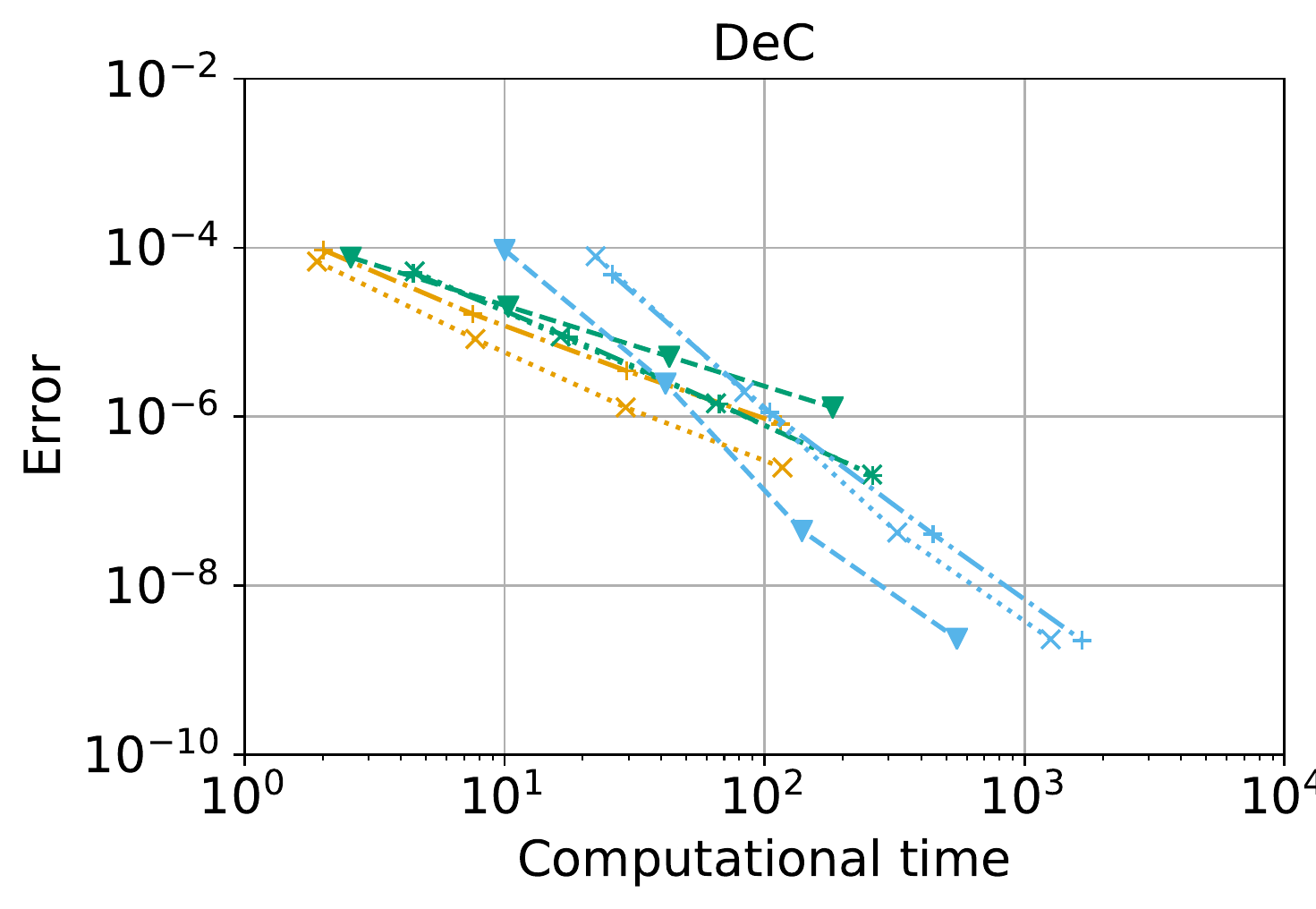} \hfill
		\includegraphics[height=0.21\textheight,trim={10mm 0 0 0}, clip]{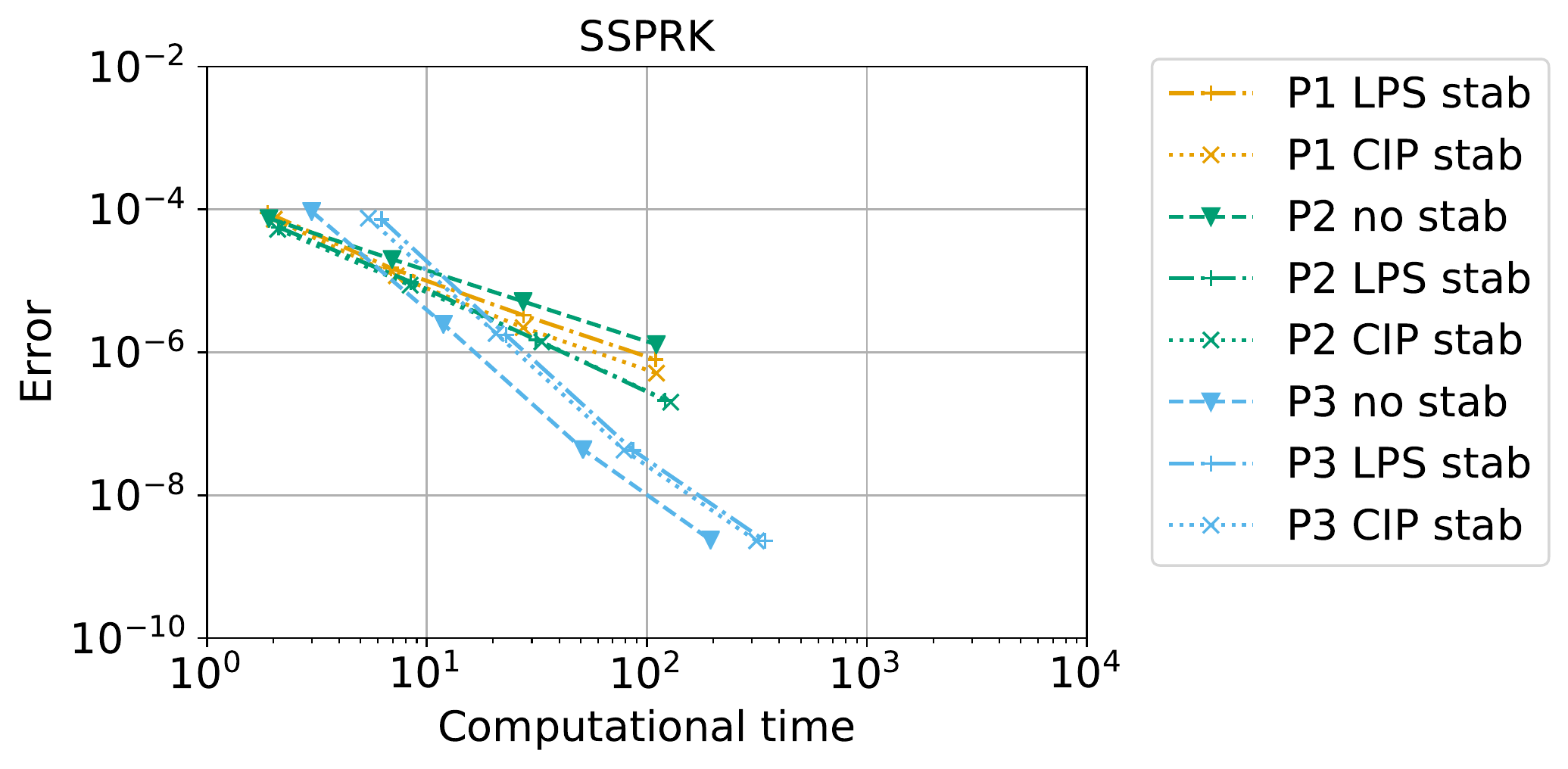}\\\vspace{2mm}
		\textit{Basic} elements\\ \vspace{2mm}
		\includegraphics[height=0.21\textheight,trim={0 0 61mm 0}, clip]{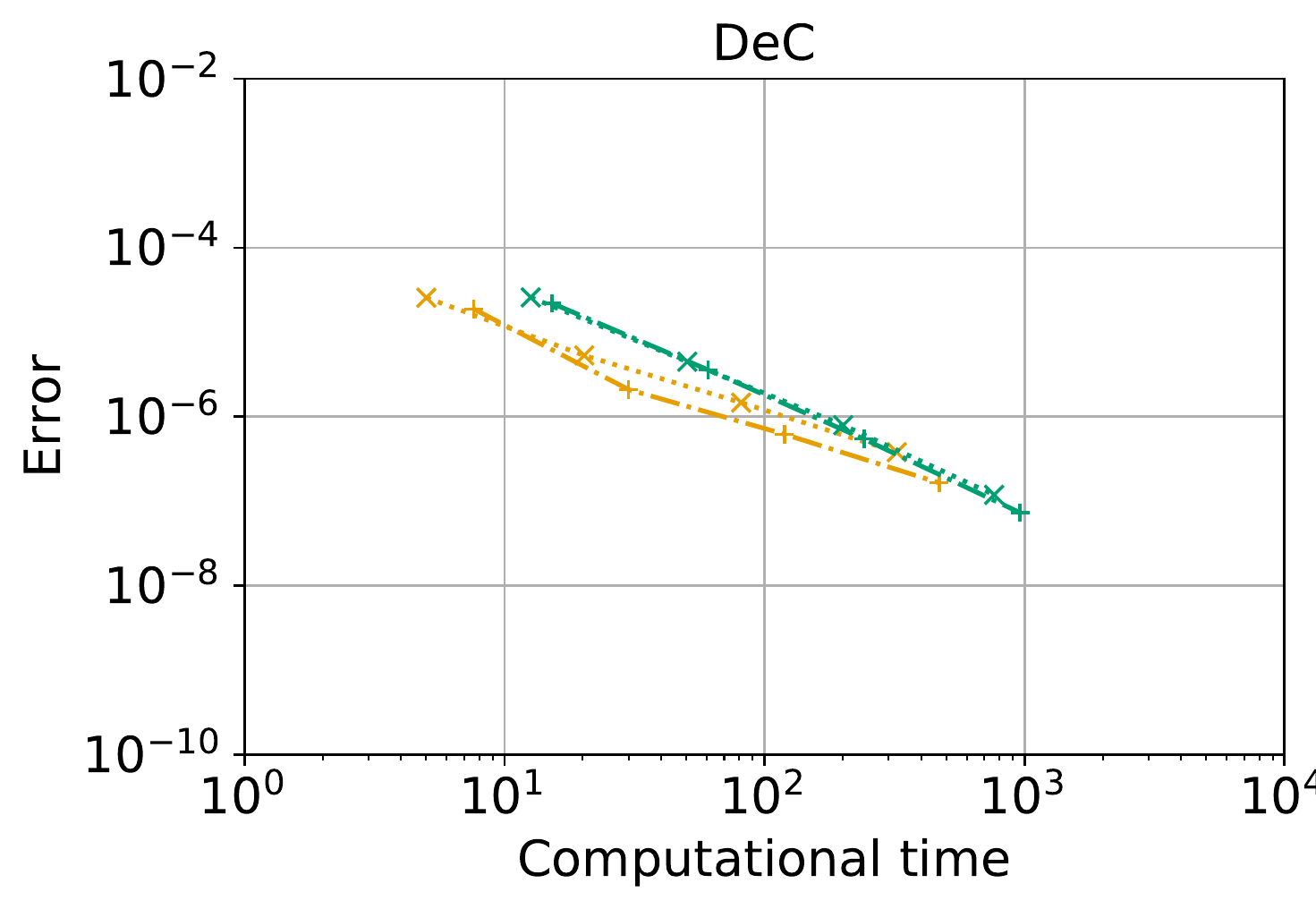}\hfill
		\includegraphics[height=0.21\textheight,trim={10mm 0 0 0}, clip]{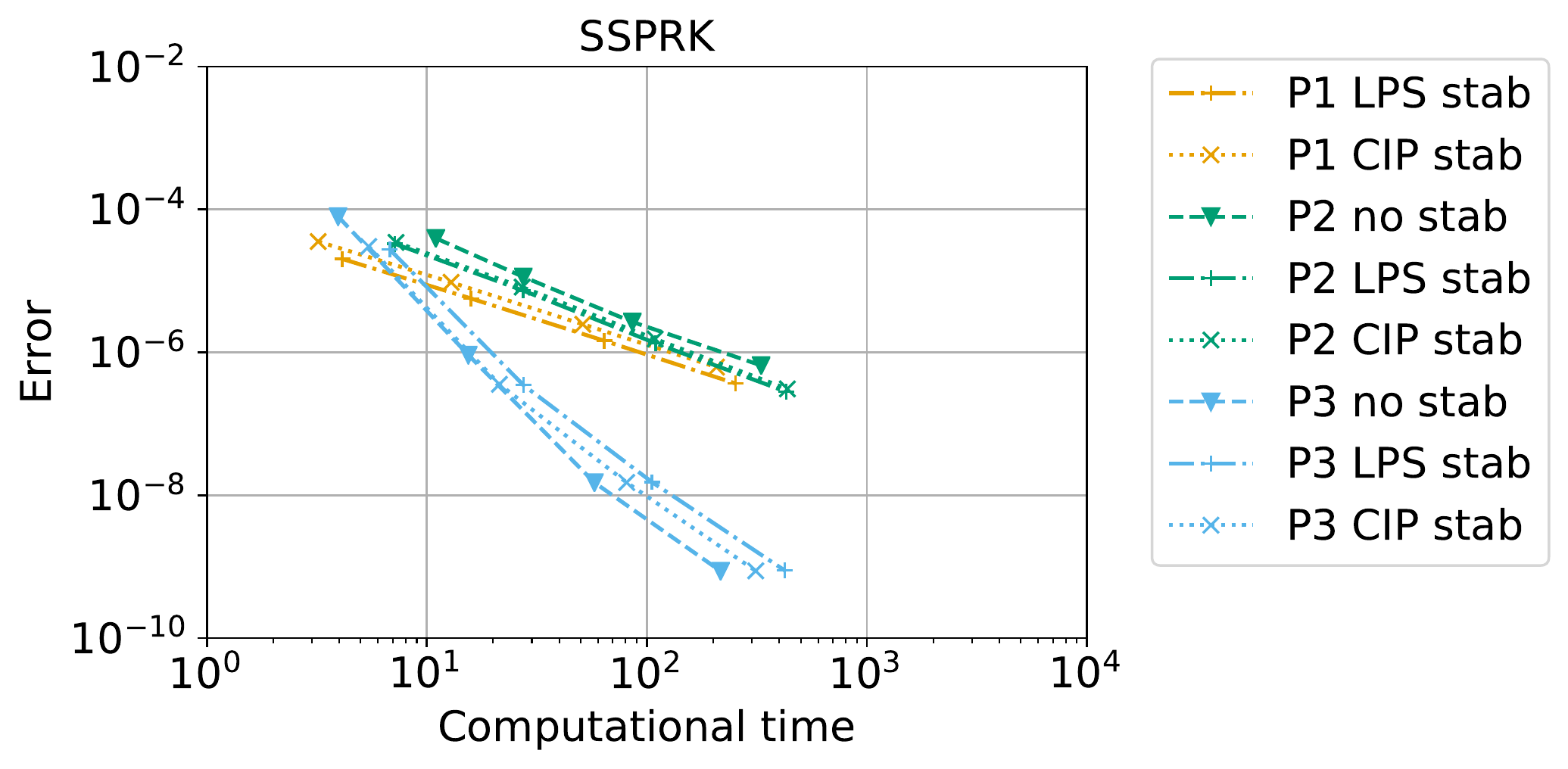}\\\vspace{2mm}
		\textit{Bernstein} elements\\\vspace{2mm}
		\includegraphics[height=0.21\textheight,trim={0 0 61mm 0}, clip]{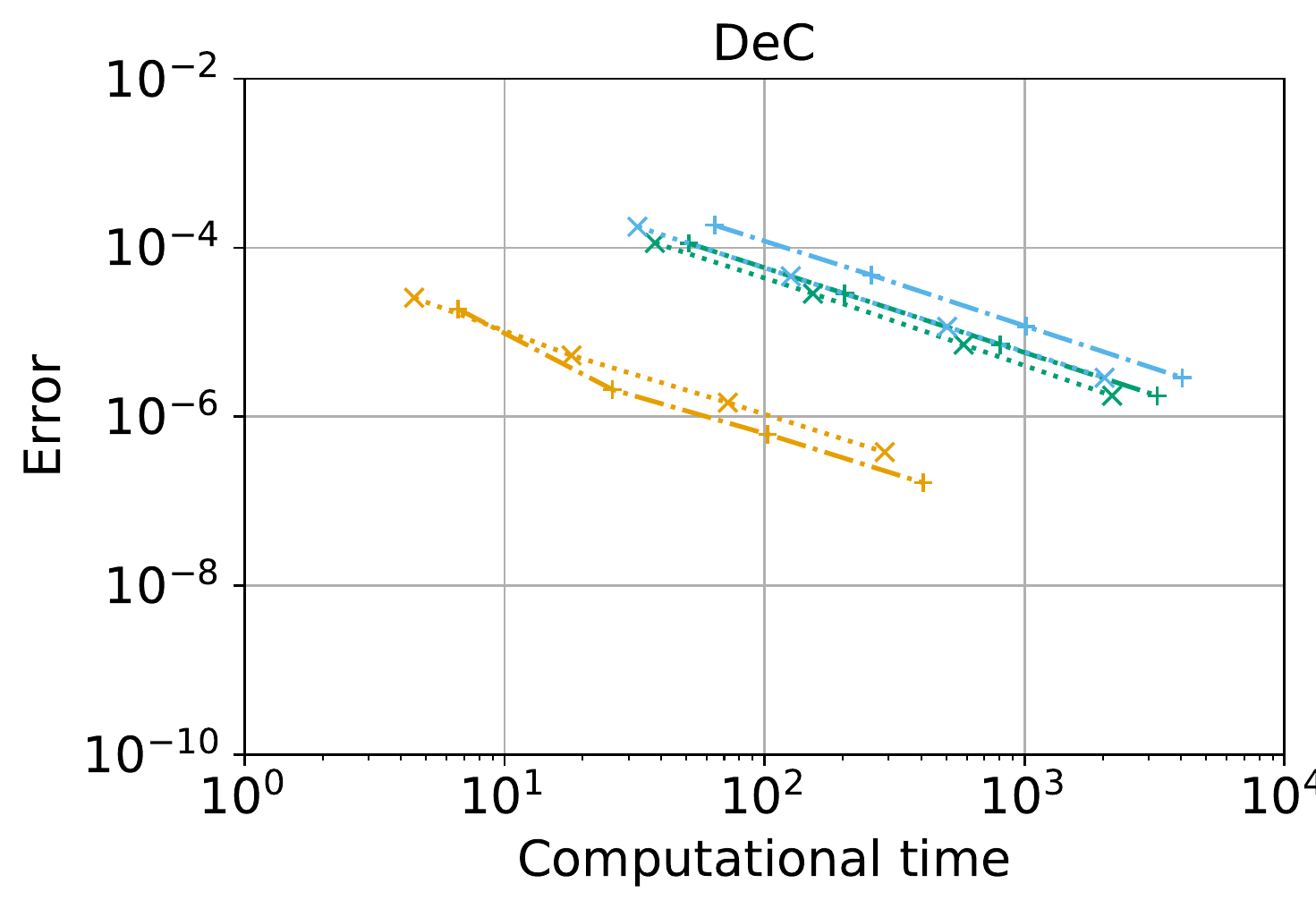}\hfill
		\includegraphics[height=0.21\textheight,trim={10mm 0 0 0}, clip]{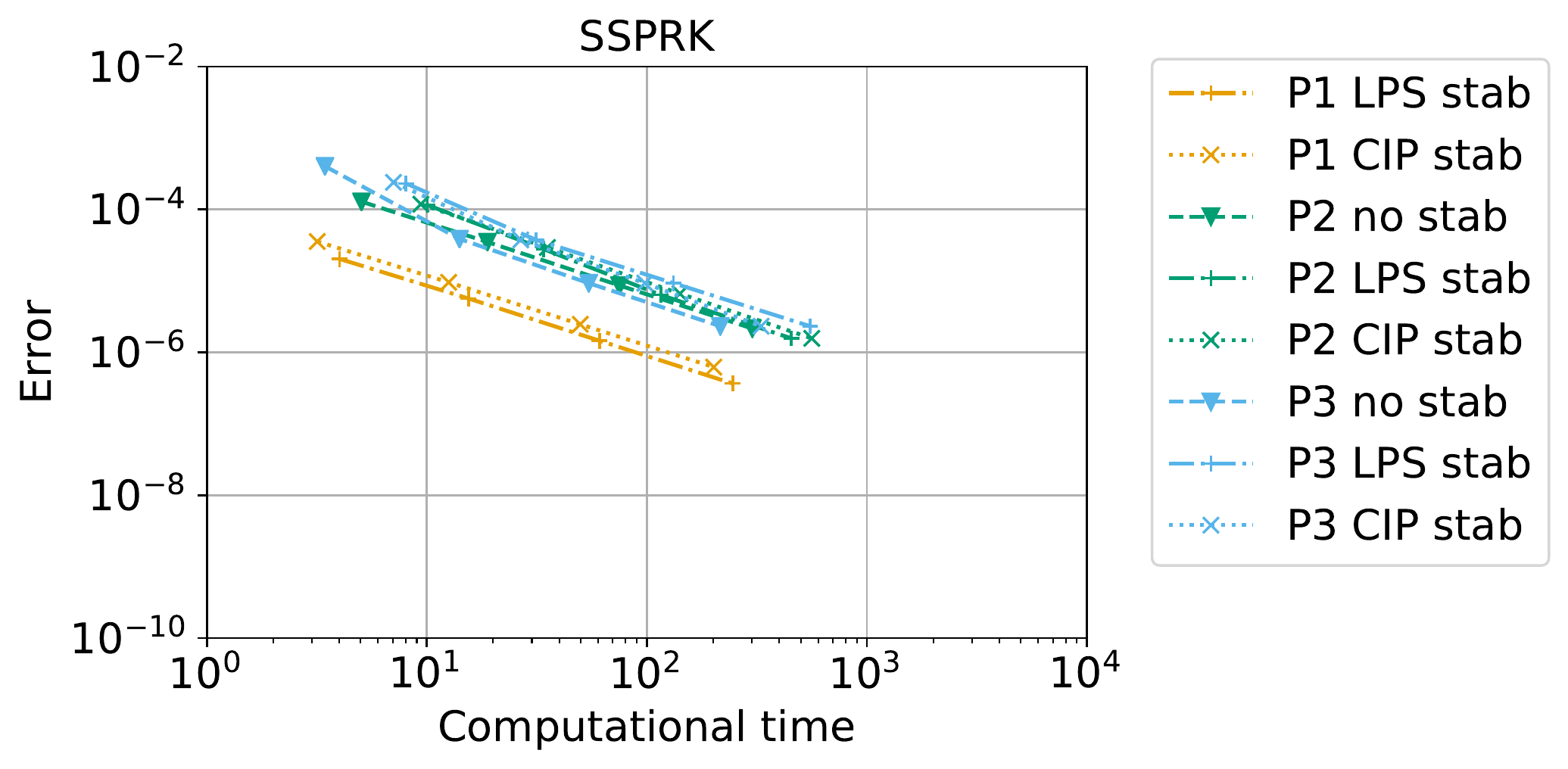}
	\end{center}
	\caption{Error for Shallow Water equations \eqref{eq_GN2} with respect to computational time for all elements and stabilization techniques: DeC on the left, SSPRK on the right}
	\label{fig:timeVsErrorSW}
\end{figure}

\revMR{The results obtained are similar to those of the other cases.  The convergence rates are at least the expected ones with {\it cubature} elements while we still see problems with
DeC and \textit{basic }elements in the fourth order case, as well as with {\it Bernstein} polynomials for both  $\mathbb{P}_2$ and $\mathbb{P}_3$. On the other hand,
some superconvergence is measured in the $\mathbb{P}_3$ case with both \textit{cubature} and \textit{basic} elements. This  creates an even larger bias in the error-cpu time plots, \cref{fig:timeVsErrorSW}, in favor of these higher
polynomial degrees.}


\section{Conclusion}\label{sec:conclusion}
In summary, we propose a comparison of high order continuous Galerkin methods with stabilization techniques for hyperbolic problems. 
On the linear advection equation, we perform a Fourier analysis on the spatial discretization, then a von Neumann analysis on the space--time discretization given by each combination of stabilization, time discretization and finite elements. 
This provides reliable parameters and CFL conditions for all the mentioned methods that can be used both in the linear advection case and in nonlinear problems, as the Burgers' and shallow water simulations showed.

\revMR{The Fourier analysis is limited to one dimensional problems (or structured multidimensional meshes), so the main ongoing development  is the verification 
of the properties of the methods studied in a multidimensional setting 
%
 based on the approximation choices}
  suggested   e.g. in \cite{article_fekete,article_cubature_2001,article_cubature_2006} and references therein. 

\section*{Acknowledgment}
This work was performed within the Ph.D. project of Sixtine Michel: ``Evaluation of coastal and urban submersion risks'', supported by INRIA and the BRGM, co-funded by in INRIA--Bordeaux Sud--Ouest and the Conseil R\'egional de la Nouvelle Aquitaine.
Mario Ricchiuto and Davide Torlo have been supported by team CARDAMOM in INRIA--Bordeaux Sud--Ouest.
Davide Torlo and R\'emi Abgrall have been supported by the Swiss National Foundation grant No 200020\_175784. 

\appendix
\section{Time schemes}\label{app:timeCoefficients}
In this appendix we introduce the time integration coefficients used in this work, to make the study fully reproducible. In  \cref{tab:Butcher} there are the RK coefficients, in \cref{tab:ButcherSSPRK} the SSPRK coefficients and in \cref{tab:DeCcoeff} the DeC coefficients.
\begingroup
\setlength{\tabcolsep}{6pt} 
\renewcommand{\arraystretch}{1.2} 
\begin{table}[h!]
	\centering
\begin{tabular}{ c c  c }
	\begin{tabular}{| c |  c c  | }
		\hline 
		\multicolumn{3}{|c|}{\textit{RK2}} \\ \hline
		 $\alpha$  & 1   &  \\ \hline
		$\beta$ & $\frac{1}{2}$ & $\frac{1}{2}$  \\ \hline
	\end{tabular} &
	\begin{tabular}{|c| c c c  | }
		\hline 
		\multicolumn{4}{|c|}{\textit{RK3}} \\ \hline
		$\alpha$  & $\frac{1}{2}$   &     & \\ 
		& -1 & 2  & \\ \hline
		$\beta$ & $\frac{1}{6}$ & $\frac{2}{3}$ & $\frac{1}{6}$ \\ \hline
	\end{tabular} &
	\begin{tabular}{| c | c  c c c | }
		\hline 
		\multicolumn{5}{|c|}{\textit{RK4}} \\ \hline
		$\alpha$ &  $\frac{1}{2}$   & & &   \\ 
		& 0 & $\frac{1}{2}$  & & \\ 
		& 0 & 0    & 1 &   \\ \hline 
		$\beta$ & $\frac{1}{6}$ & $\frac{1}{3}$ & $\frac{1}{3}$ & $\frac{1}{6}$ \\ \hline
	\end{tabular}
\end{tabular}
\caption{Butcher Tableau of RK methods}\label{tab:Butcher}
\end{table}
\begin{table}[h!]
\begin{center}
		\begin{tabular}{| c  c c  |c c c  | }
		\hline 
		\multicolumn{6}{|c|}{\textit{SSPRK(3,2)} by \cite{shu-1988}} \\ \hline
		\multicolumn{3}{|c|}{$\gamma$} & \multicolumn{3}{|c|}{$\mu$} \\ \hline
		1 &  &        &    $\frac{1}{2}$  &   &   \\ 
		0 & 1 &   &   0 & $\frac{1}{2}$  & \\ 
		$\frac{1}{3}$ & 0 & $\frac{2}{3}$ &   0 & 0 & $\frac{1}{3}$ \\ \hline
		\multicolumn{6}{|c|}{CFL = 2.} \\ \hline
	\end{tabular} \qquad
	\begin{tabular}{| c  c c c  |  c c c c | }
		\hline 
		\multicolumn{8}{|c|}{\textit{SSPRK(4,3)} by \cite[Page 189]{Ruuth-2006}} \\ \hline
		\multicolumn{4}{|c|}{$\gamma$} & \multicolumn{4}{|c|}{$\mu$} \\ \hline
		1 &   &       &      & $\frac{1}{2}$ &  &     & \\ 
		0 & 1 &       &      & 0 & $\frac{1}{2}$ & &  \\ 
		$\frac{2}{3}$ & 0 & $\frac{1}{3}$ &      & 0 & 0 & $\frac{1}{6}$     & \\ 
		0   & 0 &  0  & 1    & 0 & 0 & 0    & $\frac{1}{2}$  \\ \hline
		\multicolumn{8}{|c|}{CFL = 2.} \\ \hline
	\end{tabular}\\
	\begin{tabular}{| c c c c c |}
	\hline 
	\multicolumn{5}{|c|}{\textit{SSPRK(5,4)} by \cite[Table 3]{Ruuth-2006} } \\ \hline
	\multicolumn{5}{|c|}{$\gamma$}   \\ \hline
	1 &   &  &     &                                                \\ 
	0.444370493651235 & 0.555629506348765 &   &    &                   \\ 
	0.620101851488403 & 0                 & 0.379898148511597 & &       \\ 
	0.178079954393132 & 0 &  0  & 0.821920045606868 &                 \\ 
	0 & 0 & 0.517231671970585 & 0.096059710526147 & 0.386708617503269  \\ \hline
	\multicolumn{5}{|c|}{$\mu$} \\ \hline
	0.391752226571890 &  & &   & \\
	0 & 0.368410593050371 & & &   \\
	0 & 0 & 0.251891774271694  &    &\\
	0 & 0 & 0    & 0.544974750228521  &\\
	0 & 0 & 0 & 0.063692468666290 & 0.226007483236906 \\ \hline
	\multicolumn{5}{|c|}{CFL = 1.50818004918983} \\ \hline
	\end{tabular}
	\caption{Butcher Tableau of SSPRK methods}\label{tab:ButcherSSPRK}
	\end{center}
\end{table}
\begin{table}[h!]
	\begin{center}
	\begin{tabular}{ |c|c|c c |}\hline
		\multicolumn{4}{|c|}{Order 2}\\ \hline
		m & $\beta^m$ & \multicolumn{2}{|c|}{$\rho^{m}_z$}\\ \hline
		1 & 1       & $\frac{1}{2}$ & $\frac{1}{2}$ \\ \hline
	\end{tabular}
	\quad
	\begin{tabular}{ |c|c|c c c |}\hline
		\multicolumn{5}{|c|}{Order 3}\\ \hline
		m & $\beta^m$ & \multicolumn{3}{|c|}{$\rho^{m}_z$}\\ \hline 
		1 & $\frac{1}{2}$      & $\frac{5}{24}$ & $\frac{1}{3}$  & $-\frac{1}{24}$ \\ 
		2 & 1         & $\frac{1}{6}$ & $\frac{2}{3}$  & $\frac{1}{3}$\\ \hline
	\end{tabular}
	\quad
	\begin{tabular}{ |c|c|c c c c|}\hline
		\multicolumn{6}{|c|}{Order 4}\\ \hline
		m & $\beta^m$ & \multicolumn{4}{|c|}{$\rho^{m}_z$}\\ \hline 
		1 & $\frac{1}{3}$      & $\frac{1}{8}$ & $\frac{19}{72}$  & $-\frac{5}{72}$  & $\frac{1}{72}$ \\ 
		2 & $\frac{2}{3}$       & $\frac{1}{9}$ & $\frac{4}{9}$  & $\frac{1}{9}$  & 0\\ 
		3 & 1         & $\frac{1}{8}$ & $\frac{3}{8}$  & $\frac{3}{8}$ & $\frac{1}{8}$\\ \hline
	\end{tabular}
	\caption{DeC coefficients for equispaced subtimesteps. }\label{tab:DeCcoeff}
	\end{center}
\end{table}
\endgroup

\section{Fourier analysis, spatial and temporal eigenanalysis}\label{app:fourier_full_disc}
\begin{figure}[h!]
	\centering
	\title{Without any stabilization method} \\
	\includegraphics[width=0.24\textwidth]{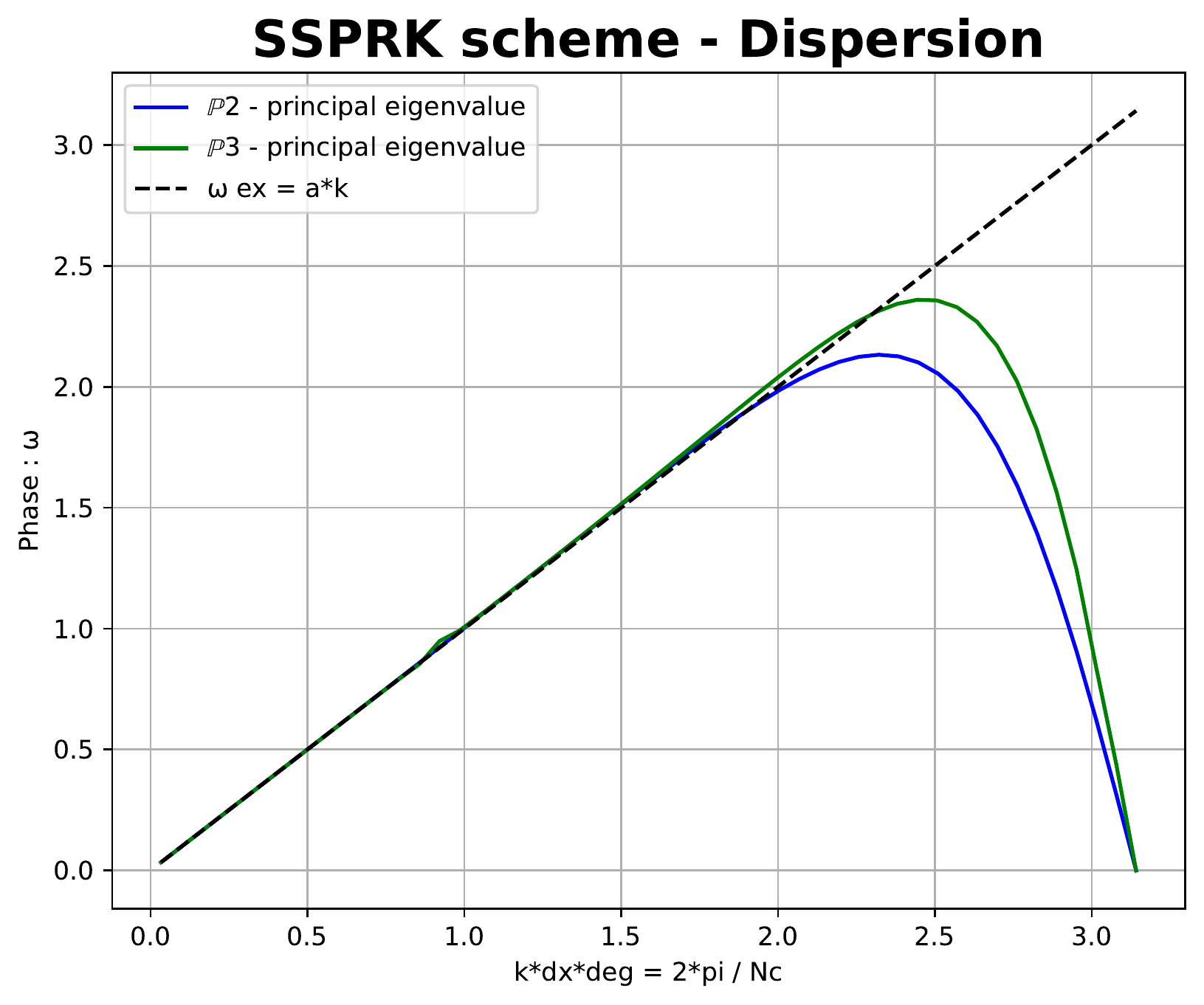}
	\includegraphics[width=0.24\textwidth]{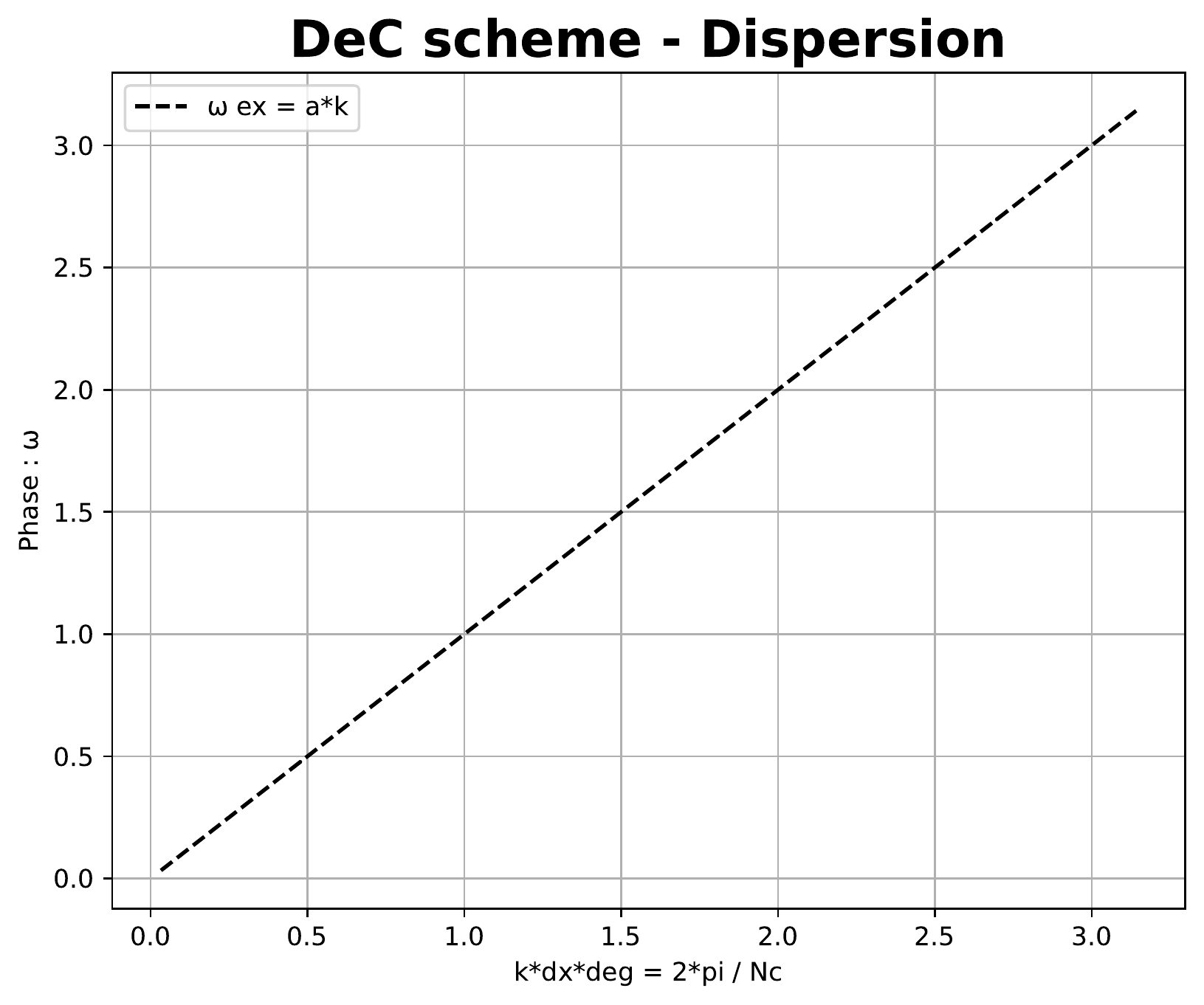}
	\includegraphics[width=0.24\textwidth]{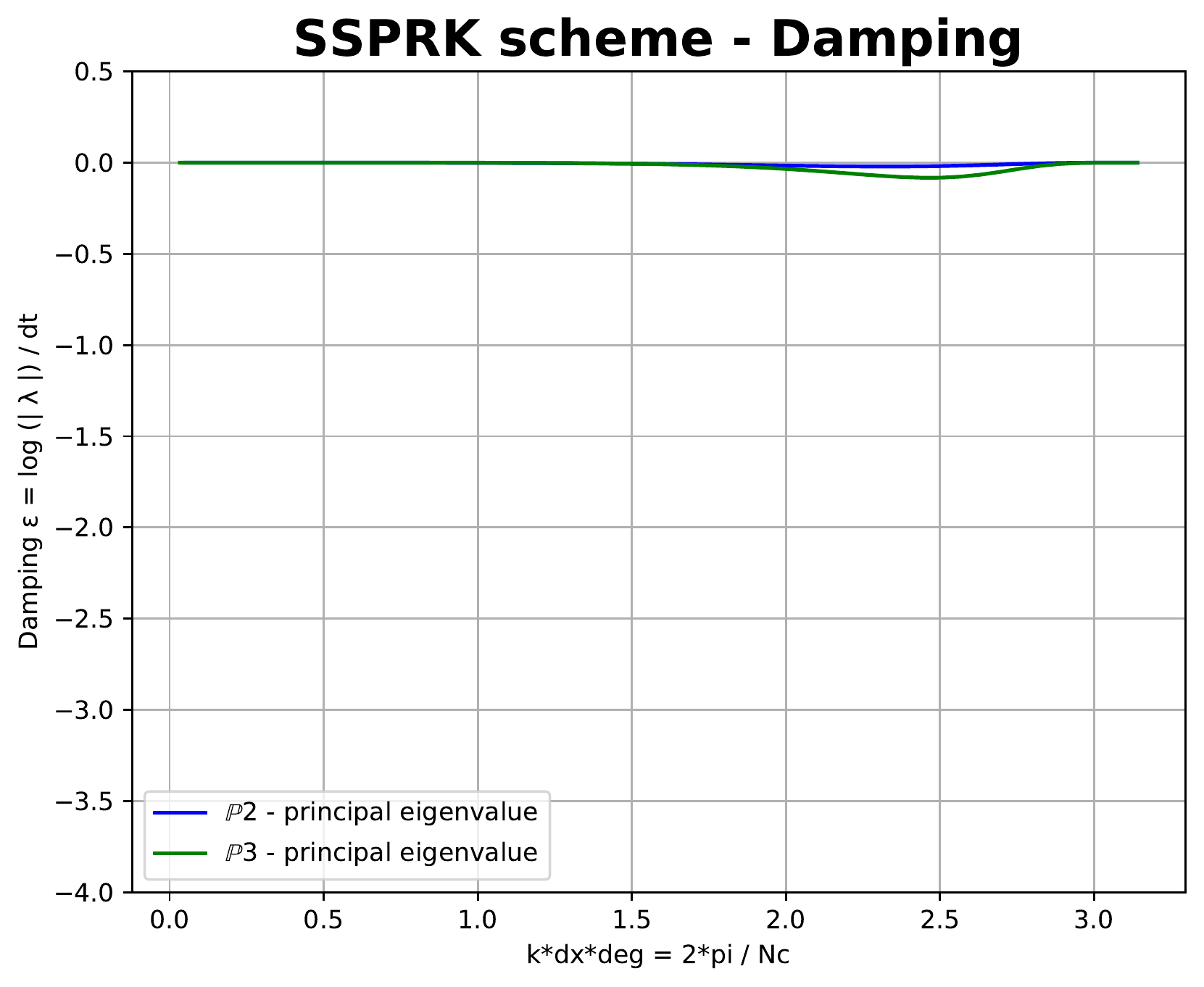}
	\includegraphics[width=0.24\textwidth]{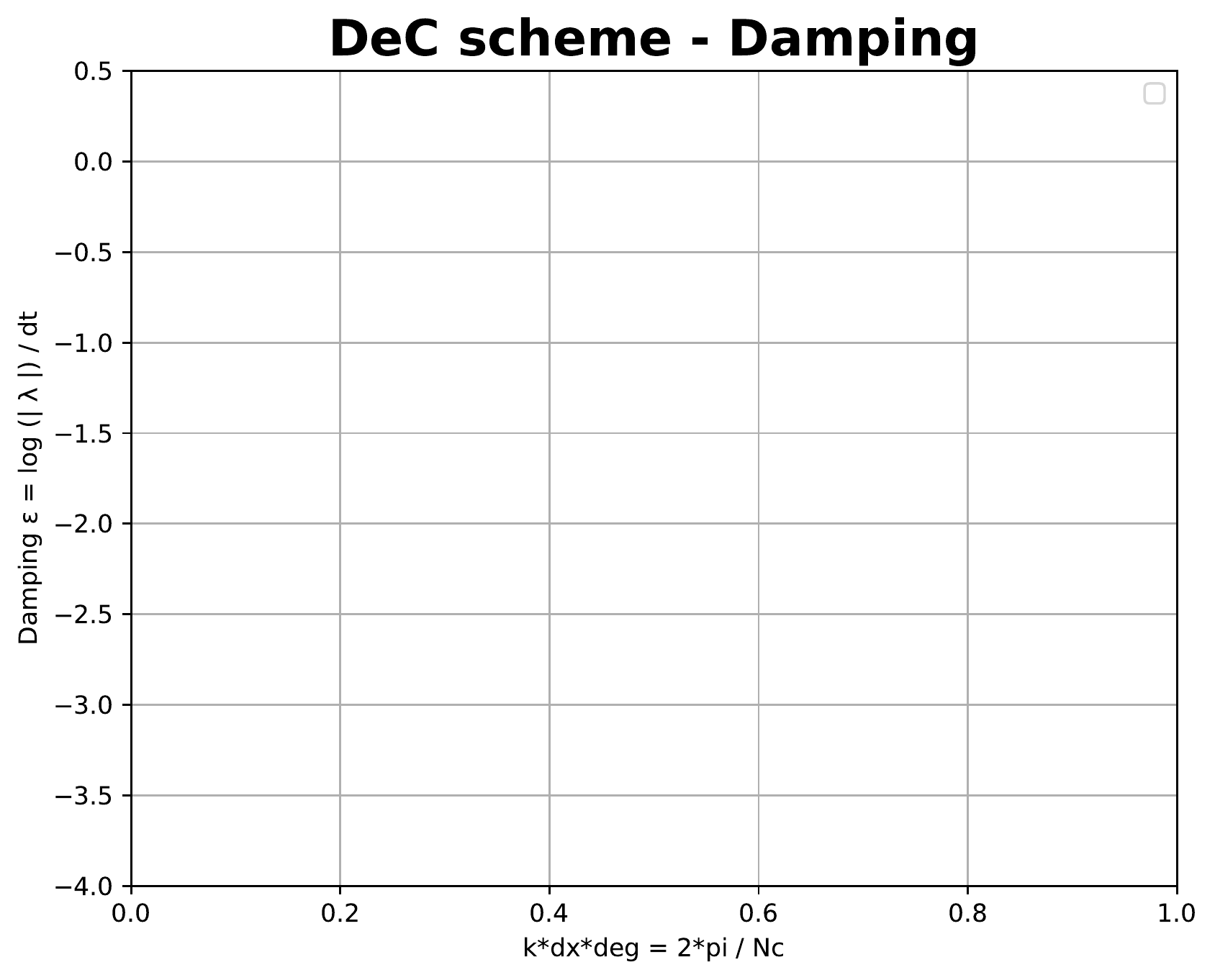}	
	\\
	\title{Using the SUPG stabilization method} \\
	\includegraphics[width=0.24\textwidth]{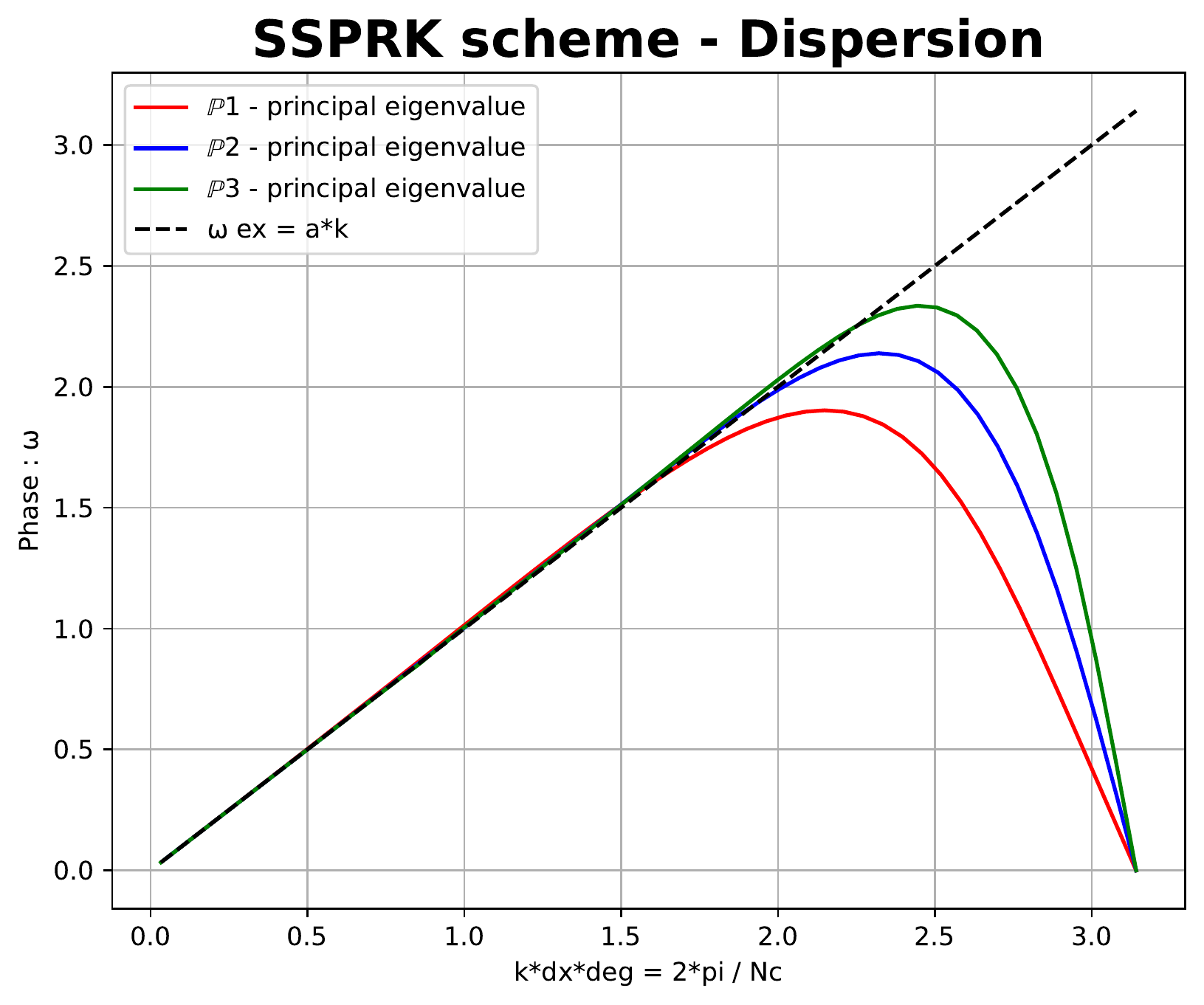}
	\includegraphics[width=0.24\textwidth]{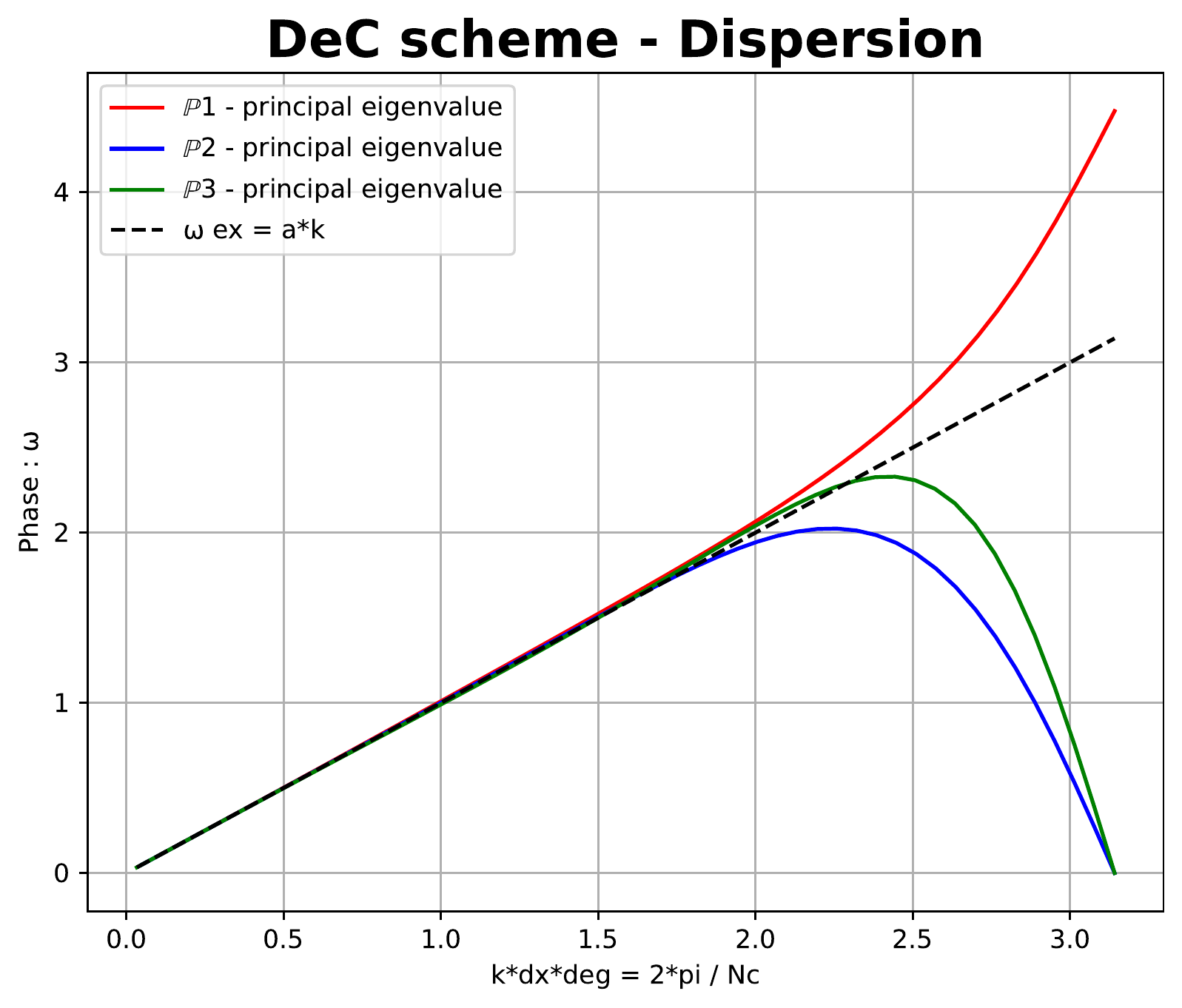}
	\includegraphics[width=0.24\textwidth]{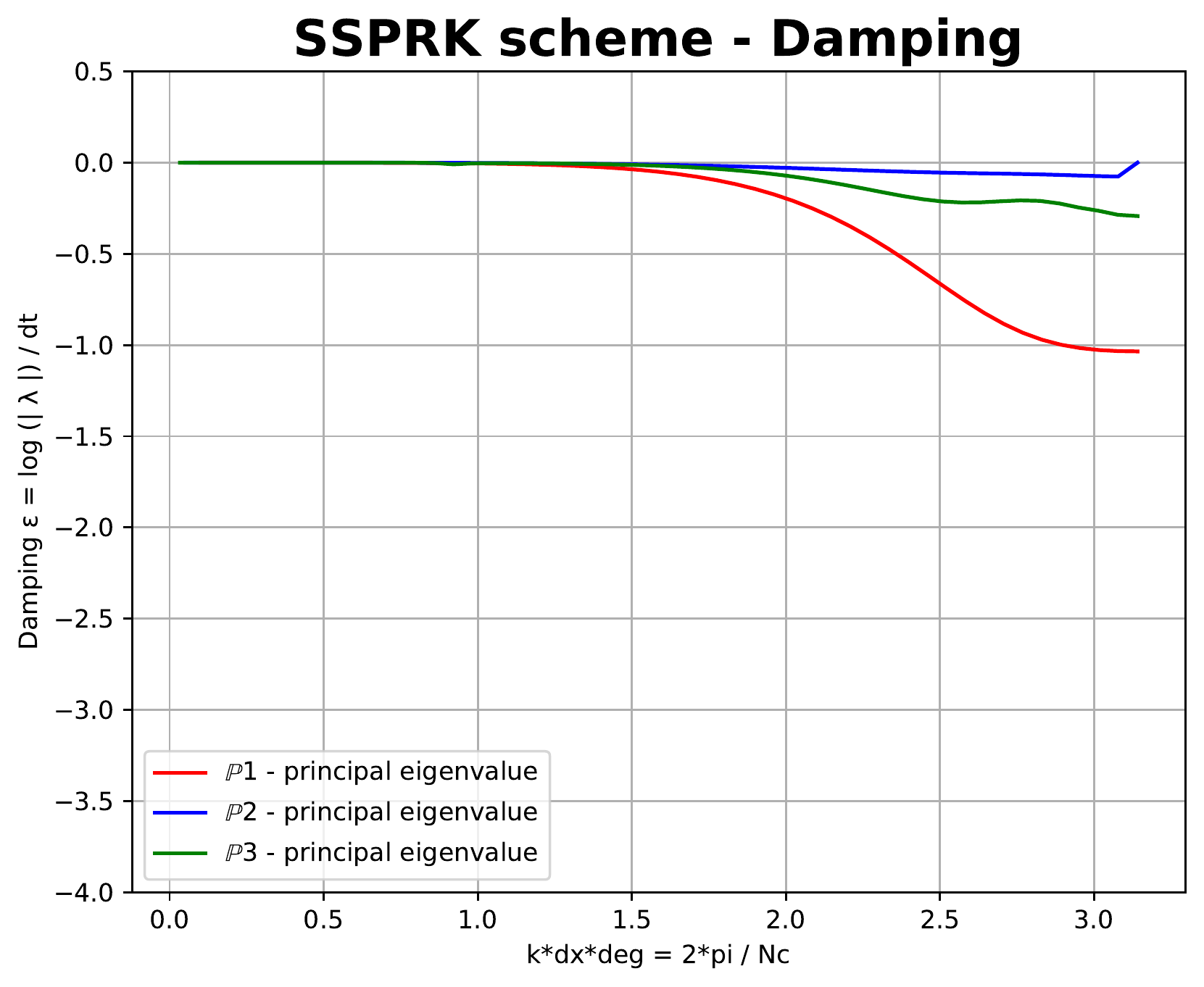}
	\includegraphics[width=0.24\textwidth]{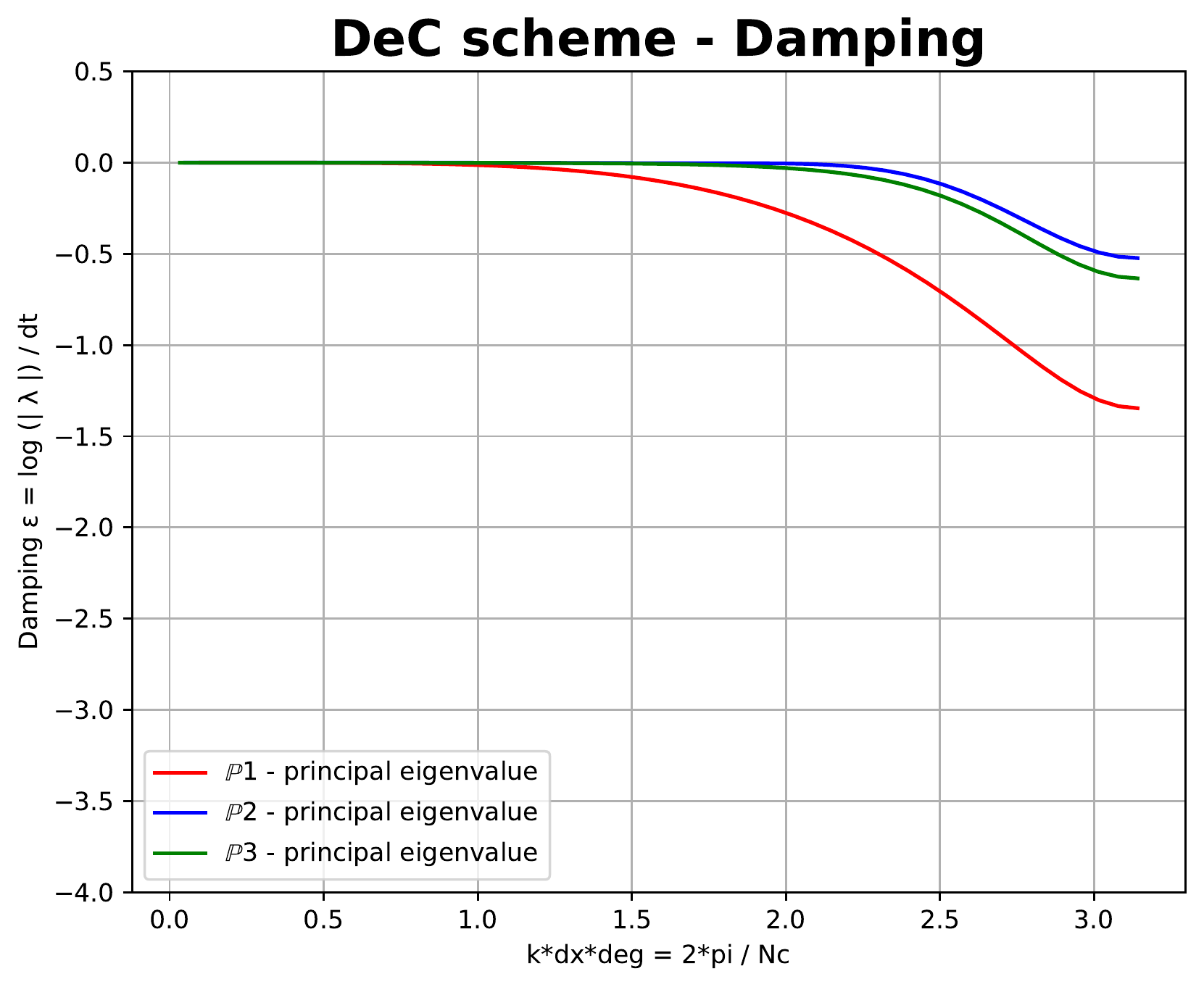}	
	\\
	\title{Using the LPS stabilization method} \\
	\includegraphics[width=0.24\textwidth]{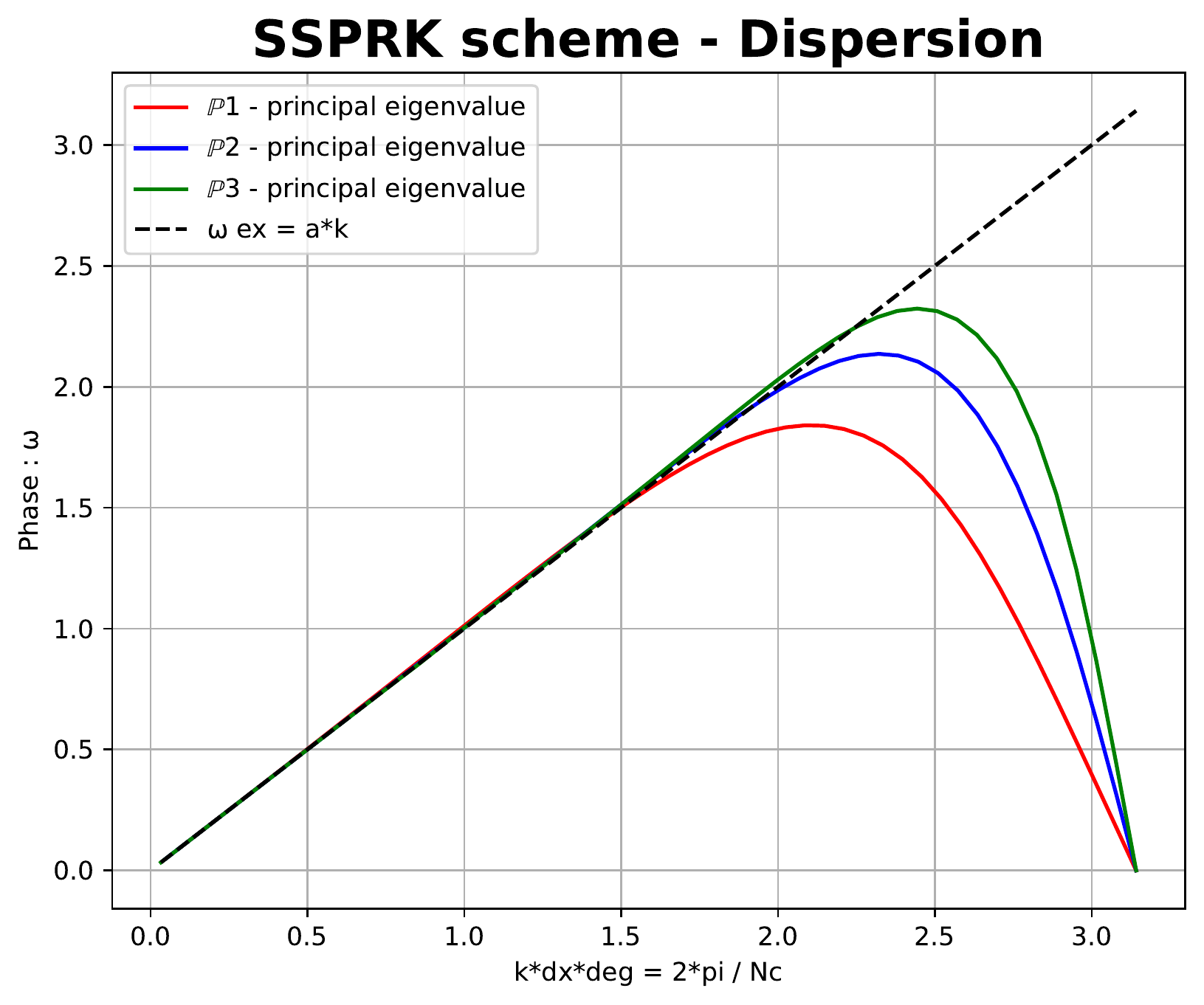}
	\includegraphics[width=0.24\textwidth]{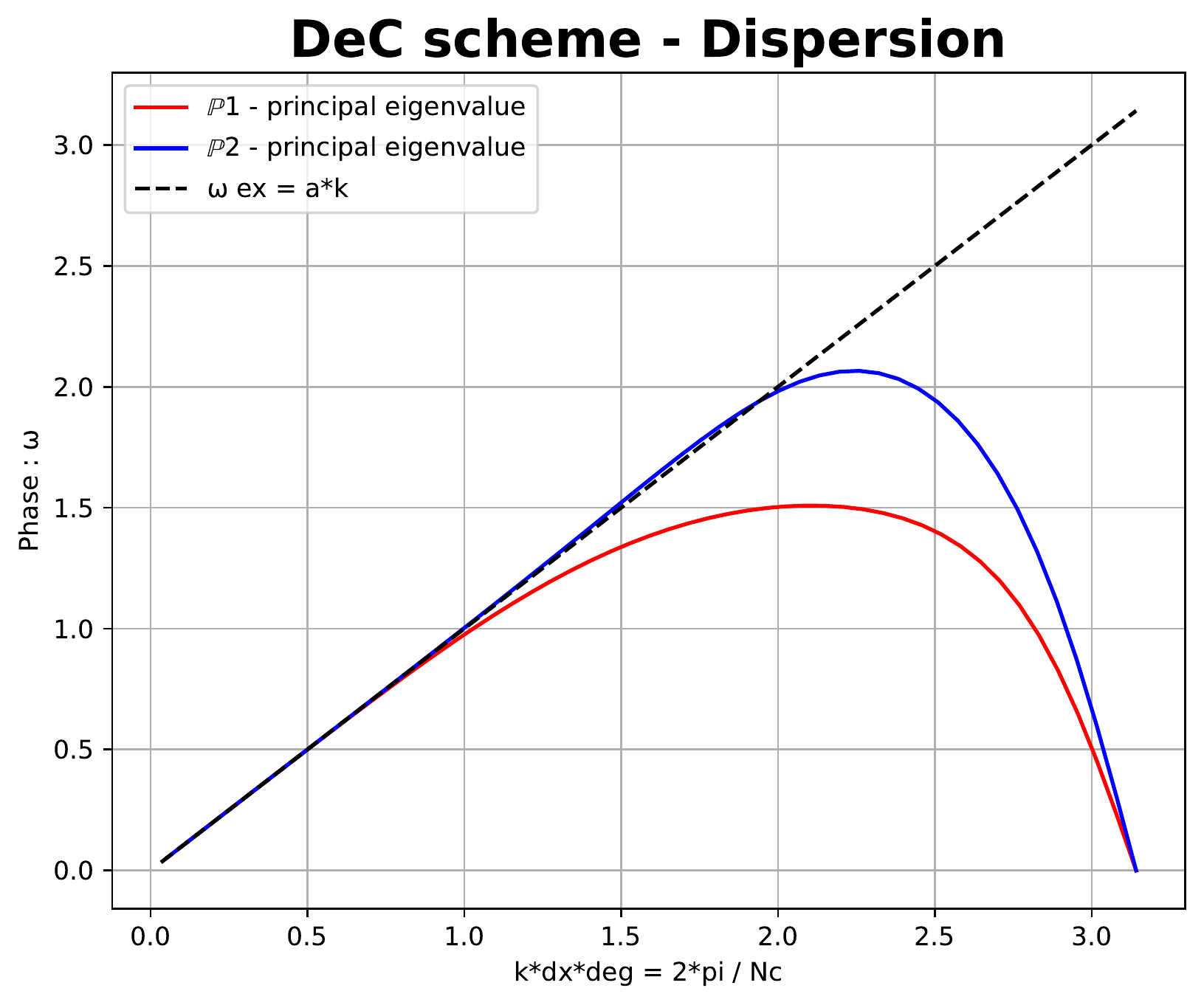}
	\includegraphics[width=0.24\textwidth]{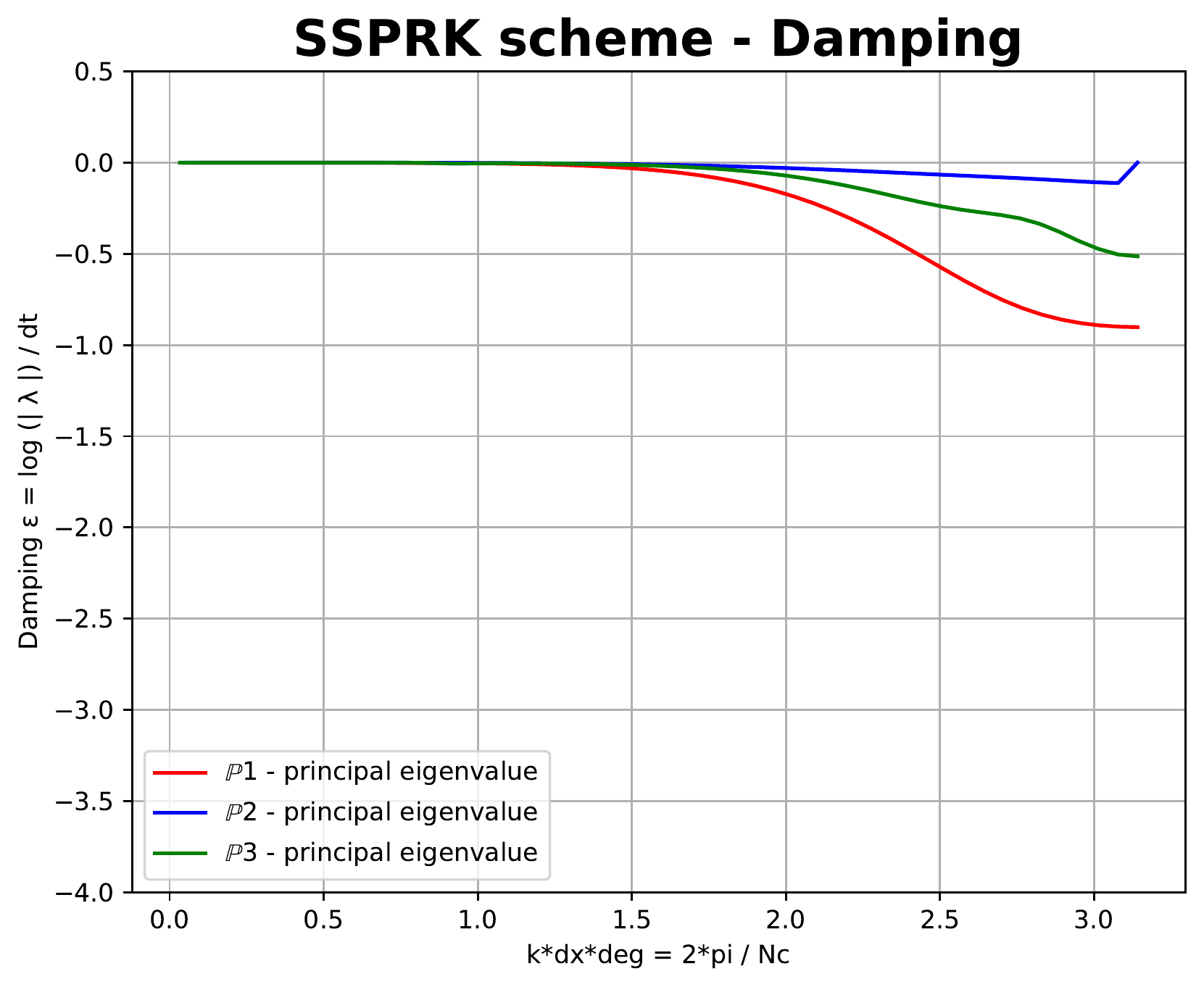}
	\includegraphics[width=0.24\textwidth]{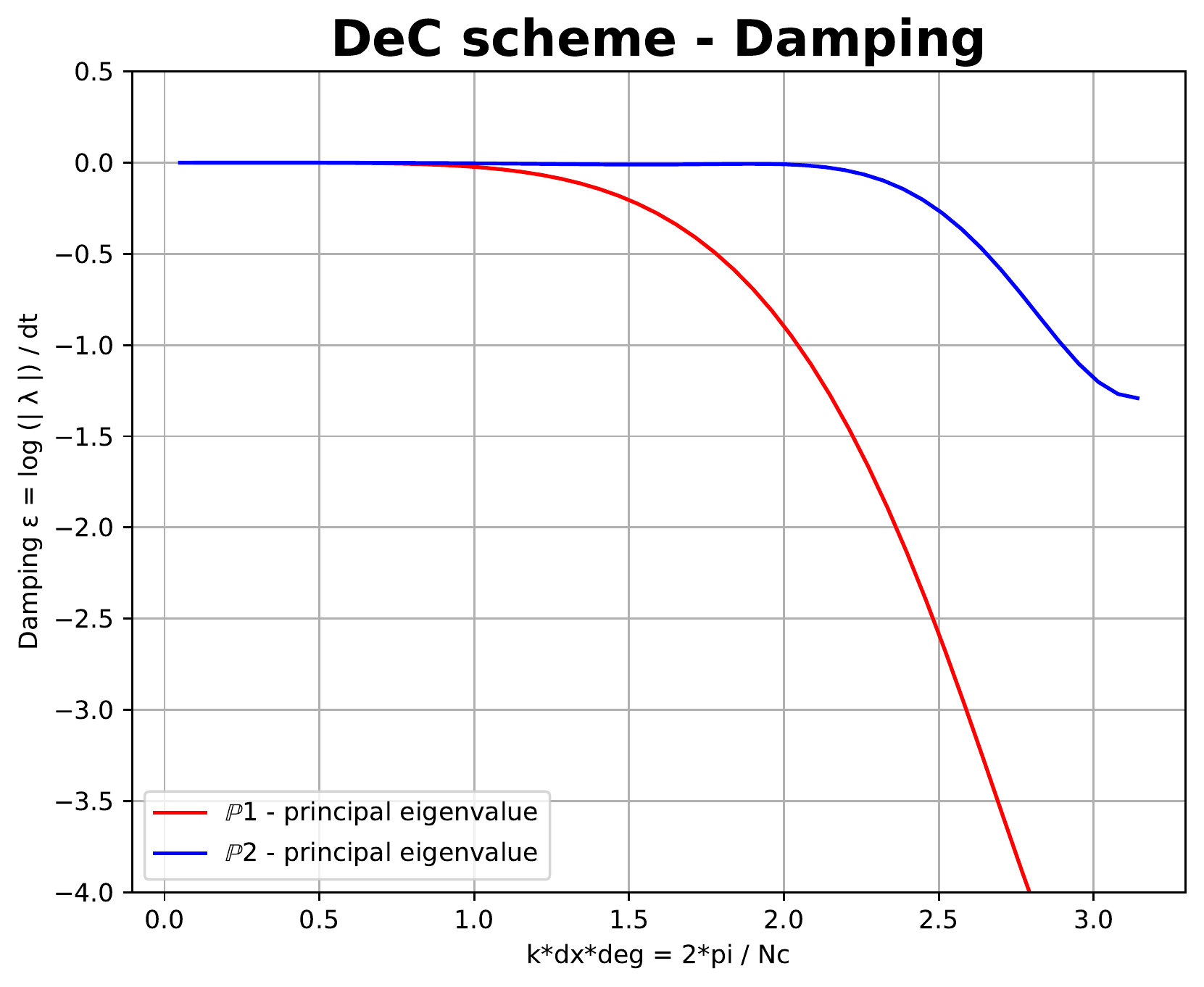}	
	\\
	\title{Using the CIP stabilization method} \\
	\includegraphics[width=0.24\textwidth]{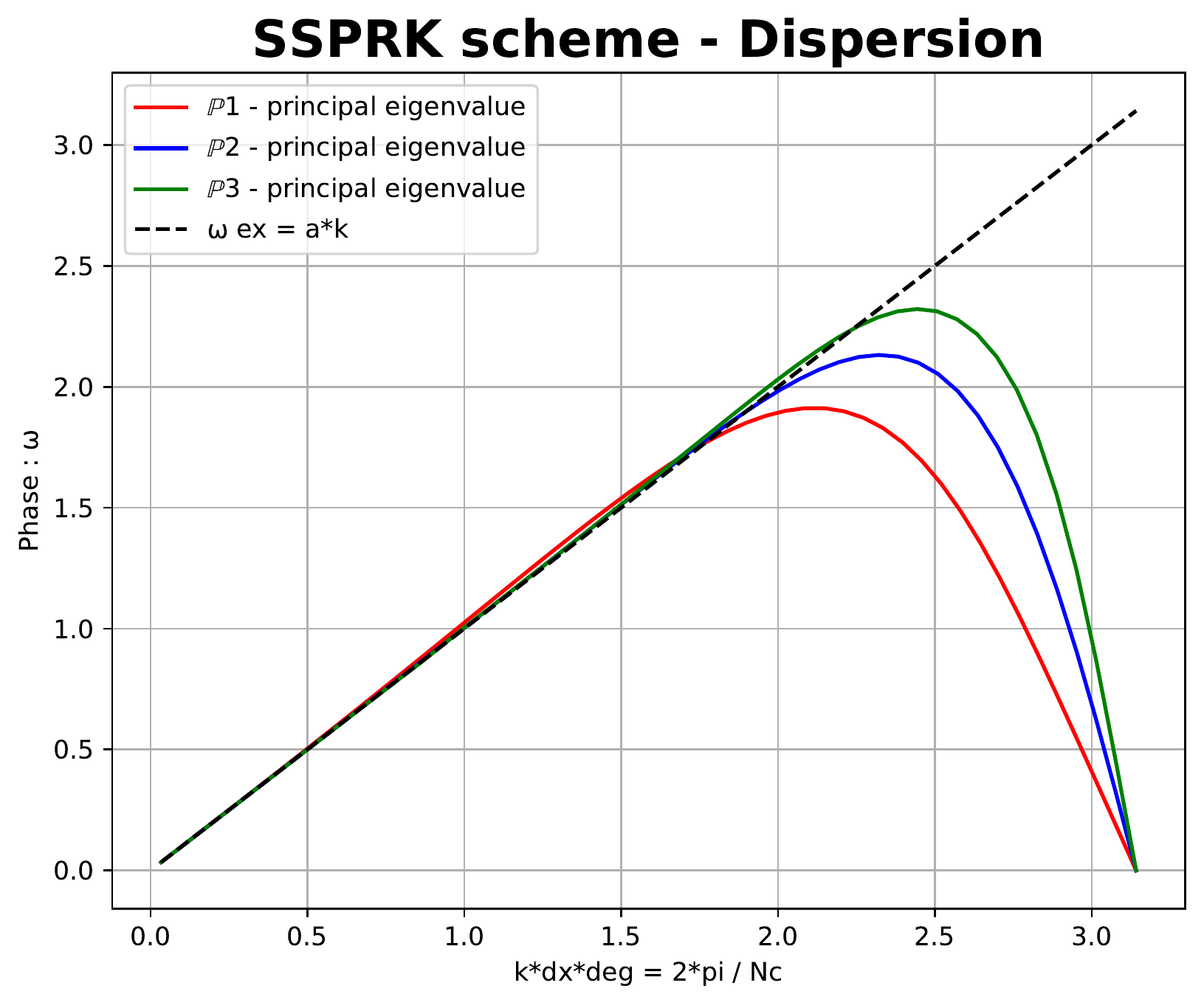}
	\includegraphics[width=0.24\textwidth]{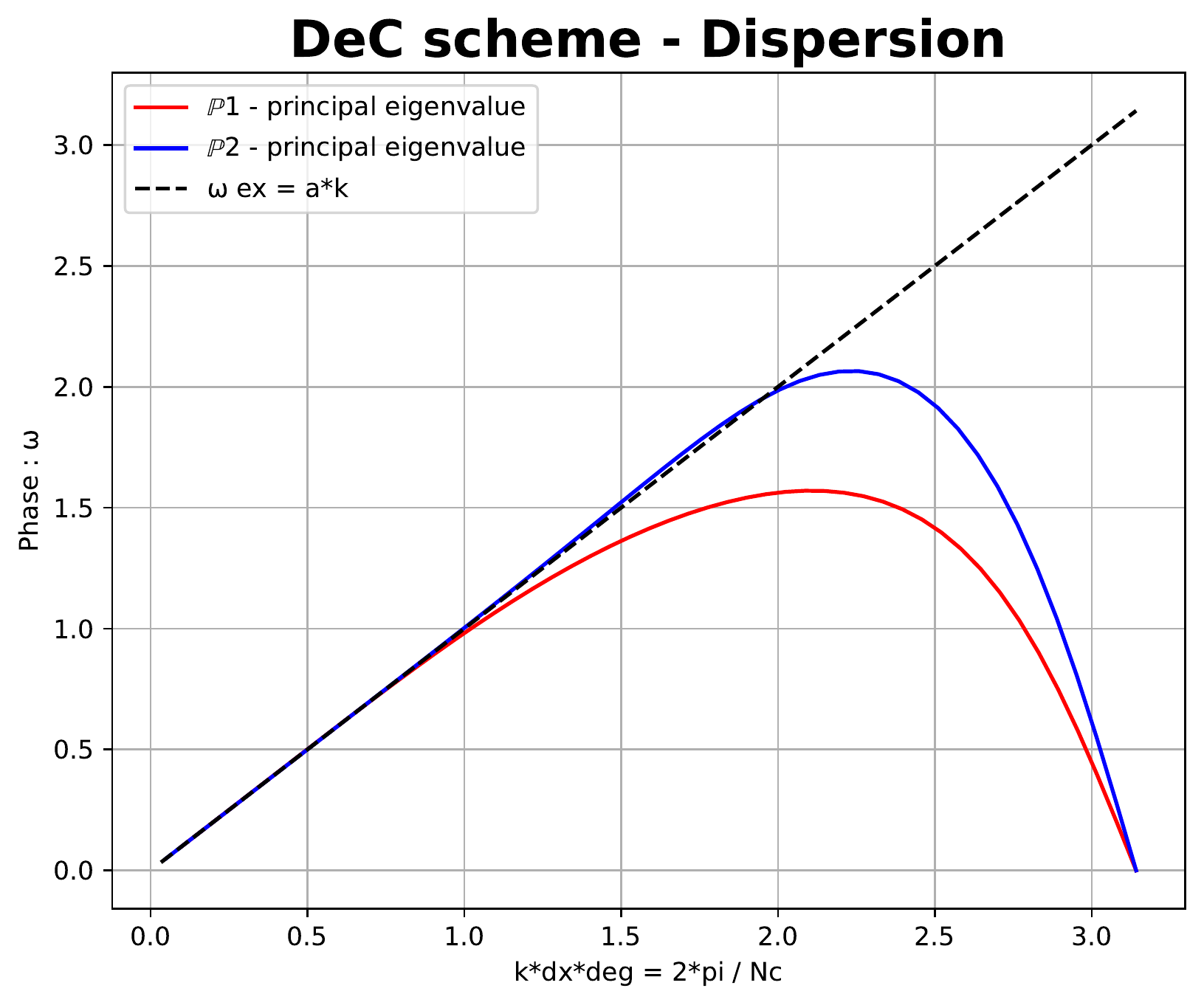}
	\includegraphics[width=0.24\textwidth]{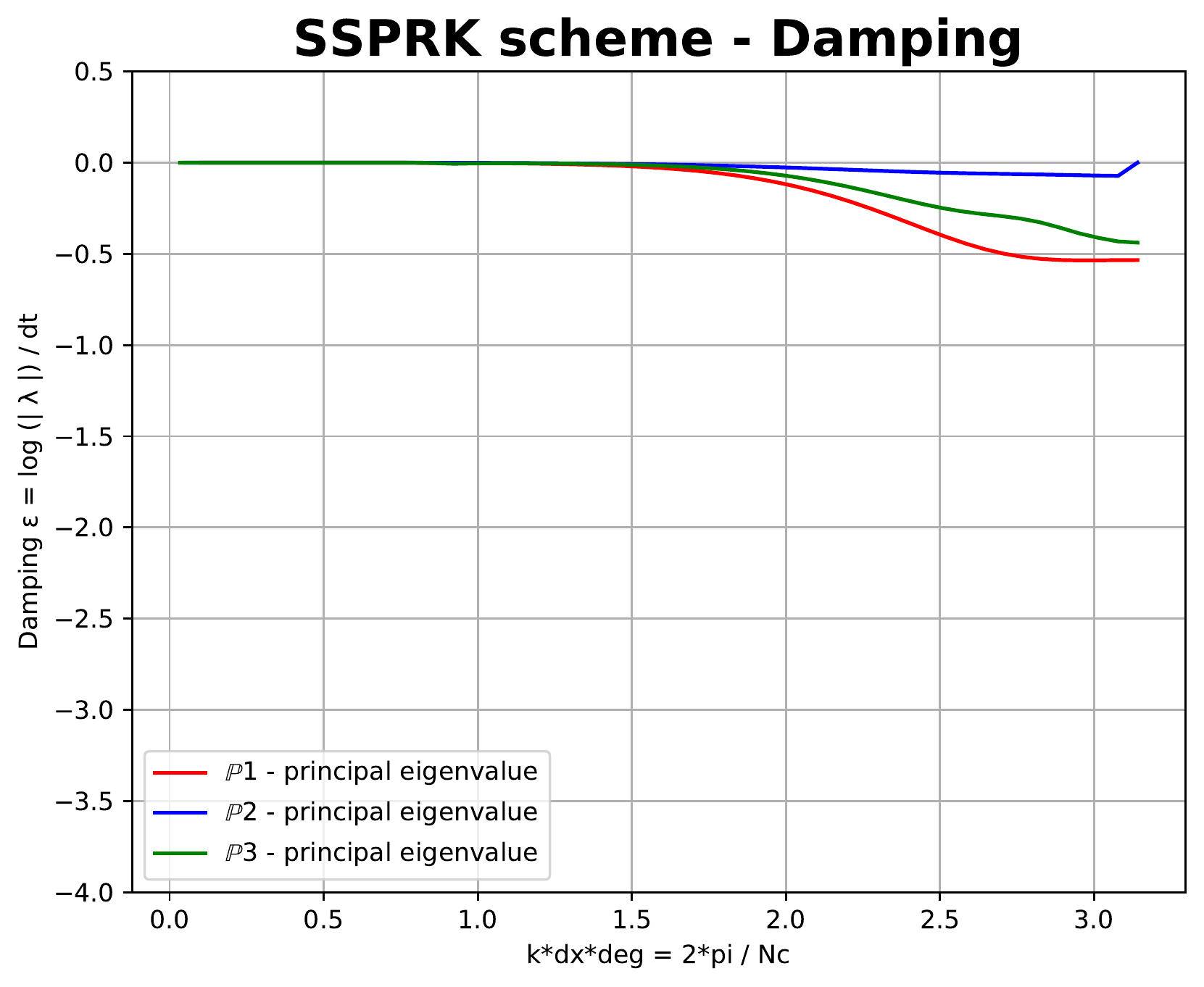}
	\includegraphics[width=0.24\textwidth]{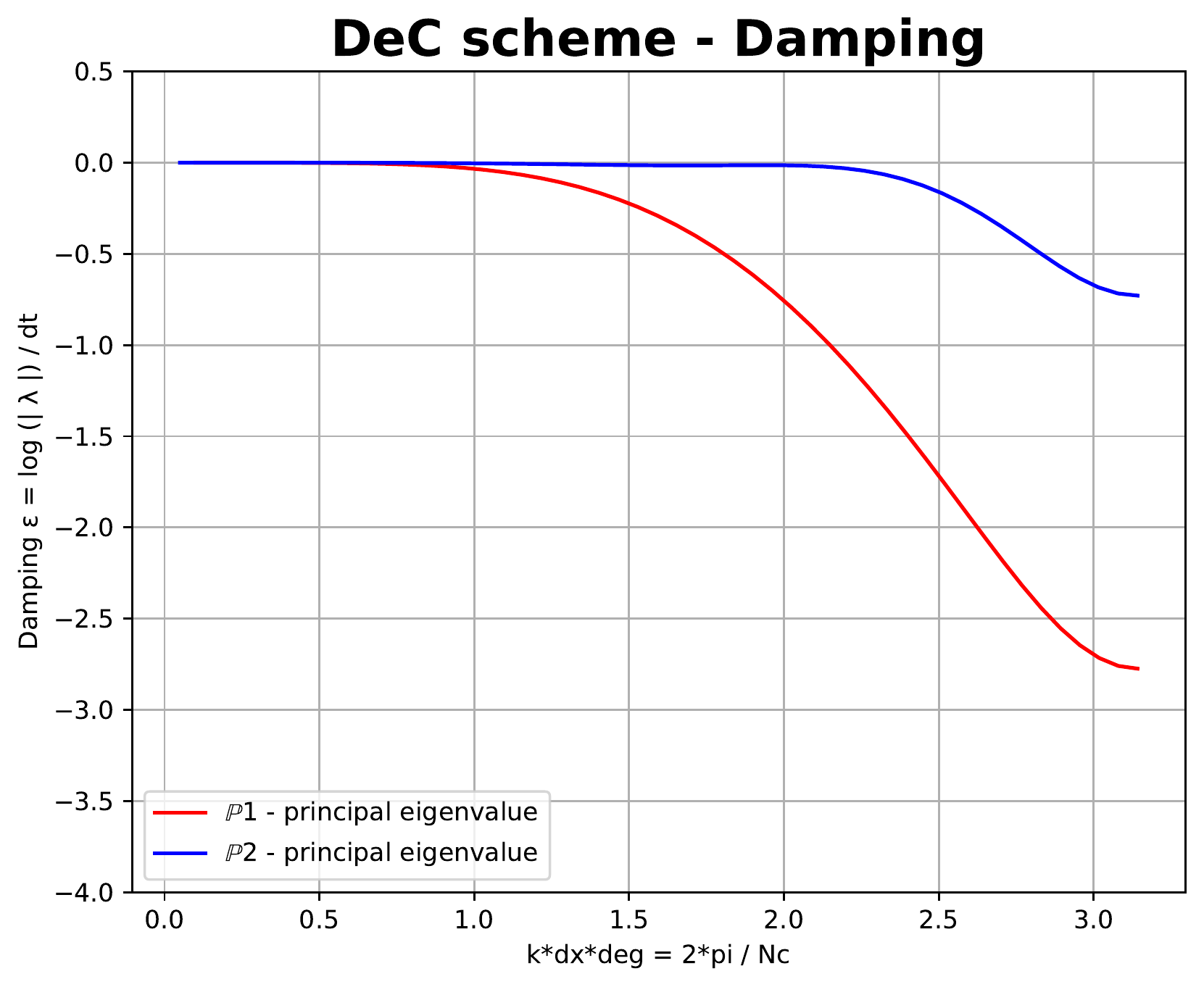}	
	\caption{Dispersion and damping coefficients for \textit{basic} elements, with DeC and SSPRK methods and all stabilization techniques}\label{fig:dispersionBasic}
\end{figure}

\begin{figure}[h!]
	\centering
	\title{Without any stabilization method} \\
	\includegraphics[width=0.24\textwidth]{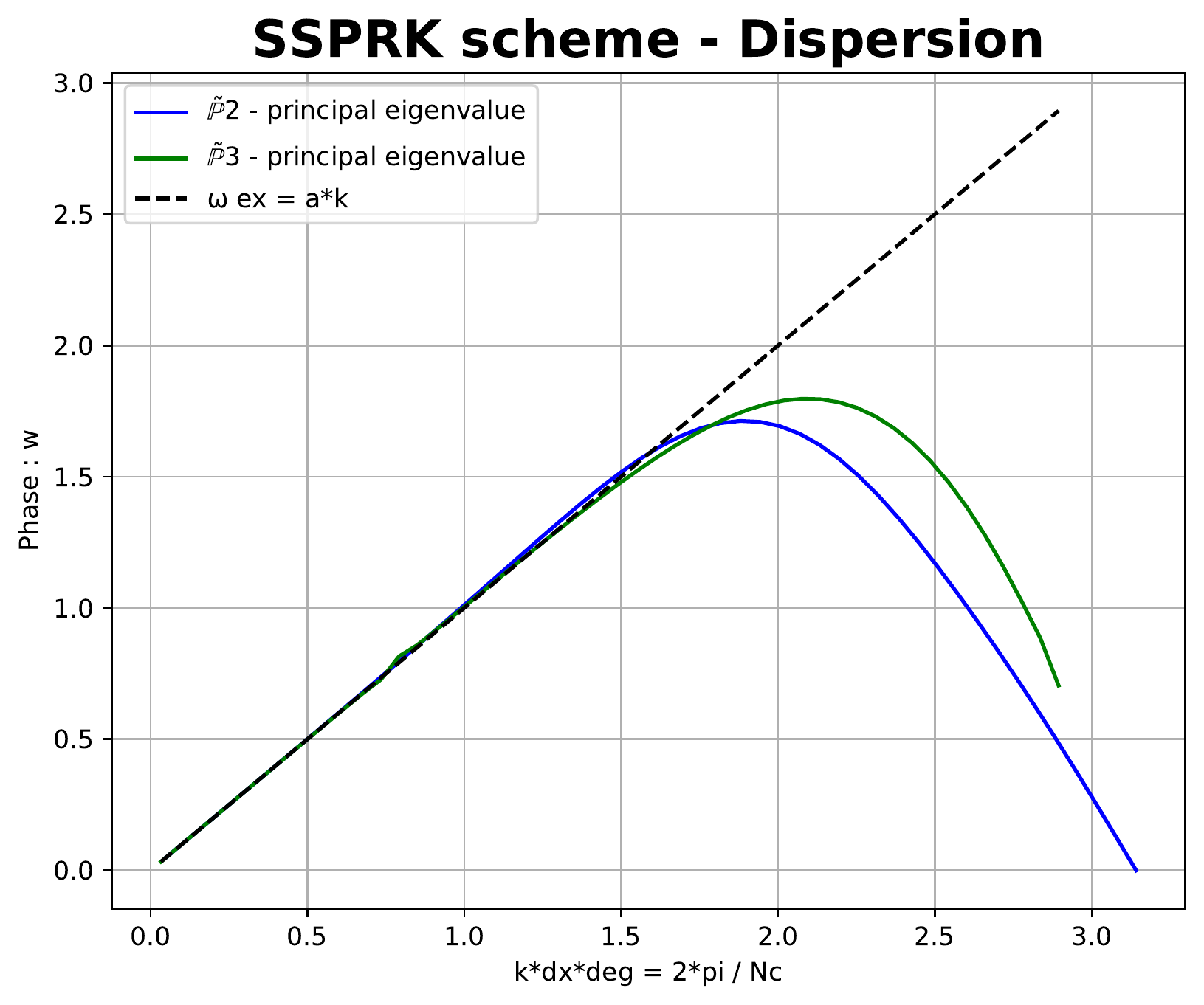}
	\includegraphics[width=0.24\textwidth]{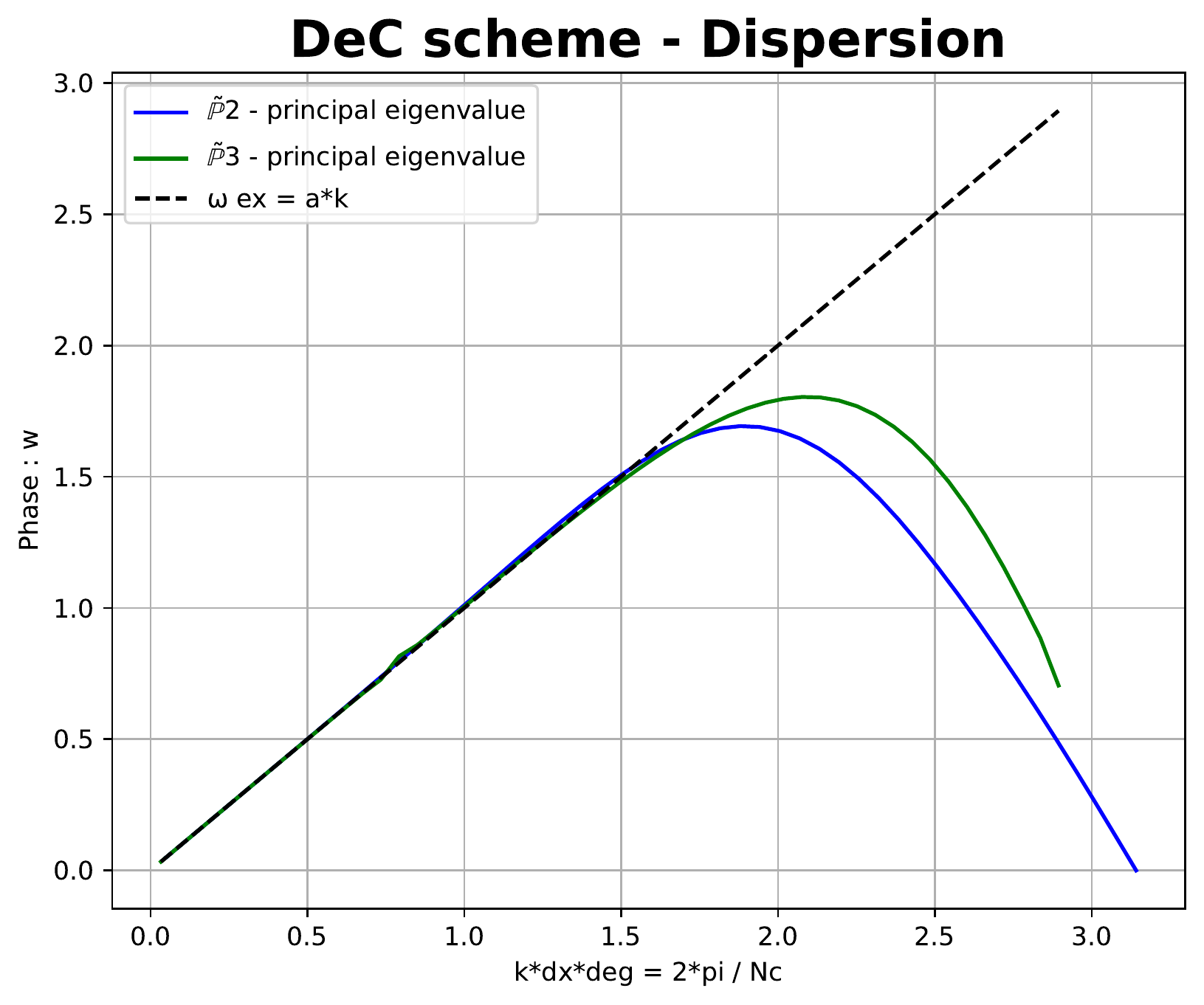}
	\includegraphics[width=0.24\textwidth]{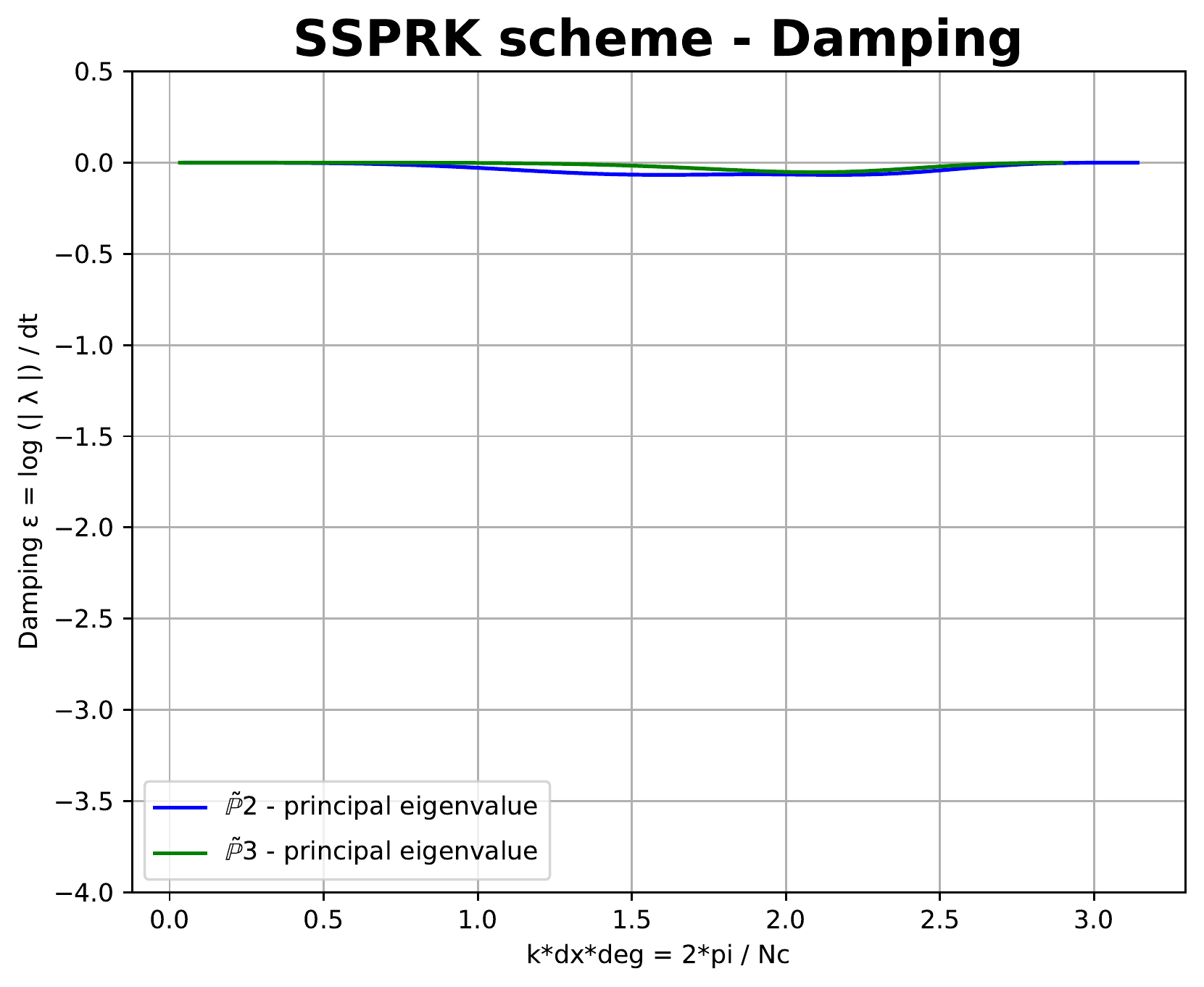}
	\includegraphics[width=0.24\textwidth]{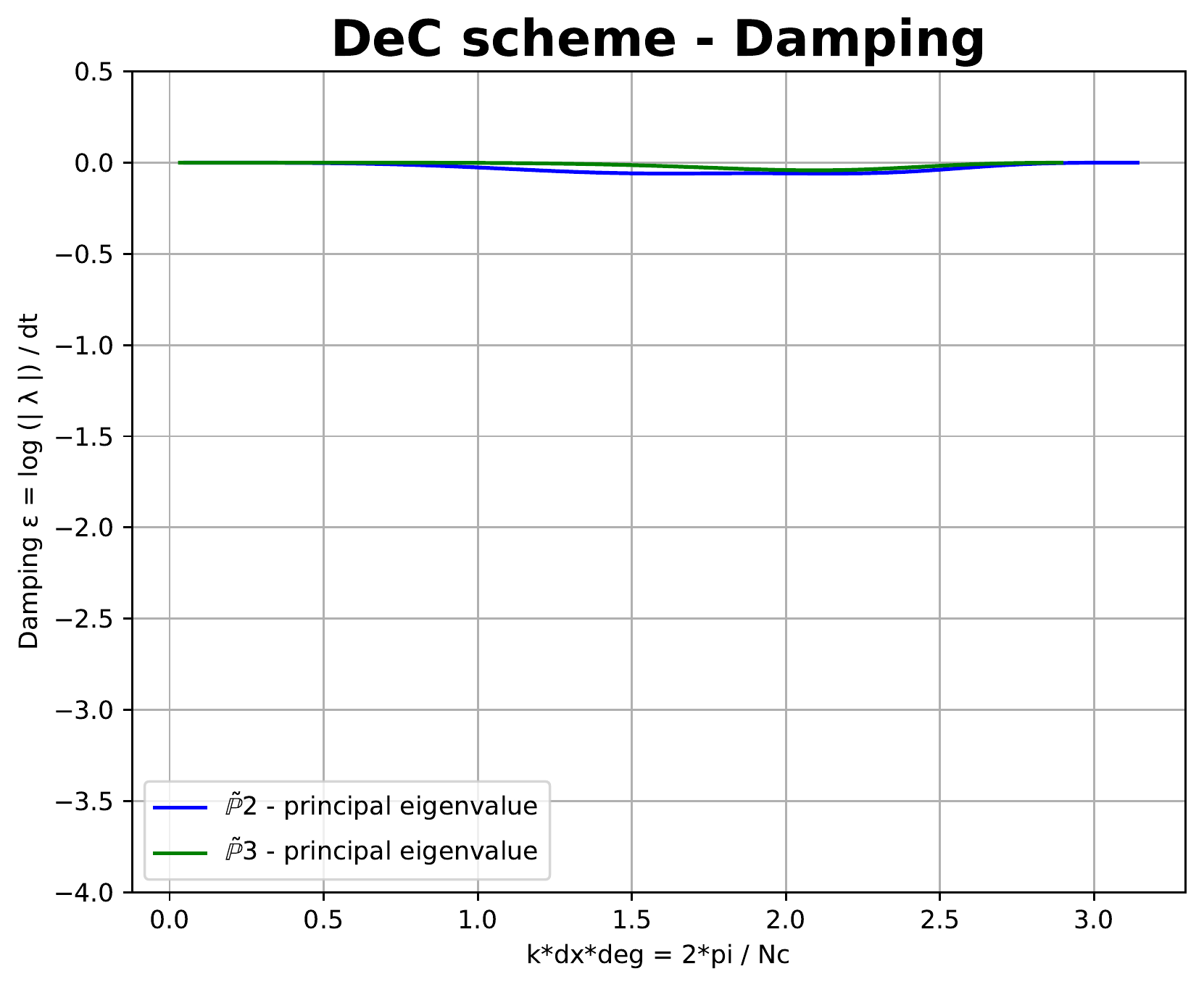}	
\\
	\centering
	\title{Using the SUPG stabilization method} \\
	\includegraphics[width=0.24\textwidth]{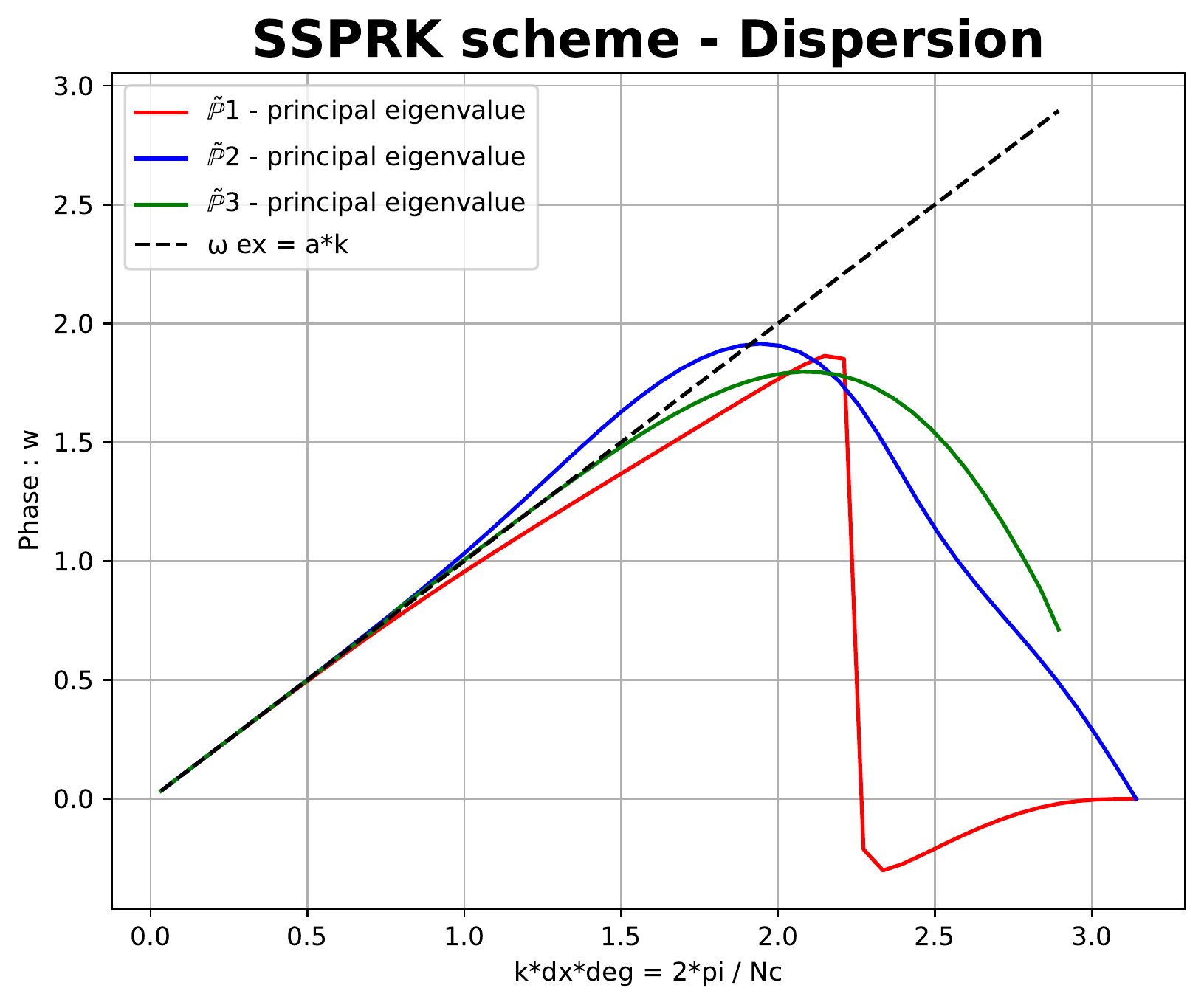}
	\includegraphics[width=0.24\textwidth]{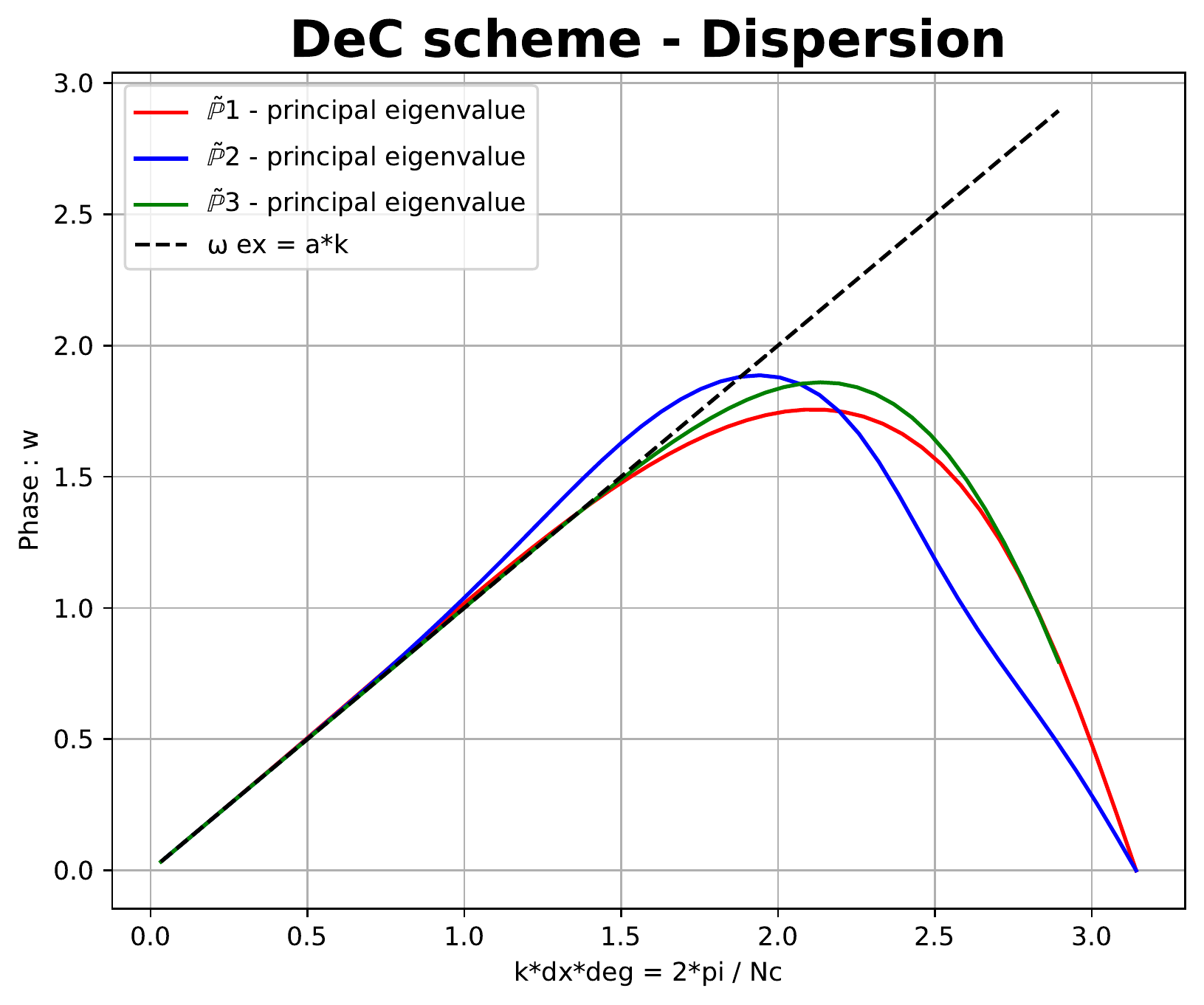}
	\includegraphics[width=0.24\textwidth]{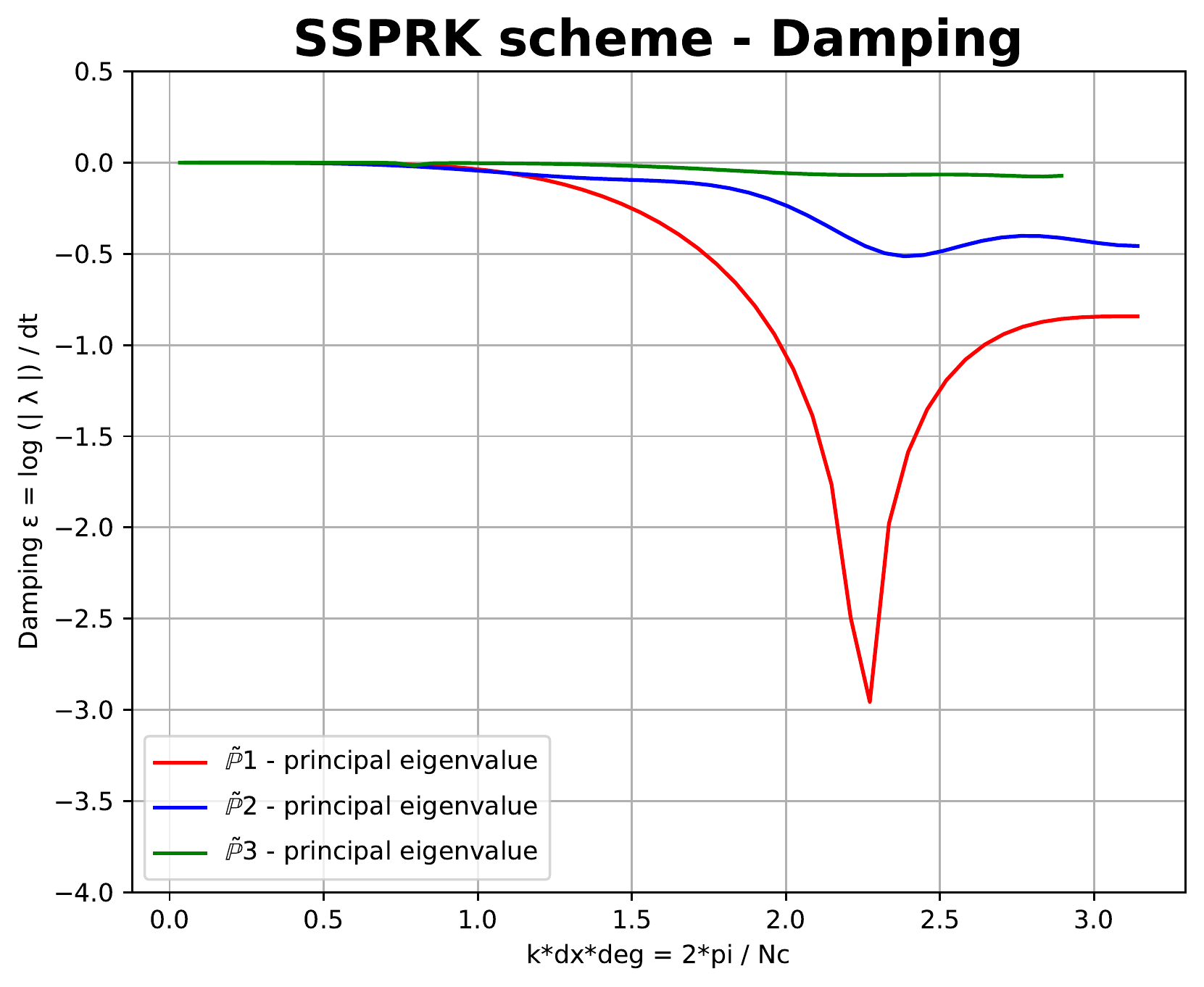}
	\includegraphics[width=0.24\textwidth]{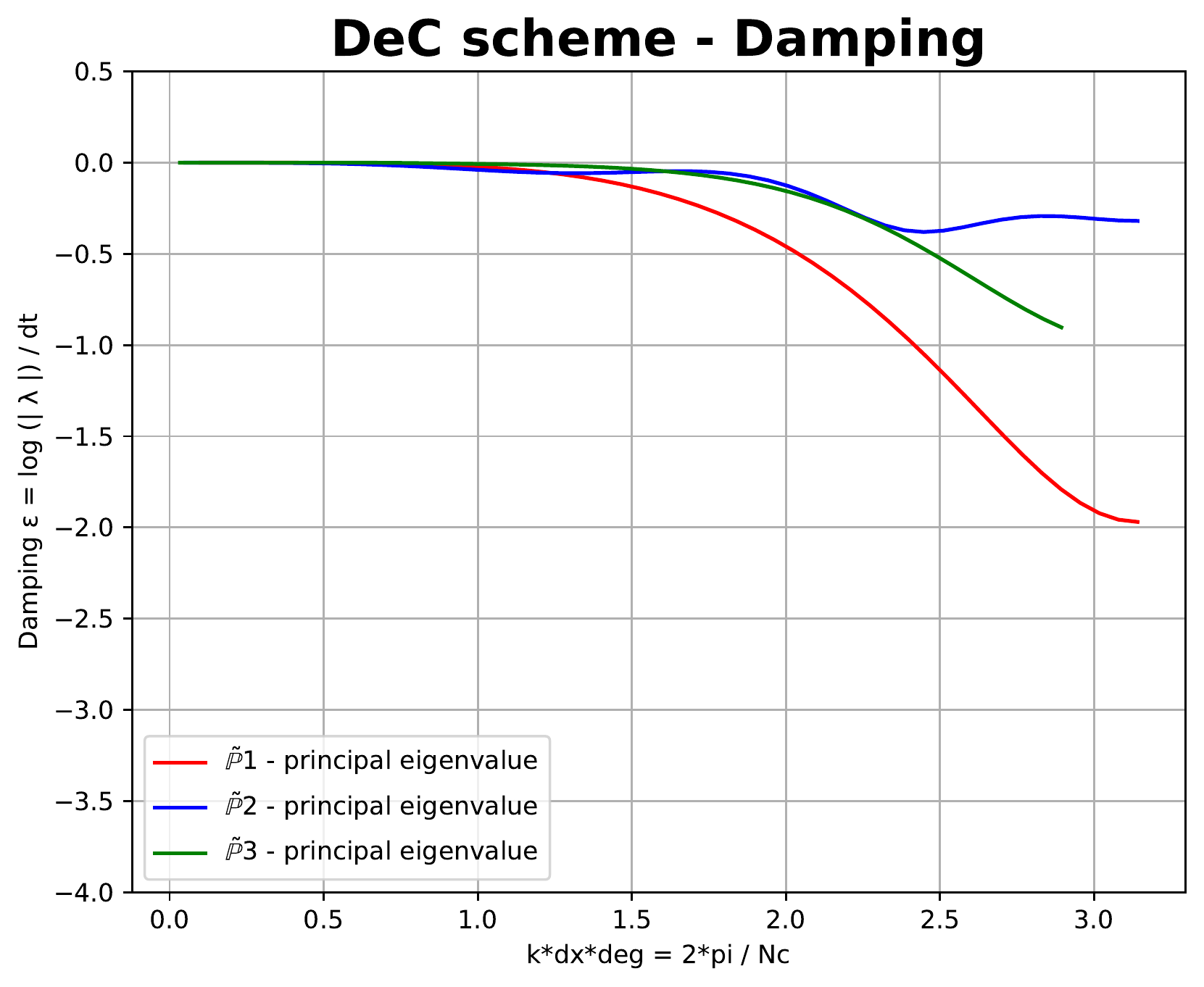}	
\\
	\centering
	\title{Using the LPS stabilization method} \\
	\includegraphics[width=0.24\textwidth]{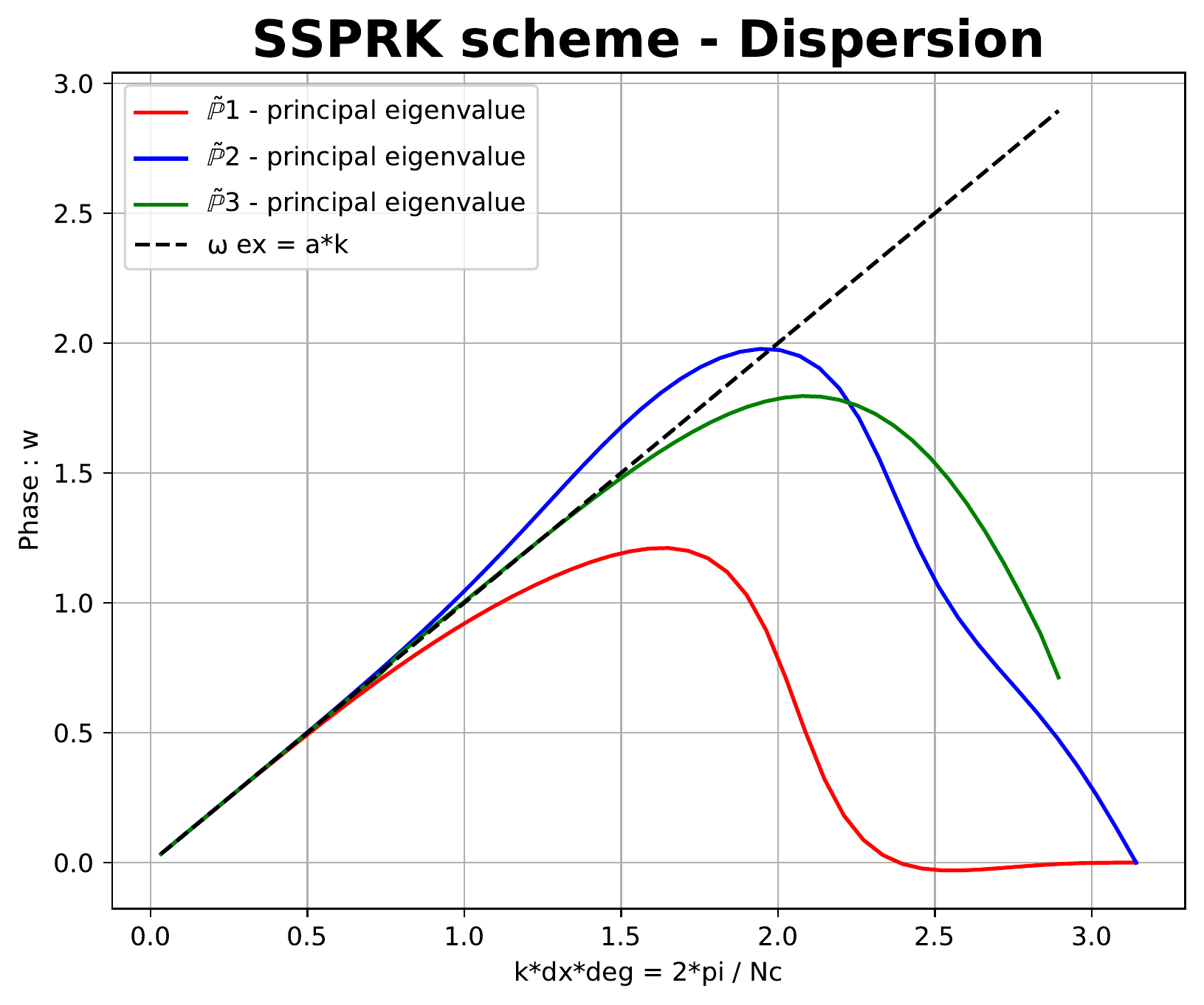}
	\includegraphics[width=0.24\textwidth]{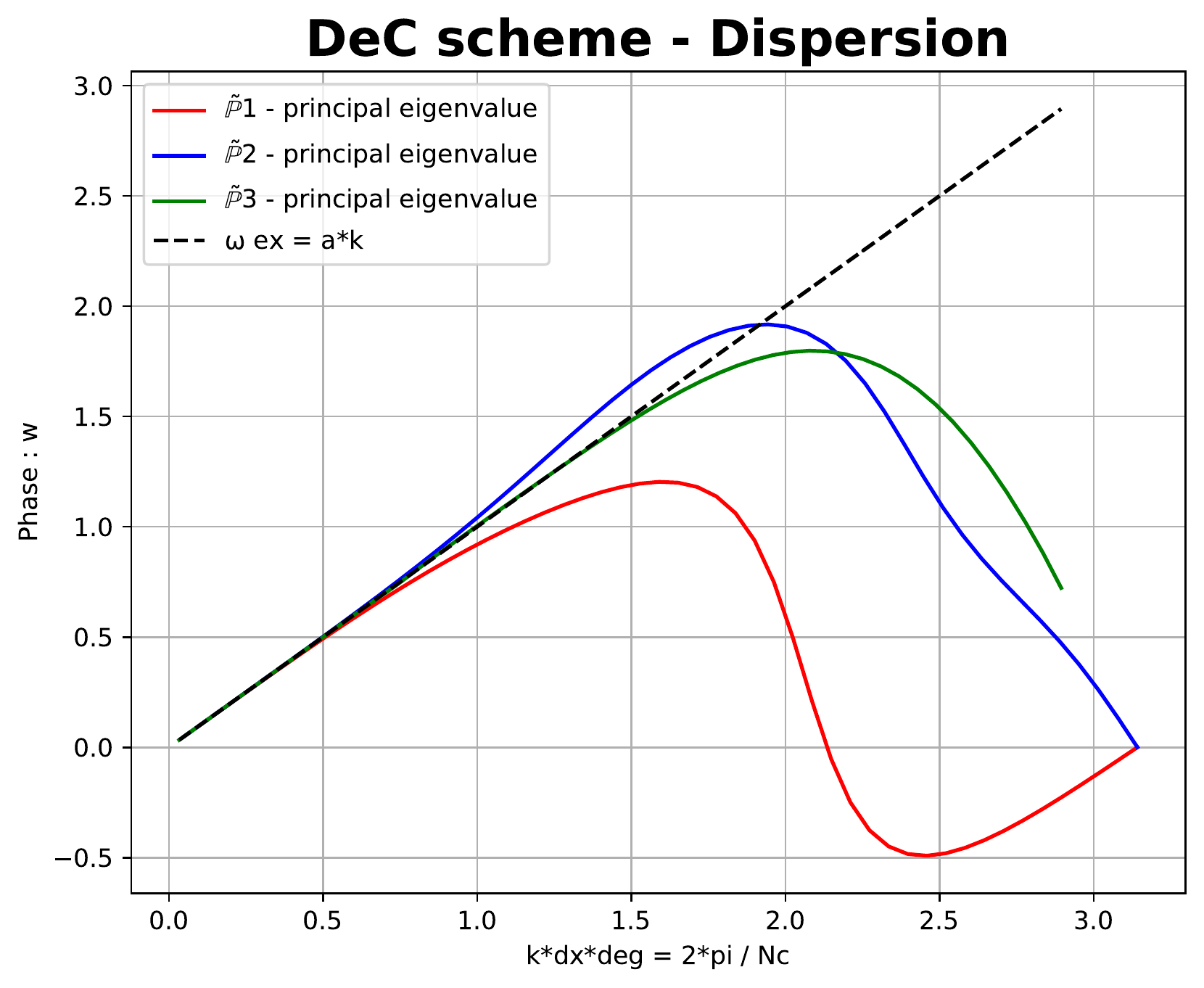}
	\includegraphics[width=0.24\textwidth]{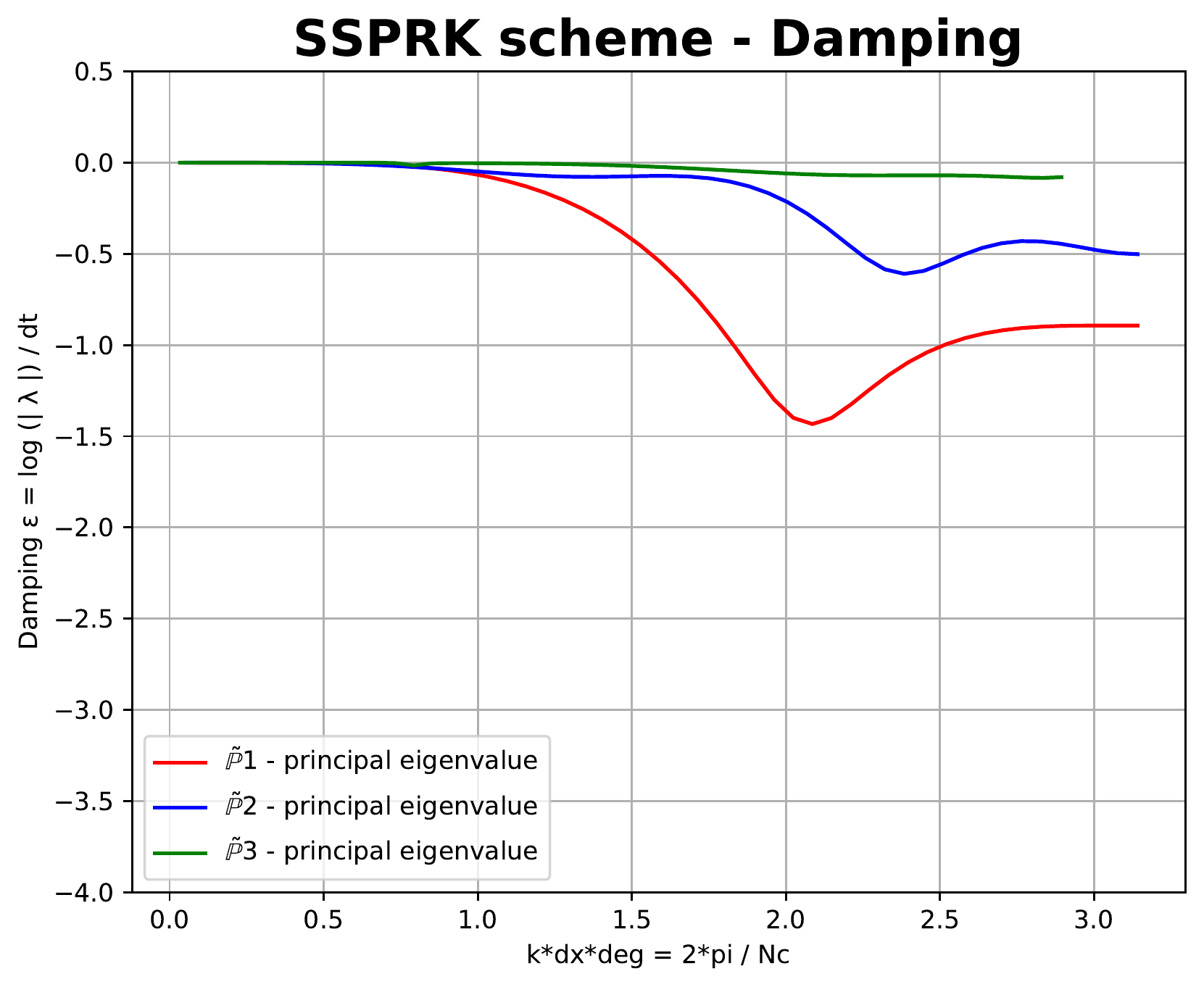}
	\includegraphics[width=0.24\textwidth]{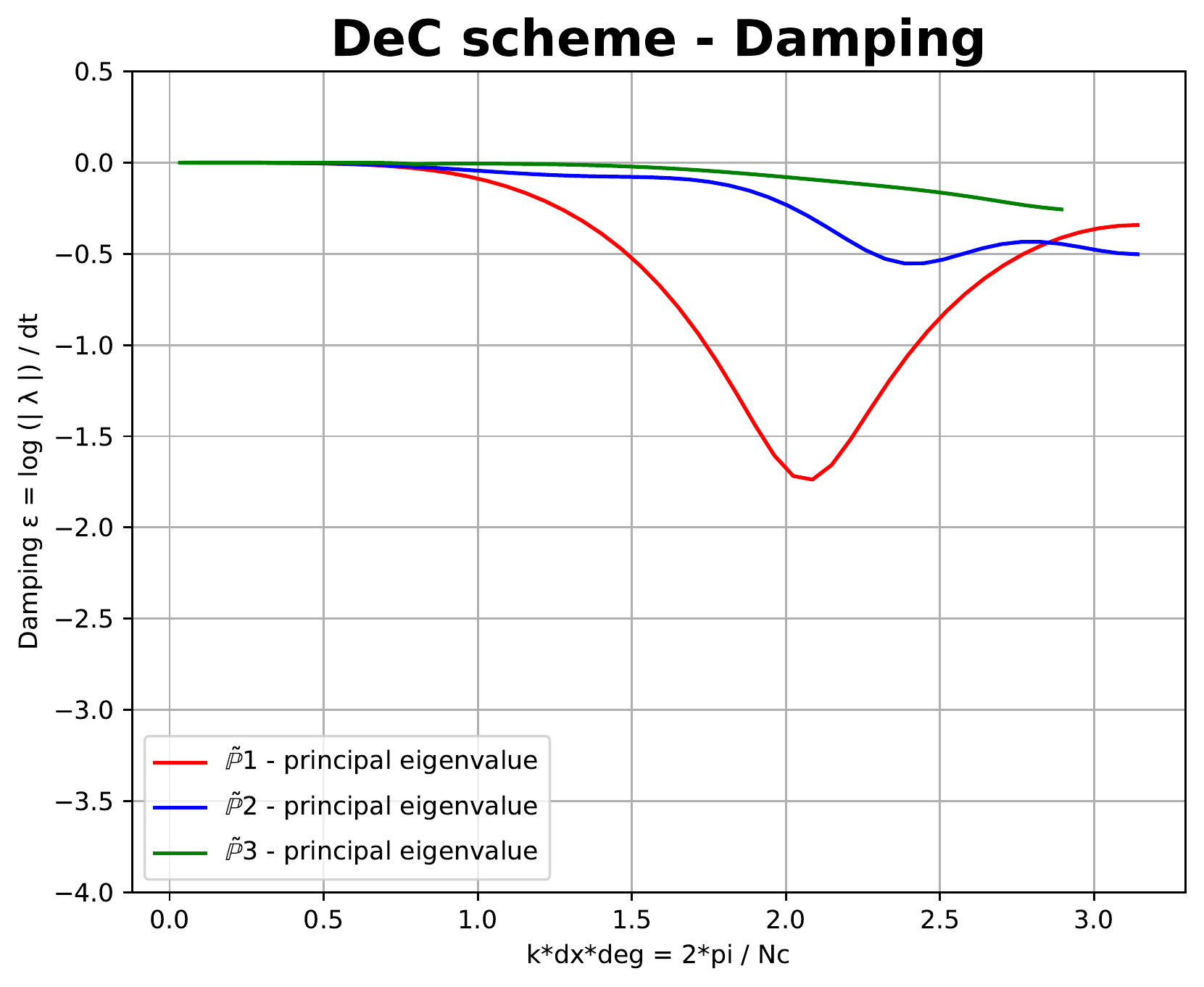}	
\\
	\centering
	\title{Using the CIP stabilization method} \\
	\includegraphics[width=0.24\textwidth]{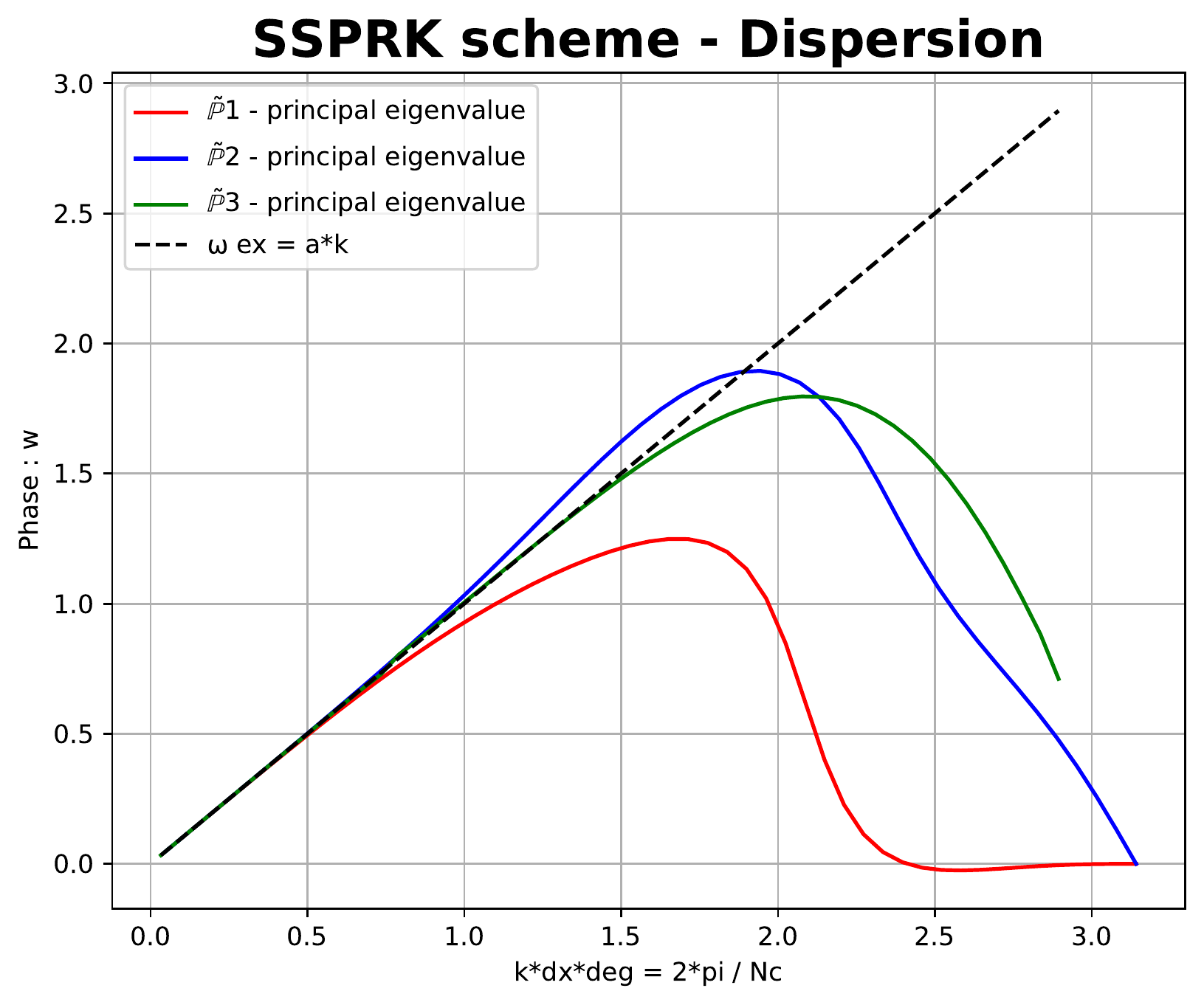}
	\includegraphics[width=0.24\textwidth]{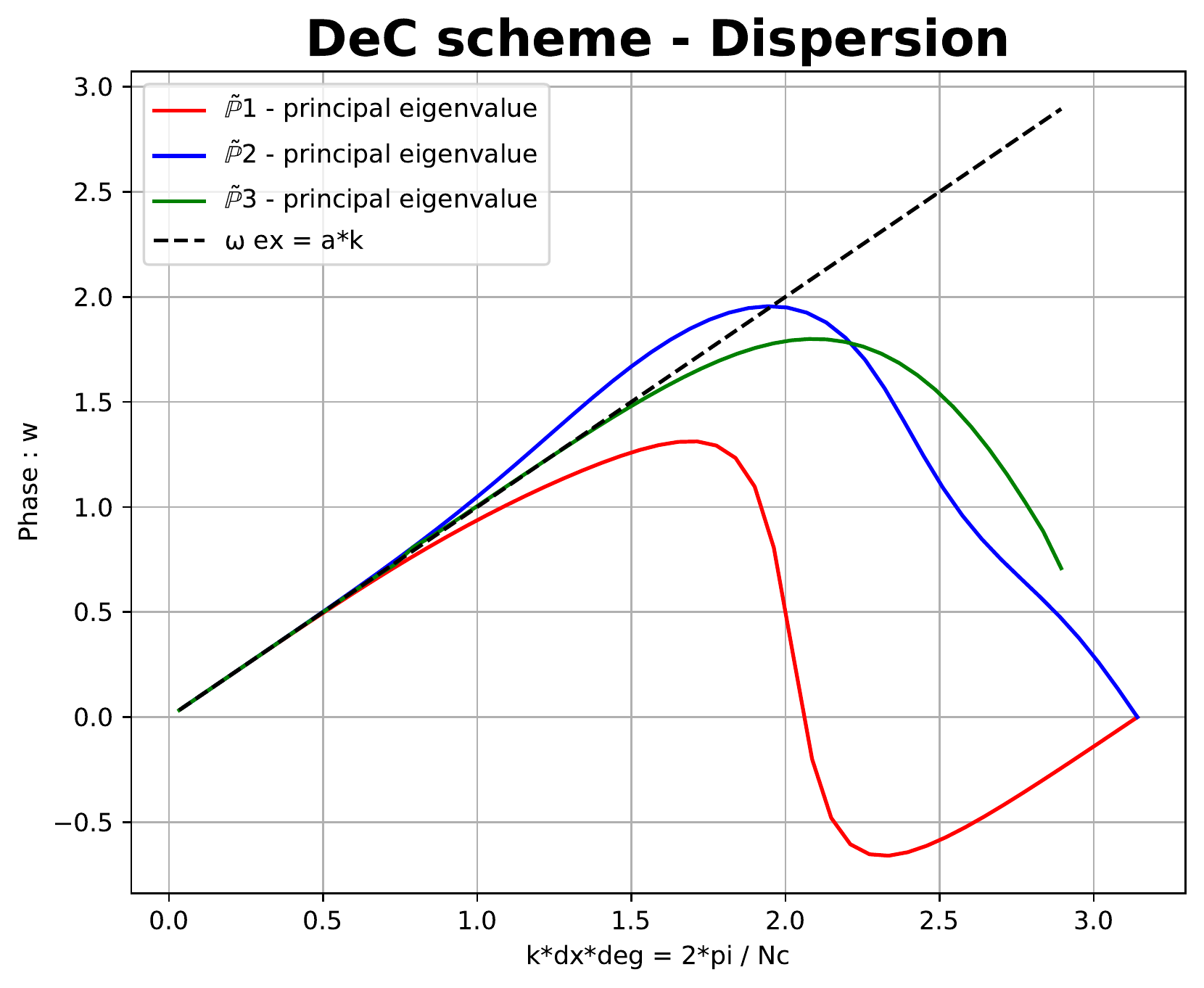}
	\includegraphics[width=0.24\textwidth]{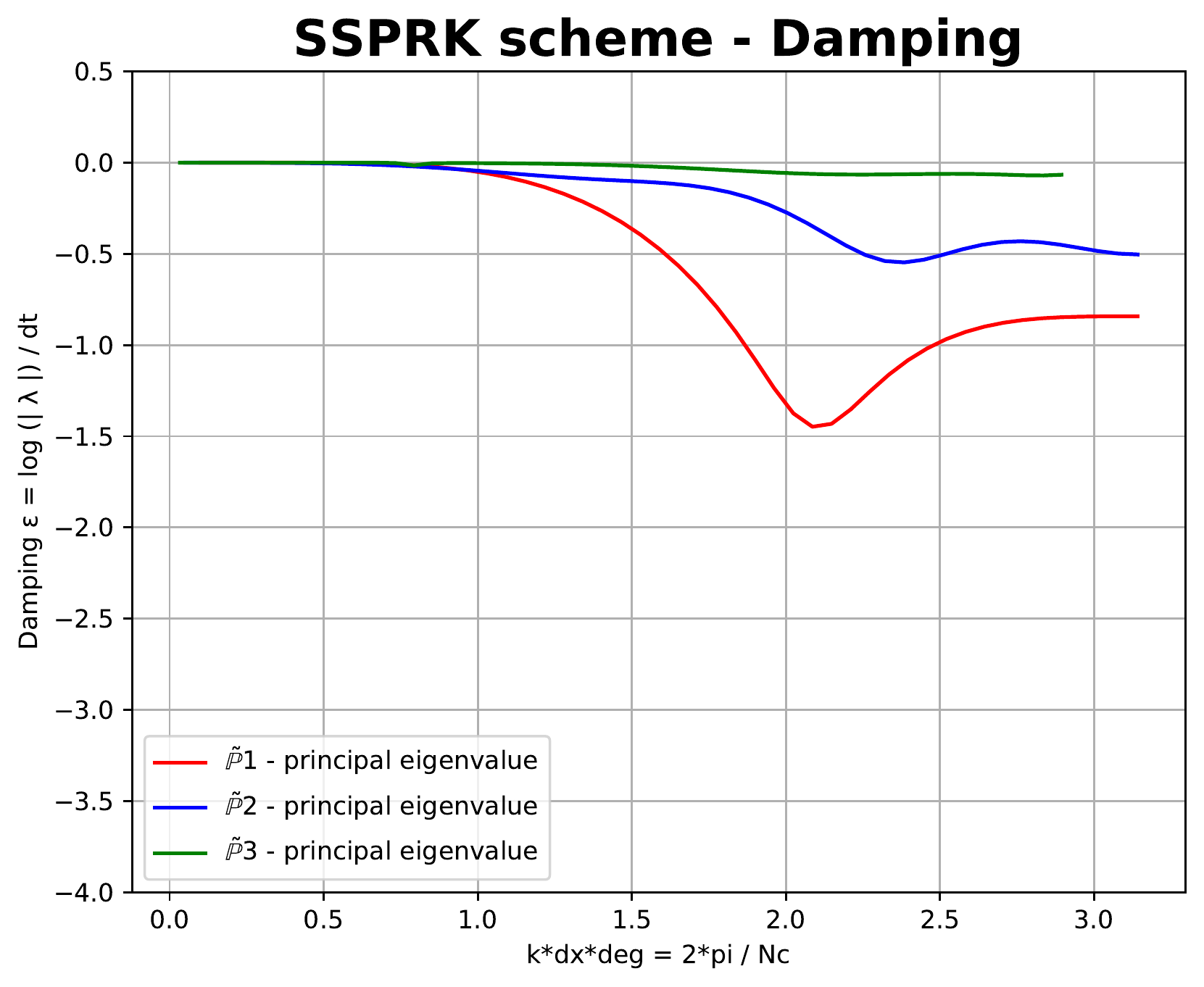}
	\includegraphics[width=0.24\textwidth]{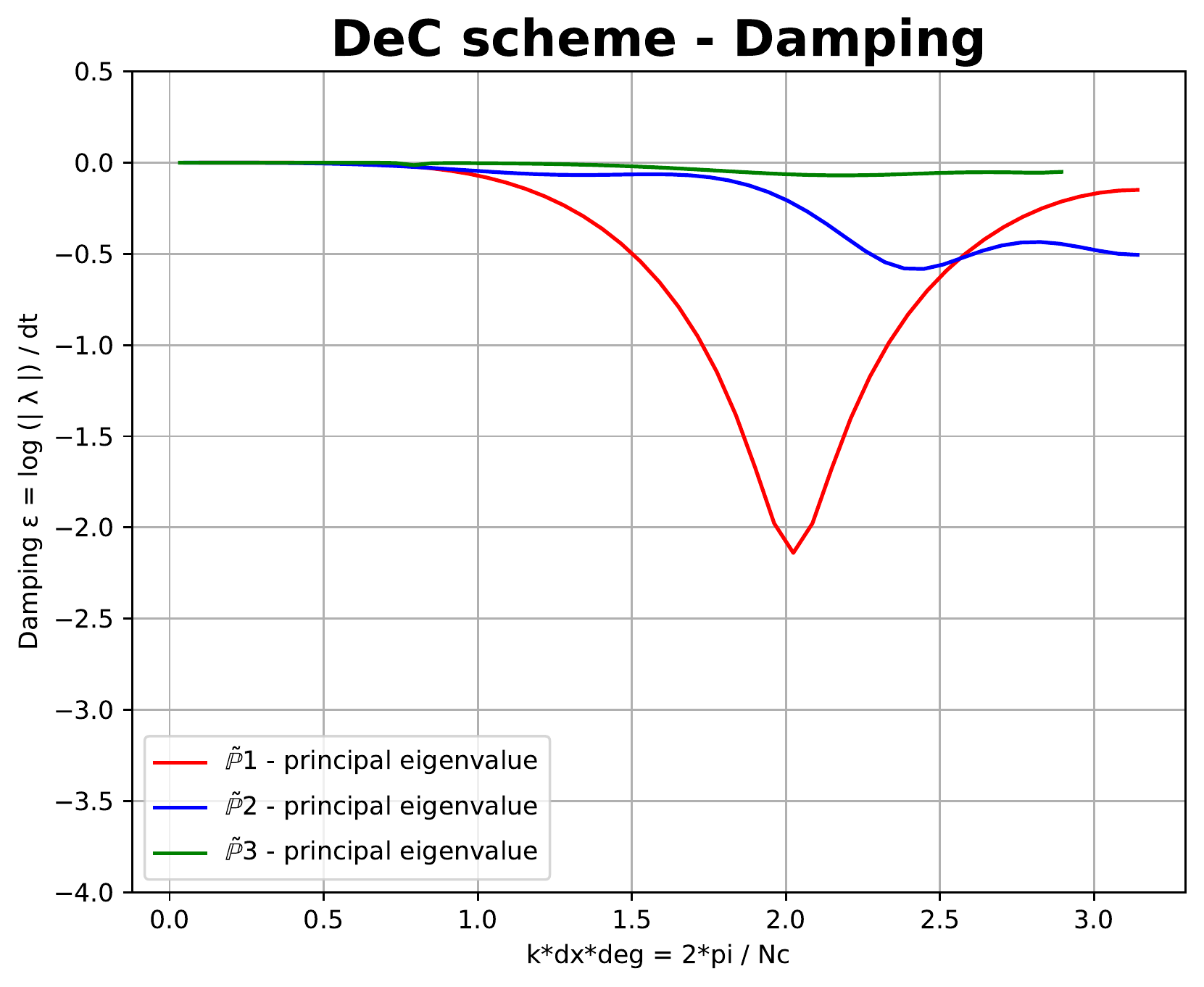}	
	\caption{Dispersion and damping coefficients for \textit{cubature} elements, with DeC and SSPRK methods and all stabilization techniques}\label{fig:dispersionCohen}
\end{figure}

\begin{figure}[h!]
	\centering
	\title{Without any stabilization method} \\
	\includegraphics[width=0.24\textwidth]{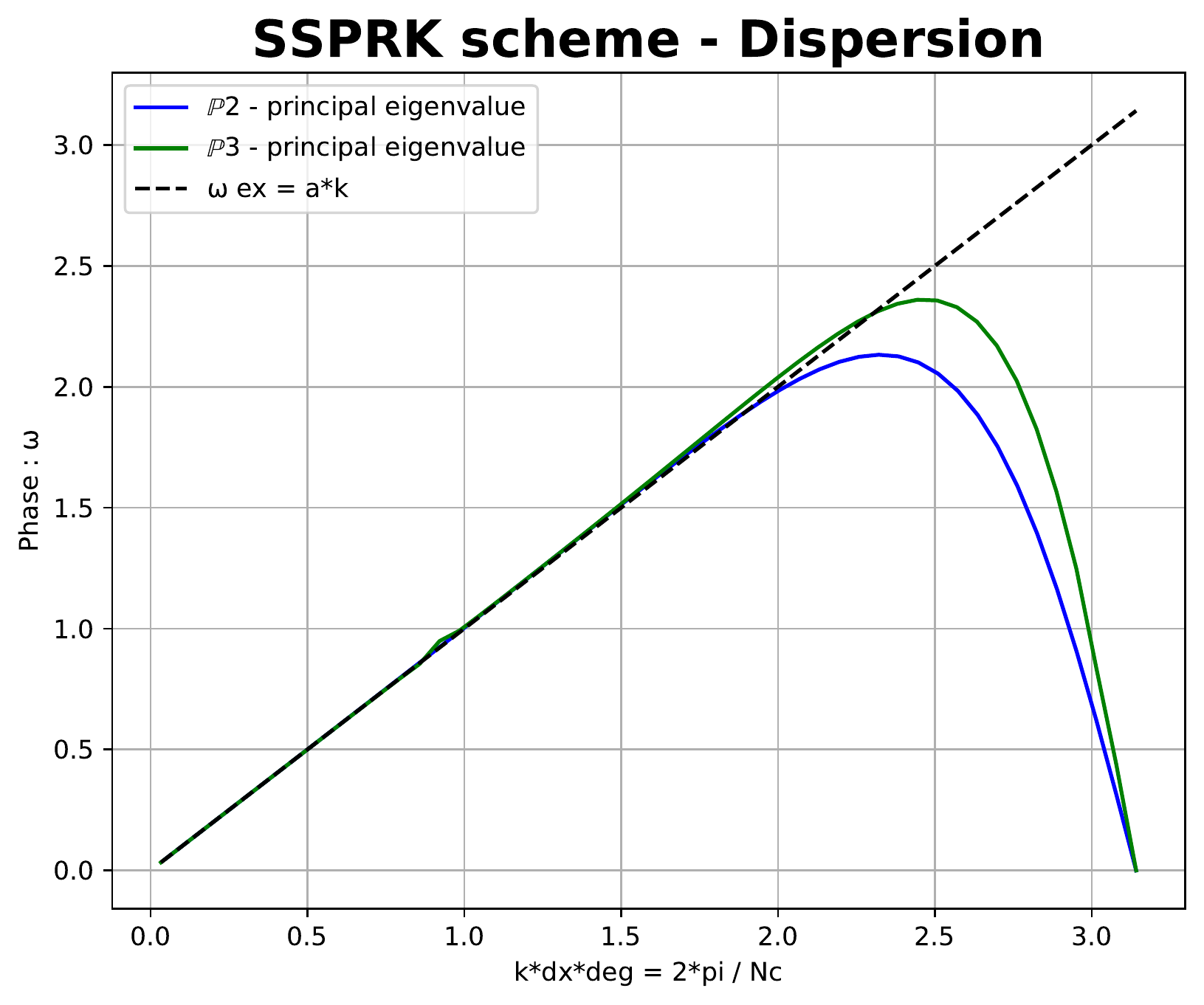}
	\includegraphics[width=0.24\textwidth]{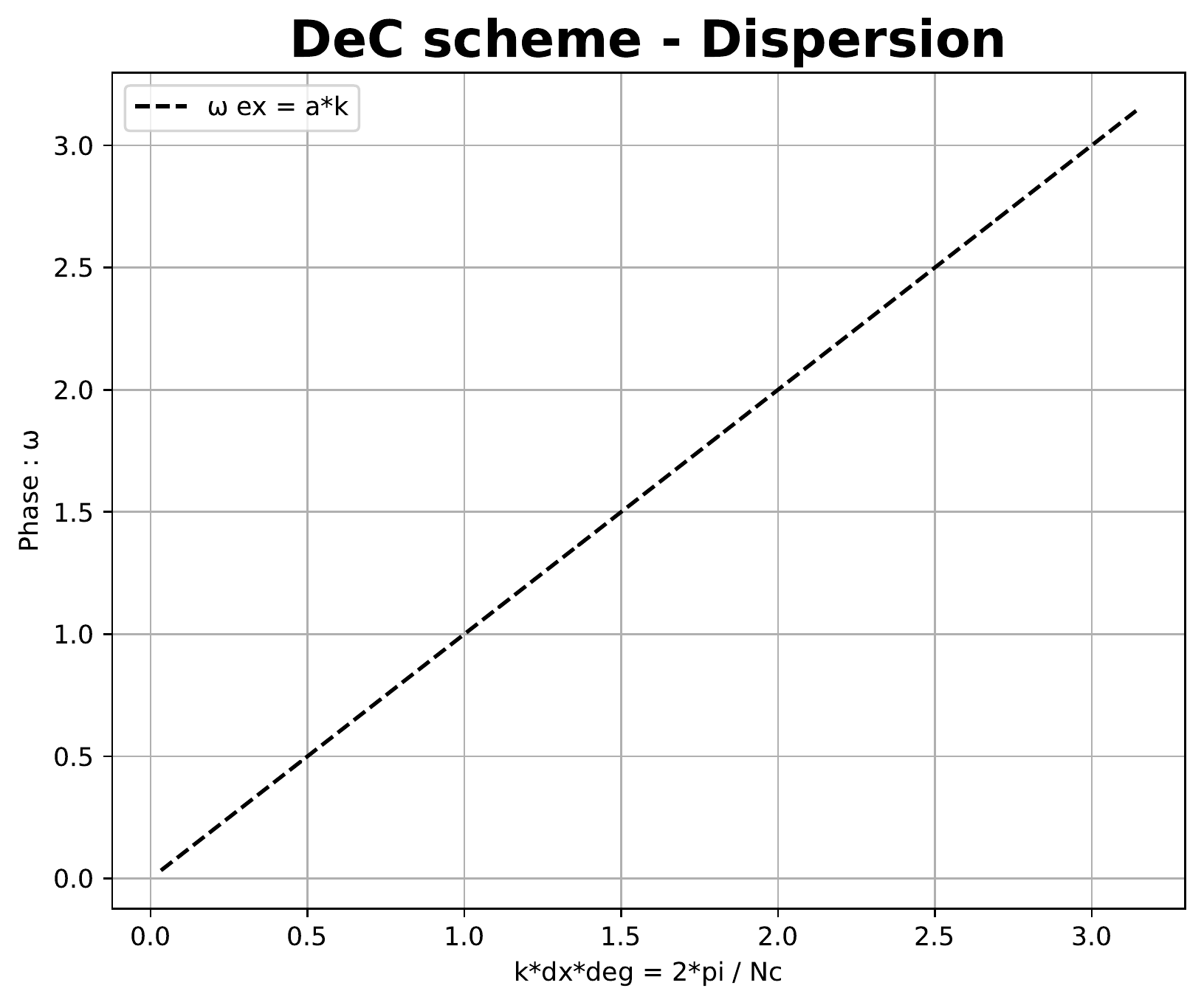}
	\includegraphics[width=0.24\textwidth]{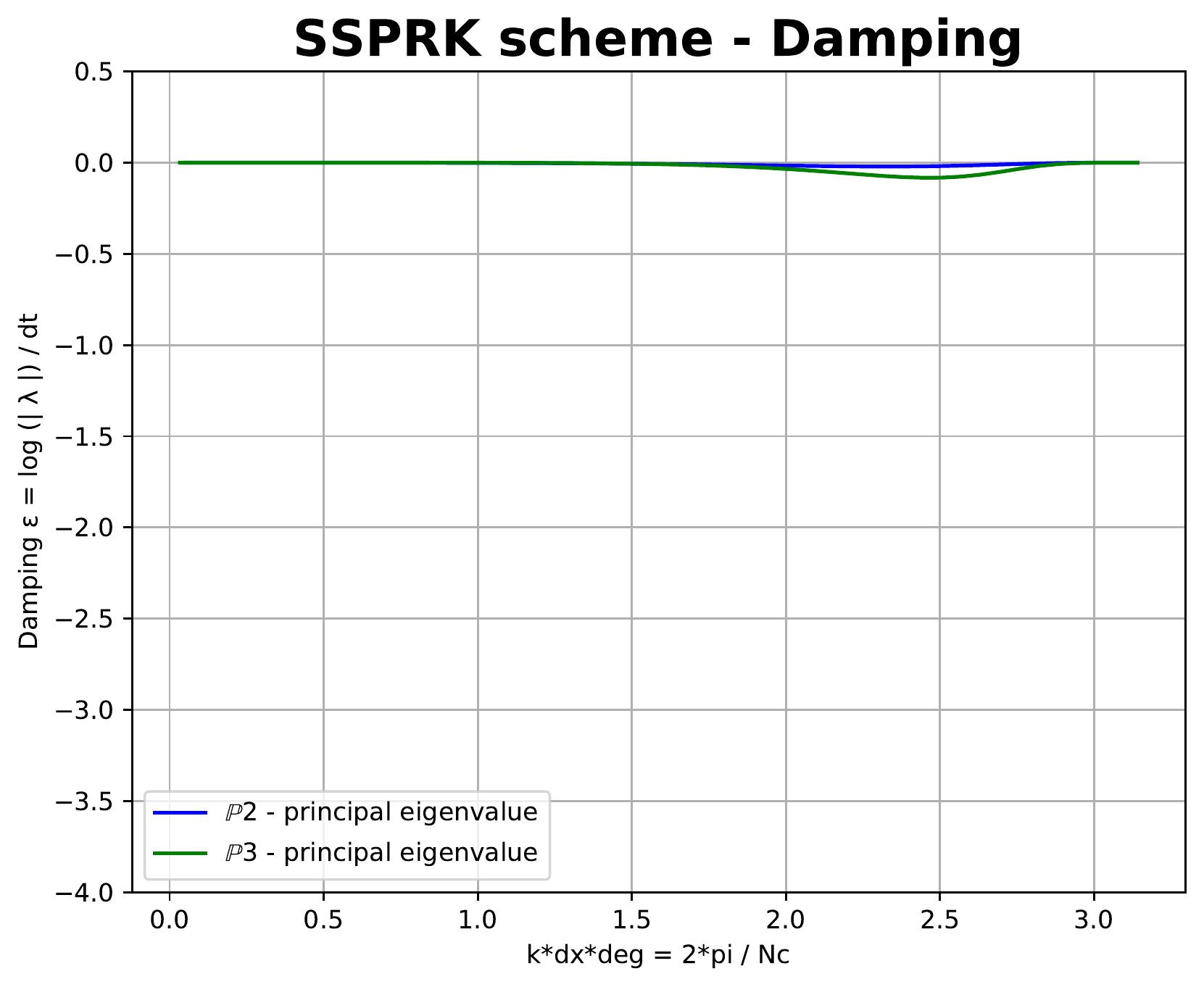}
	\includegraphics[width=0.24\textwidth]{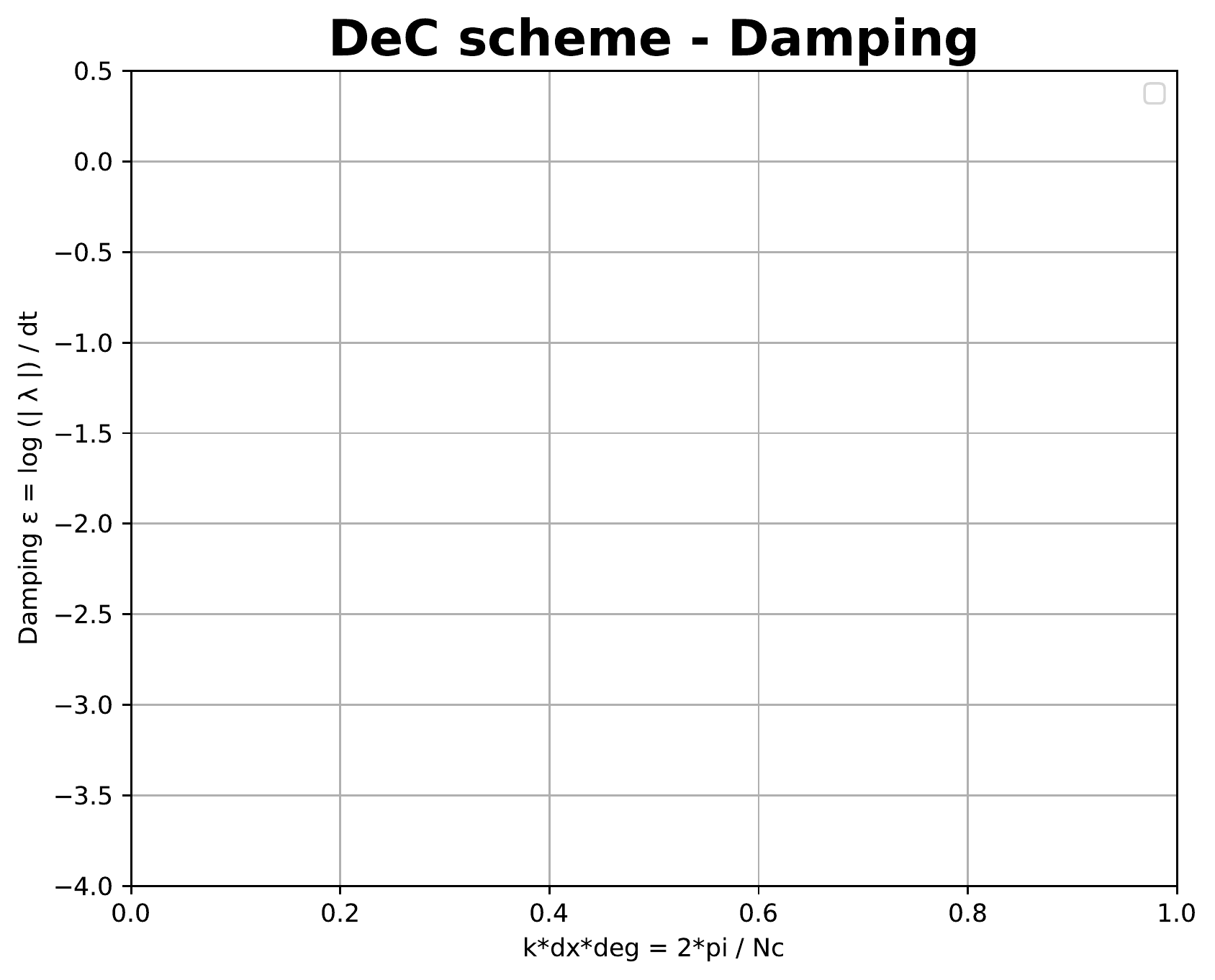}	
\\	
	\title{Using the SUPG stabilization method} \\
	\includegraphics[width=0.24\textwidth]{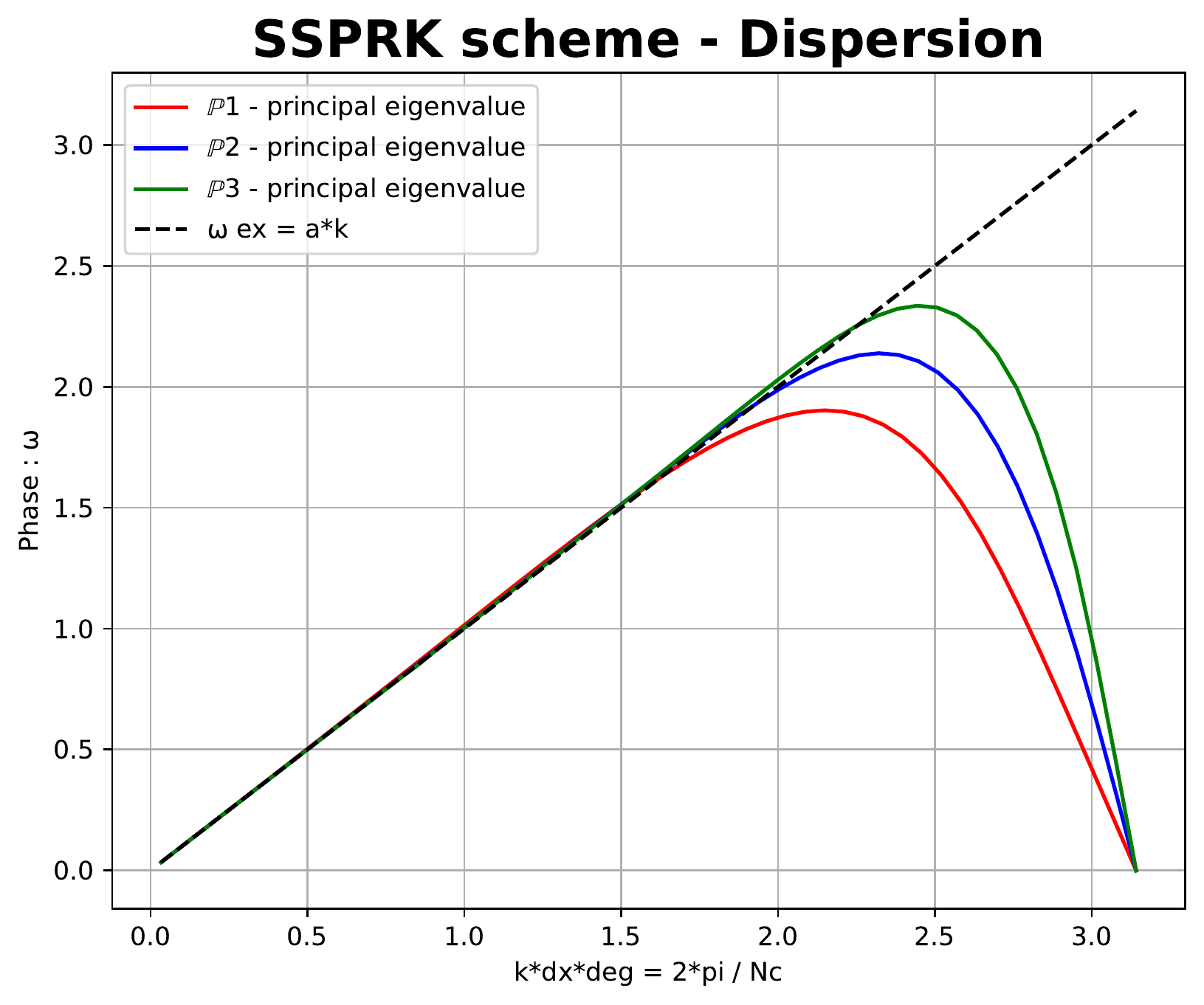}
	\includegraphics[width=0.24\textwidth]{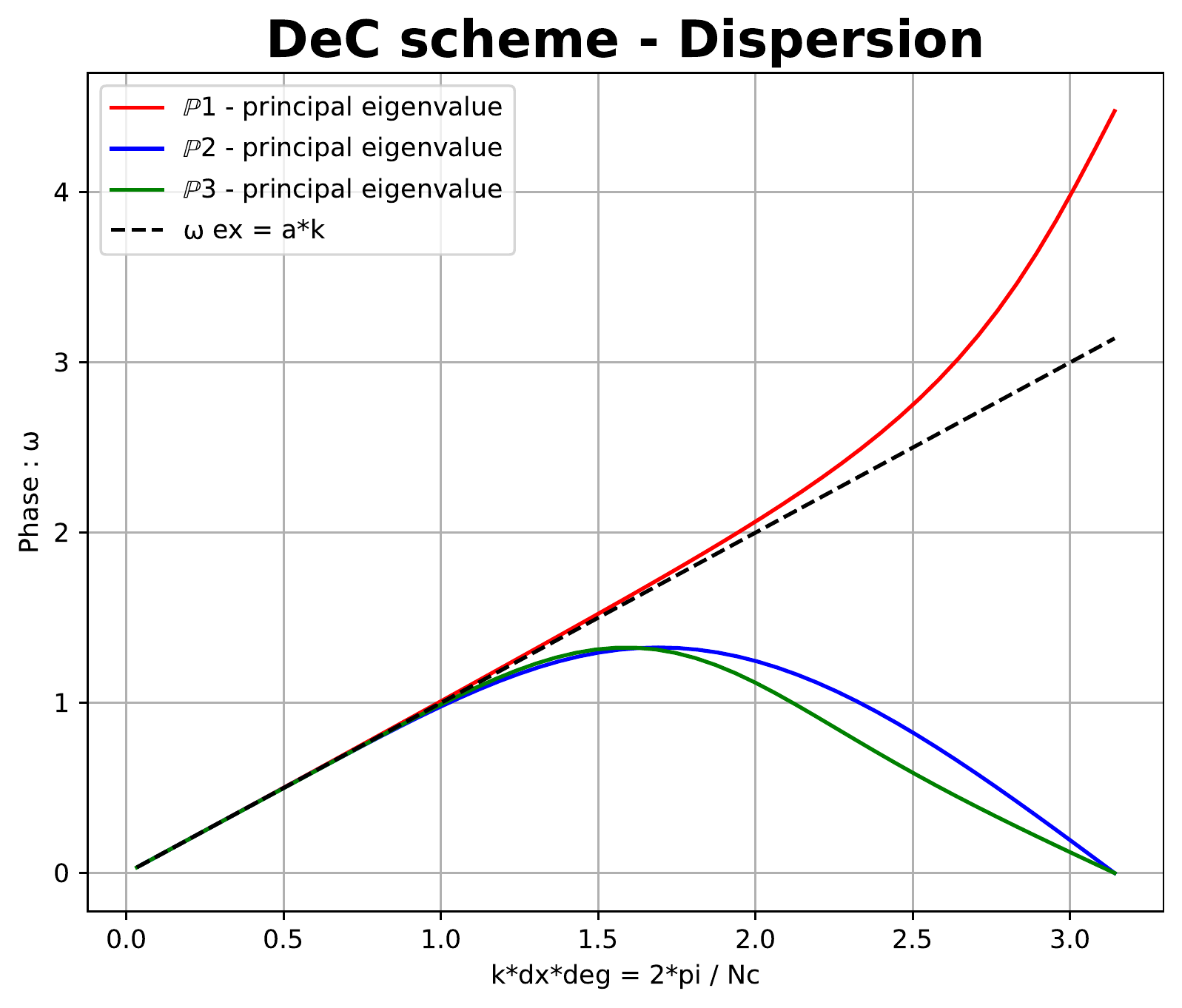}
	\includegraphics[width=0.24\textwidth]{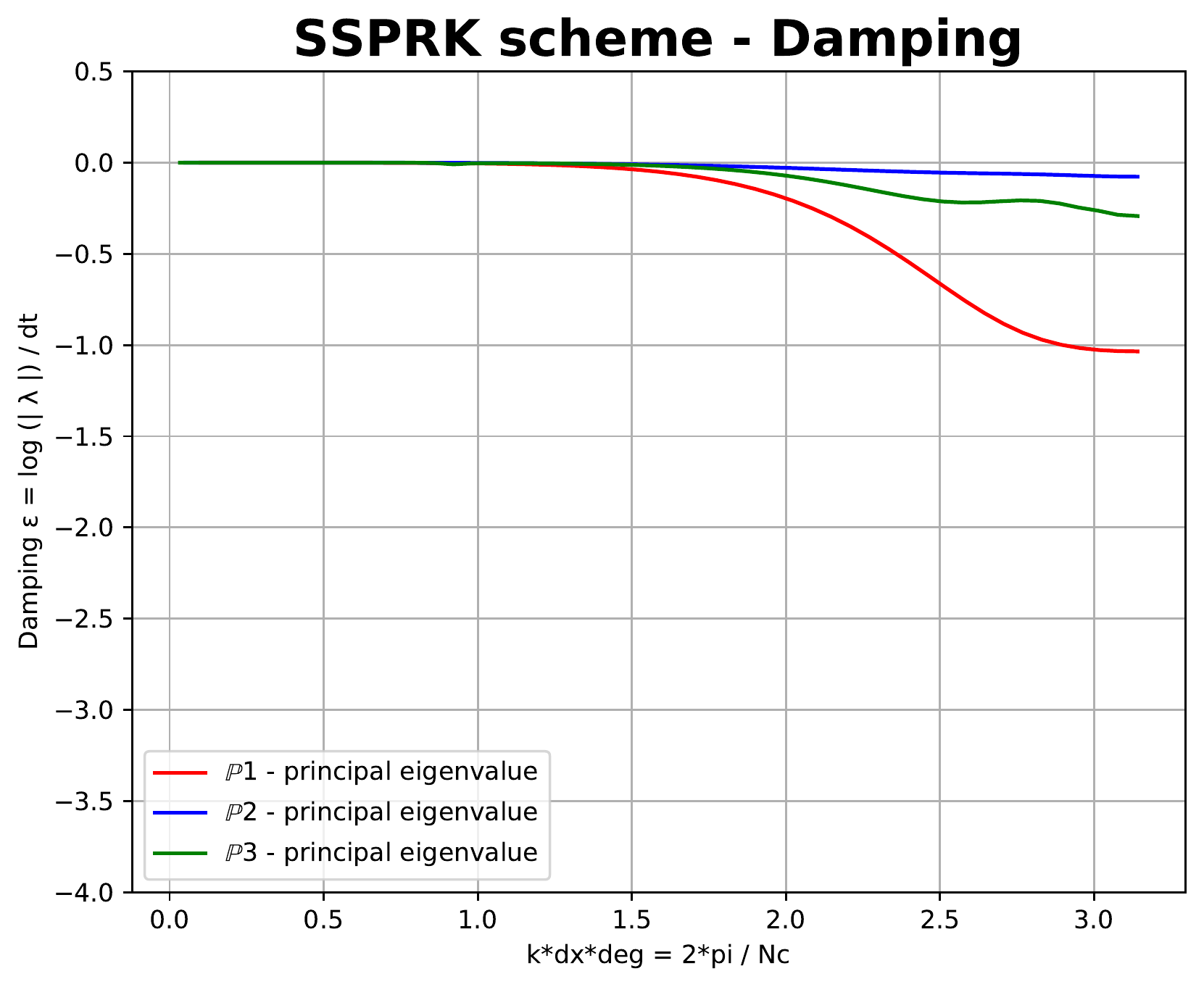}
	\includegraphics[width=0.24\textwidth]{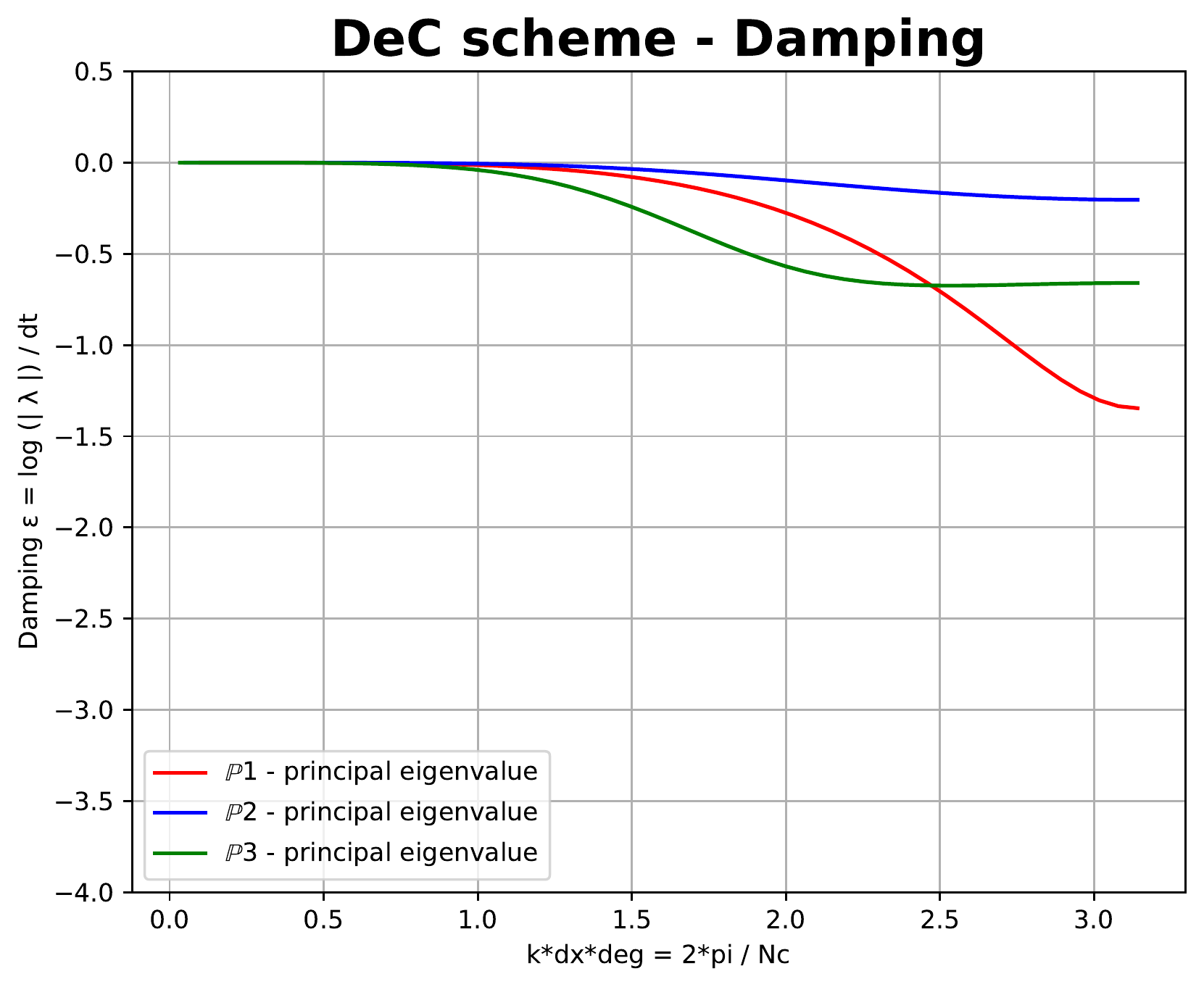}	
\\	
	\title{Using the LPS stabilization method} \\
	\includegraphics[width=0.24\textwidth]{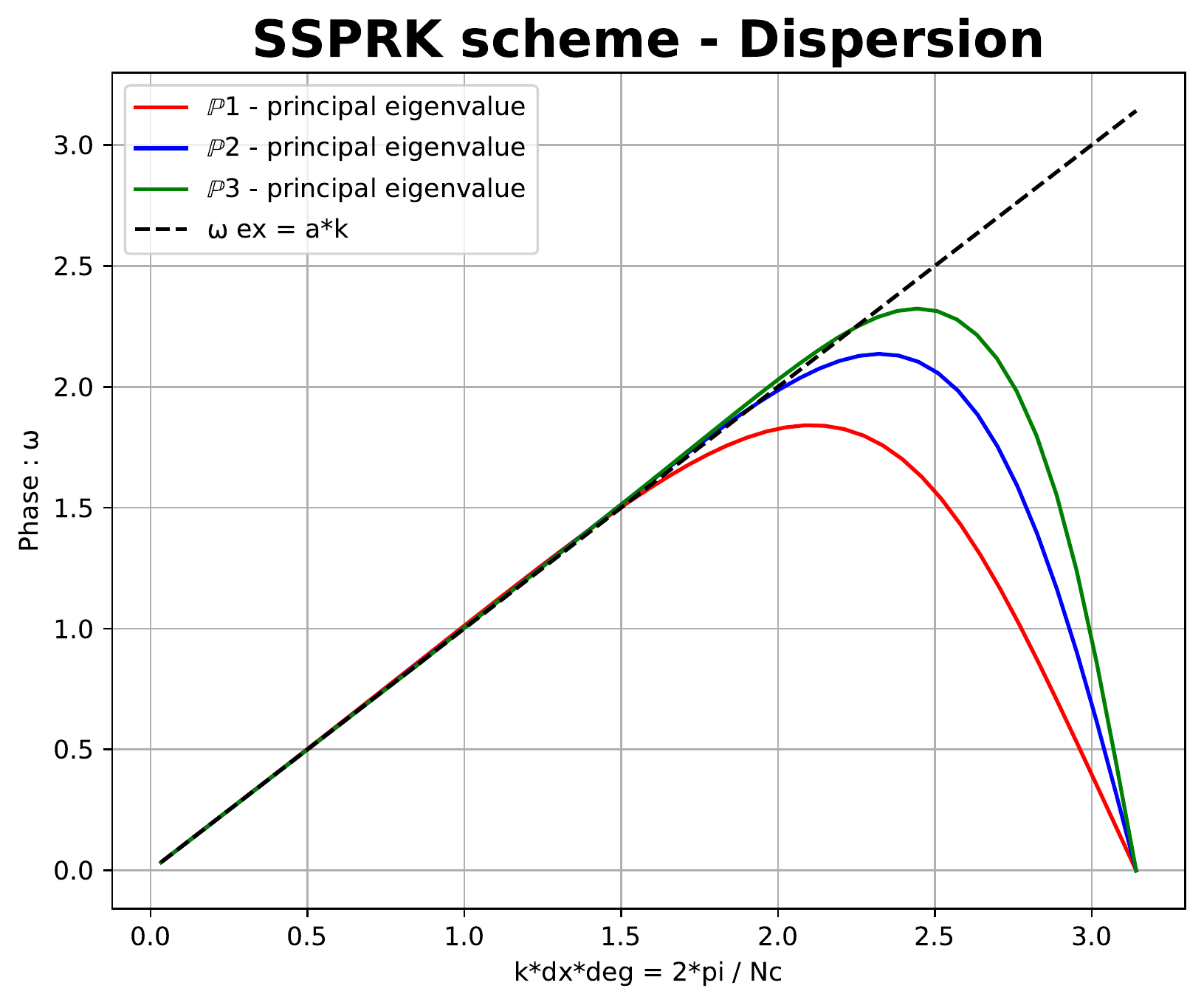}
	\includegraphics[width=0.24\textwidth]{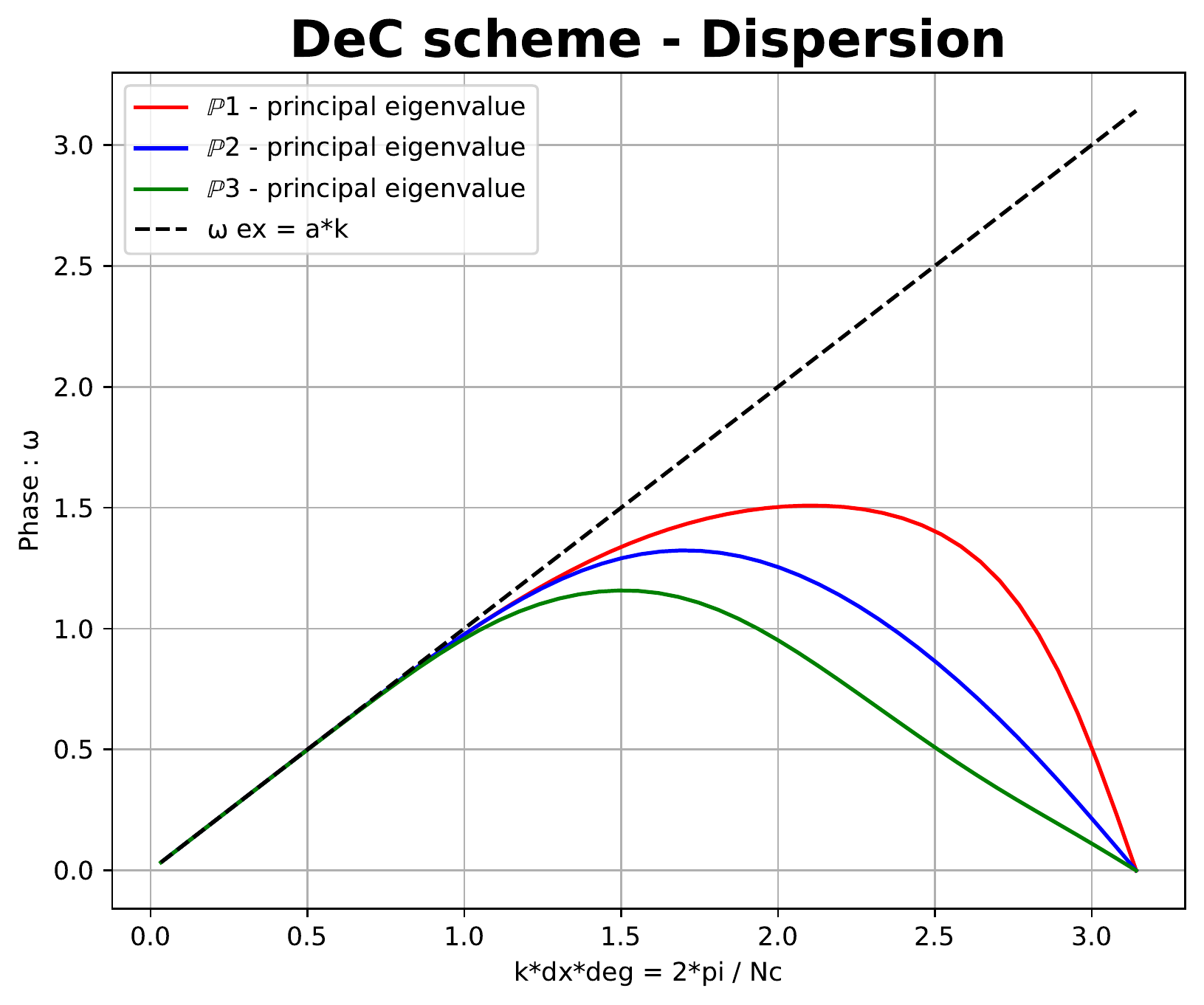}
	\includegraphics[width=0.24\textwidth]{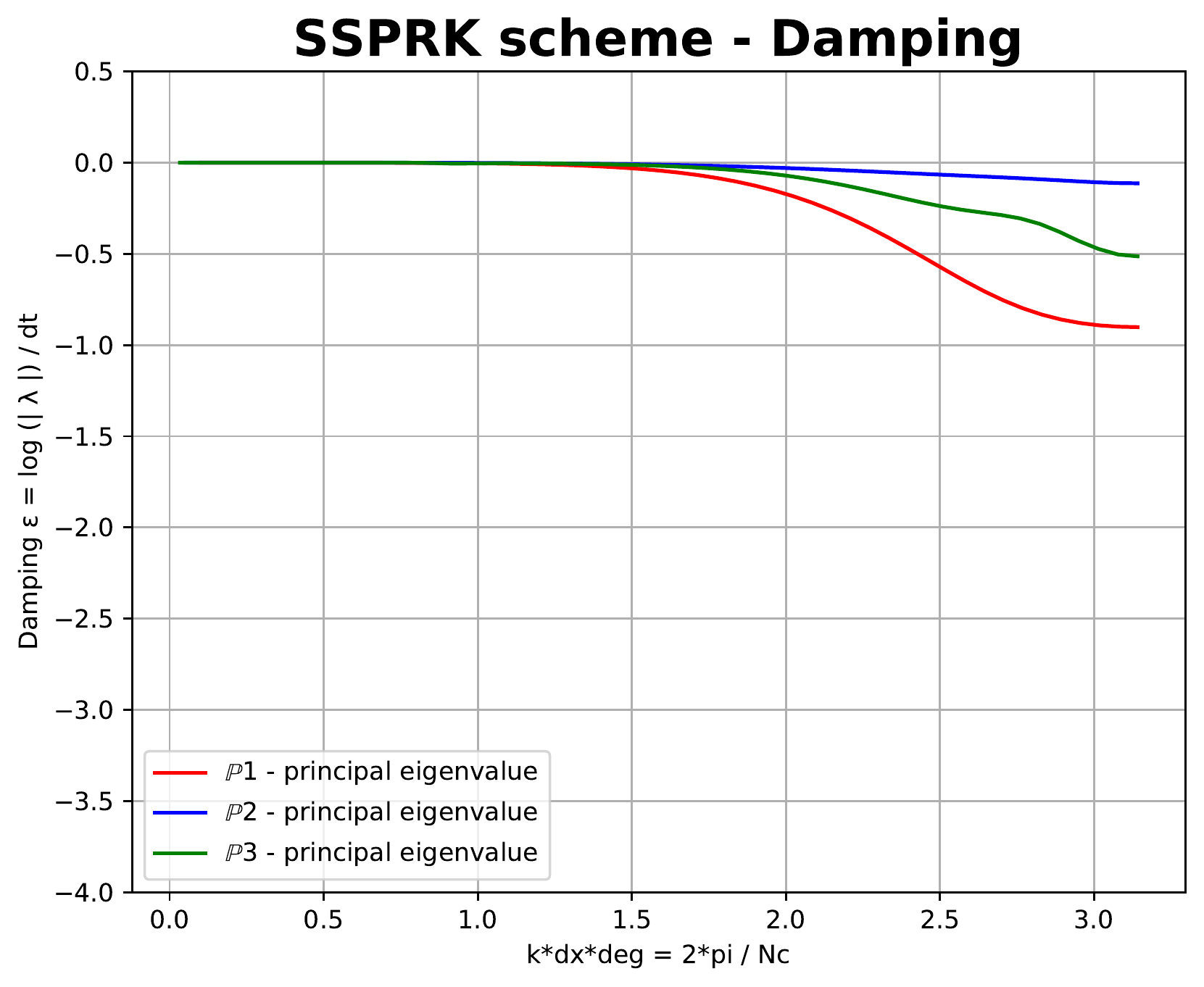}
	\includegraphics[width=0.24\textwidth]{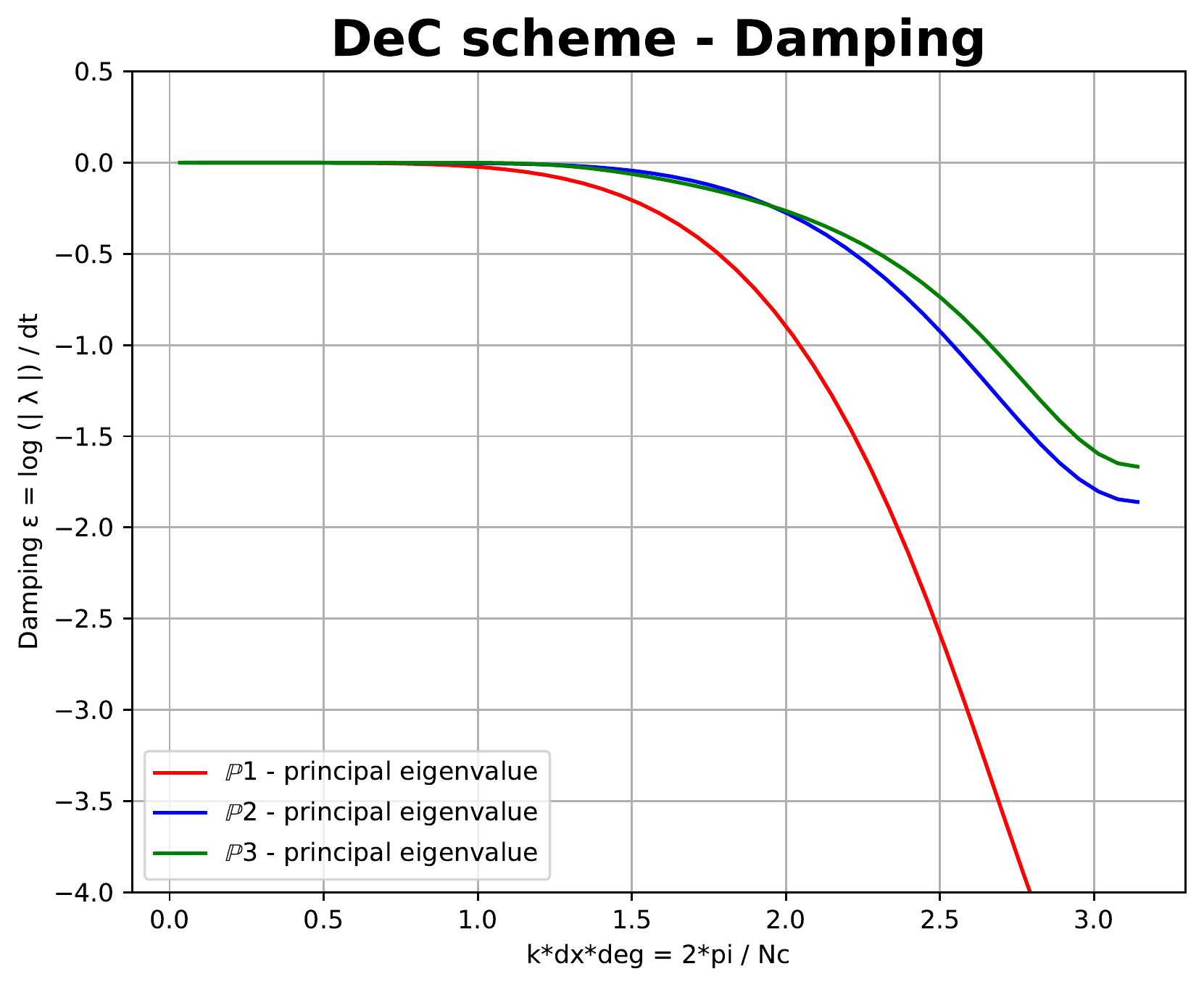}	
	\title{Using the CIP stabilization method} \\
	\includegraphics[width=0.24\textwidth]{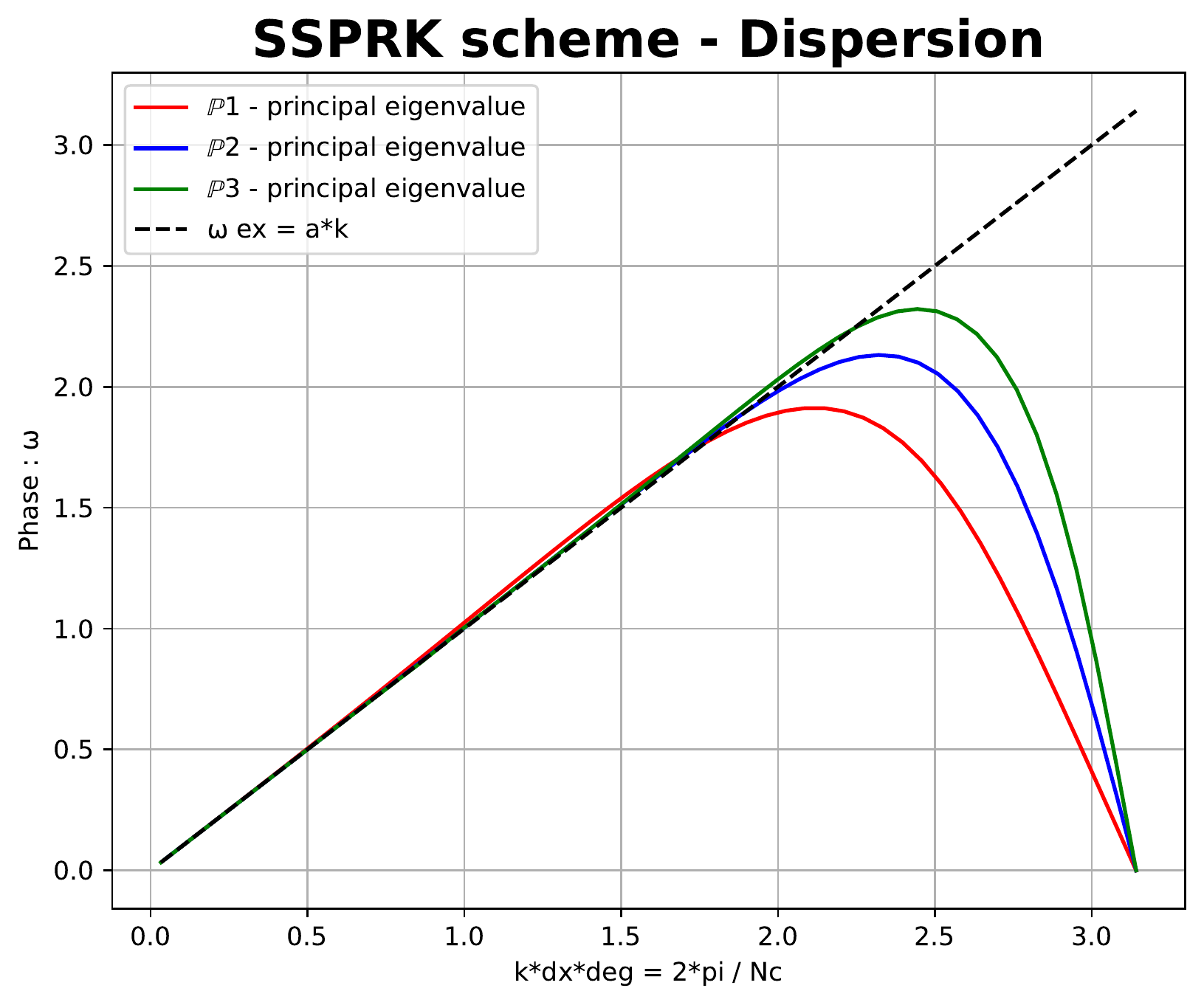}
	\includegraphics[width=0.24\textwidth]{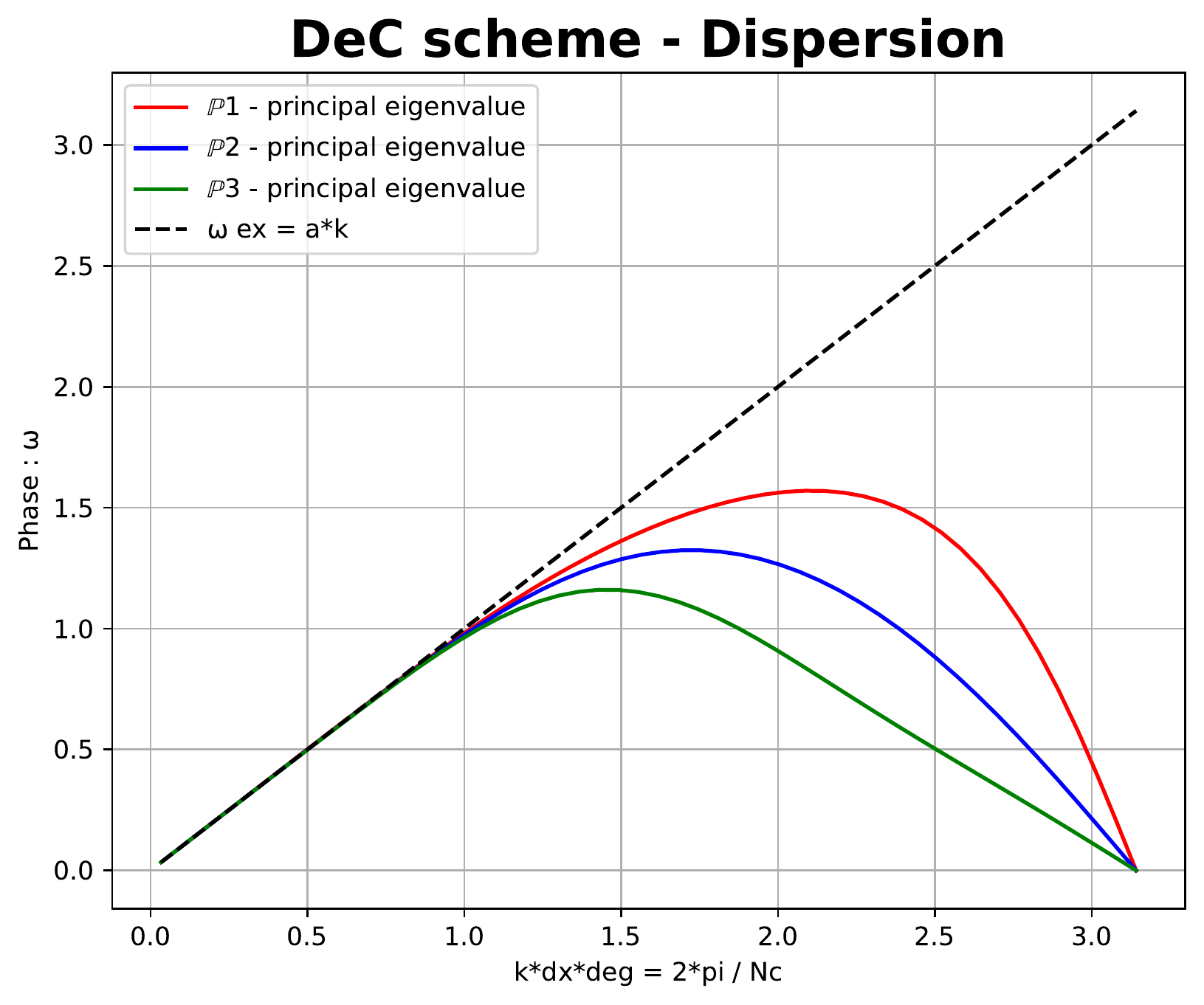}
	\includegraphics[width=0.24\textwidth]{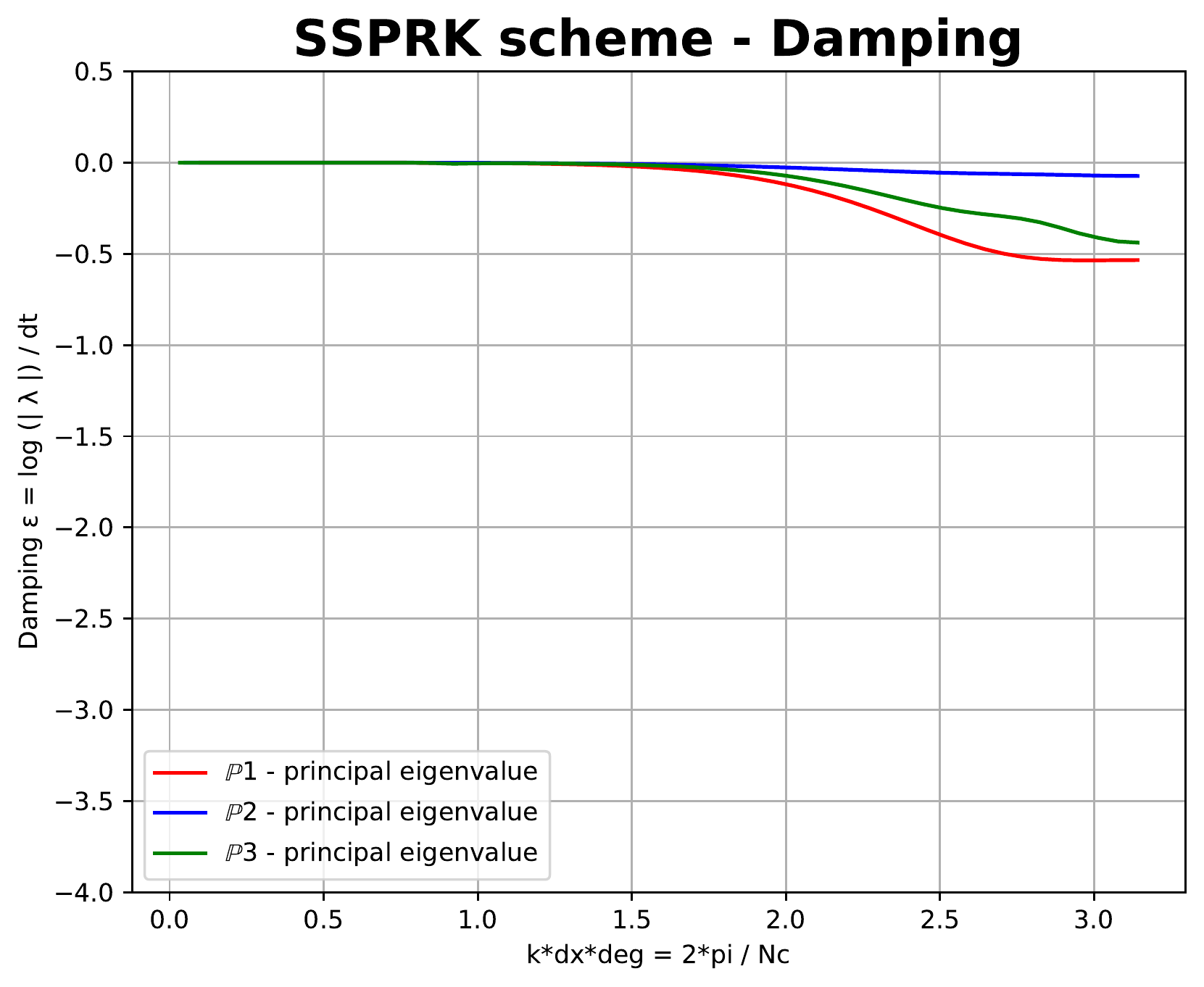}
	\includegraphics[width=0.24\textwidth]{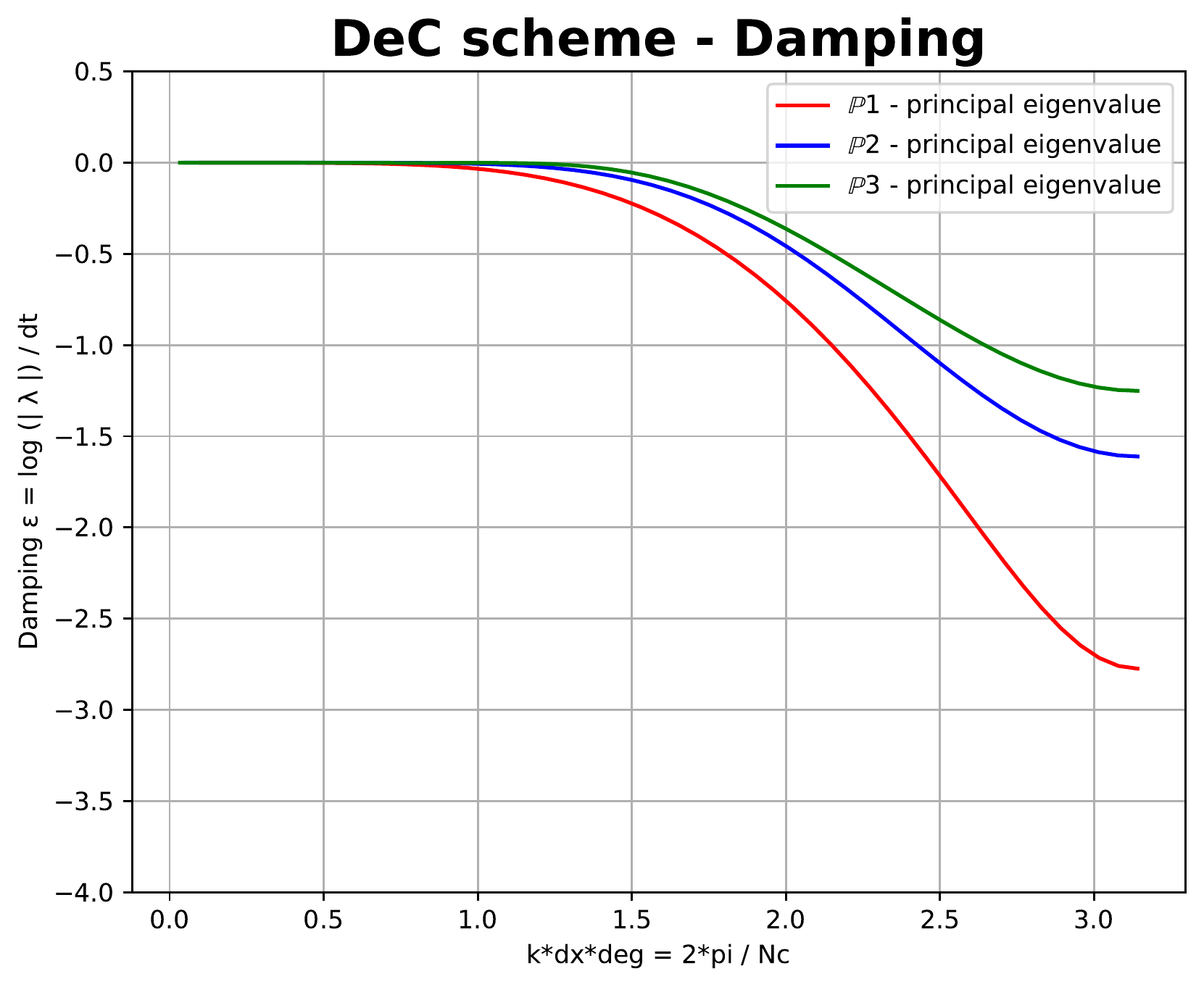}	
	\caption{Dispersion and damping coefficients for \textit{Bernstein} elements, with DeC and SSPRK methods and all stabilization techniques}\label{fig:dispersionBernstein}
\end{figure}

In this appendix we present a summary of the fully discrete Fourier analysis of \cref{sec:fourier_space_time}, comparing different time schemes (SSPRK and DeC), discretizations (\textit{basic}, \textit{cubature}, \textit{Bernstein}), and stabilization methods (LPS, CIP, SUPG). We show the phase $\omega$ and the damping $\epsilon$ coefficients using the \textit{best parameters} obtained by minimizing the relative error of the solution $\eta_u$ for each scheme in \cref{tab:dispersion_cfl-RES}.
When the scheme was unstable we did not plot the mode.
In \cref{fig:dispersionBasic} one finds the phase and the damping for \textit{basic} elements, in \cref{fig:dispersionCohen} for \textit{cubature} elements and in \cref{fig:dispersionBernstein} for \textit{Bernstein} elements.
We remark that for \textit{cubature} elements in \cref{fig:dispersionCohen}, $\Delta x_3$ is scaled differently with respect to the other orders because the point distribution is not equispaced.

In general, we can observe that the phase error increases passing from full matrix SSPRK methods to diagonal one DeC. 
This is noticeable even more for \textit{Bernstein} elements. 
\textit{Cubature }elements, which are not effected by the mass lumping, do not show this behavior, and have a dispersion error which is greater than the other lumped methods, but smaller than the other full mass matrix methods.
This step is also associated to a greater damping in the higher frequencies.


\medskip
\bibliographystyle{siam} 
\bibliography{bibliography}

\end{document}